\documentclass[11pt,a4paper]{amsart}
\usepackage{amssymb,amsthm,amsmath,amsfonts}
\usepackage[english]{babel}
\usepackage{geometry}
\usepackage{xspace}
\usepackage{enumerate}
\usepackage{bbm}
\usepackage{stmaryrd}
\usepackage{hyperref}
\usepackage{mathrsfs}
\allowdisplaybreaks[1]

\newcommand{\liml}{\lim\limits}
\newcommand{\supl}{\sup\limits}

\newcommand{\suml}{\sum\limits}

\newcommand{\pard}[2]{\frac{\partial#1}{\partial#2}}

\newcommand{\FCE}{\Fscr C_0^\infty (E)}
\newcommand{\FCb}{\Fscr C_b^\infty (E)}

\newcommand{\norm}[1]{\lVert#1\rVert}
\newcommand{\lrnorm}[1]{\left\lVert#1\right\rVert}

\newcommand{\abs}[1]{\lvert#1\rvert}
\newcommand{\lrabs}[1]{\left\lvert#1\right\rvert}

\newcommand{\lrbracket}[1]{\left\langle#1\right\rangle}

\newcommand{\lrbr}[1]{\left\langle#1\right\rangle}

\newcommand{\deli}[2]{(#1,{#2})}
\newcommand{\lrdel}[2]{\left(#1,{#2}\right)}

\newcommand{\lrdual}[2]{\sideset{_{#1'}}{_{#1}}{\mathop{\lrbracket{#2}}}}

\newcommand{\dualdel}[3]{\sideset{_{#1'}}{_{#1}}{\mathop{\deli{#2}{#3}}}}

\newcommand{\Ascr} {{\mathscr A}}
\newcommand{\Bscr} {{\mathscr B}}
\newcommand{\Cscr} {{\mathscr C}}
\newcommand{\Dscr} {{\mathscr D}}
\newcommand{\Escr} {{\mathscr E}}
\newcommand{\Fscr} {{\mathscr F}}

\newcommand{\Hscr} {{\mathscr H}}

\newcommand{\Lscr} {{\mathscr L}} 
\newcommand{\Mscr} {{\mathscr M}}
\newcommand{\Nscr} {{\mathscr N}}
\newcommand{\Pscr} {{\mathscr P}}

\newcommand{\Vscr} {{\mathscr V}}

\renewcommand\a{\alpha}
\renewcommand\b{\beta}

\newcommand{\eps}{{\varepsilon}}

\renewcommand{\phi}{\varphi}

\newcommand{\Z}{\mathbbm{Z}}
\newcommand{\R}{\mathbbm{R}}
\newcommand{\N}{\mathbbm{N}}
\newcommand{\Q}{\mathbbm{Q}}
\renewcommand{\P}{\mathbbm{P}}

\newcommand{\M}{\mathbbm{M}}
\newcommand{\C}{\mathbbm{C}}
\newcommand{\K}{\mathbbm{K}}

\newcommand{\esssup}{\operatorname{ess\;sup}}

\renewcommand{\div}{\operatorname{div}}

\newcommand{\supp}{\operatorname{supp}}

\renewcommand{\r}{\restriction}

\newcommand{\Id}{\operatorname{Id}}
\newcommand{\loc}{\textup{loc}}

\newcommand{\disjoint}{
\mathop{
\vphantom\bigcup\!
\mathop{
\smash\bigcup}
\limits^{\genfrac{}{}{0pt}{}{\displaystyle.}{}}
\vphantom\bigcup\!}}

\usepackage[foot]{amsaddr}

\numberwithin{equation}{section}

\newtheorem{theorem}{Theorem}[section]
\newtheorem{assumption}{Assumption}

\theoremstyle{definition}
\newtheorem{defi}[theorem]{Definition}

\theoremstyle{plain}
\newtheorem{lem}[theorem]{Lemma}
\newtheorem{expl}[theorem]{Example}
\newtheorem{prop}[theorem]{Proposition}
\newtheorem{rem}[theorem]{Remark}

\newtheorem{thm}[theorem]{Theorem}

\newtheorem{cor}[theorem]{Corollary}
\newtheorem*{bem*}{Bemerkung}

\newcommand{\ol}{\overline}

\usepackage{array}
\title{Convergence of non-symmetric forms with changing reference measures}
\author[J. M. T\"olle]{Jonas M. T\"olle}
\email{jonas.tolle@aalto.fi}

\address[JMT]{Aalto University\\
Department of Mathematics and Systems Analysis\\
PO Box 11100 (Otakaari 1, Espoo)\\
00076 Aalto\\
Finland}
\date{\today}


\begin{document}

\begin{abstract}
We develop a functional-analytic theory of convergence for non-symmetric generalized Dirichlet forms defined on a sequence of Hilbert spaces $H_n$ converging to a limit space $H$ in the sense of Kuwae--Shioya, with the principal applications concerning $L^2$-spaces over varying reference measures $\mu_n \to \mu$. Extending Mosco's variational convergence for Dirichlet forms to the non-symmetric setting via a functional controlling the off-diagonal values of the forms, we introduce a notion of \emph{generalized convergence of non-symmetric Dirichlet forms} and show it to be equivalent to strong convergence of the associated resolvents, to weak convergence of the co-resolvents, and to strong graph convergence of the infinitesimal generators. A version of the Trotter--Kato theorem in the Kuwae--Shioya framework yields the corresponding equivalence for $C_0$-contraction semigroups. As applications, we give sufficient conditions for the generalized convergence of elliptic and sectorial $a_{i,j}$-forms with changing speed measures, both on $\mathbb{R}^d$ and in infinite dimensions over a Gelfand triple, and deduce weak convergence of the laws of the associated right processes under tightness of the path measures.
\end{abstract}

\maketitle

{\small\tableofcontents}

\section{Introduction}\label{sec:intro}

In this work the author examines convergence problems of non-symmetric forms defined on varying Hilbert spaces.
The aim is to provide necessary and sufficient conditions for convergence of the associated resolvents and semigroups.
The convergence of the associated processes is also considered.

\emph{A motivating example.} Let $d\ge 3$, let $t_n\to t>0$, and let $\mu_n=\sigma_n dx$ and $\mu=\sigma dx$ be the centered Gaussian measures on $\R^d$ with covariance matrices $t_n\Id$ and $t\Id$, respectively. Let $J$ be a non-zero skew-symmetric $d\times d$ matrix and $c\in\R$, and consider the non-symmetric forms
\[\Escr^n(u,v)=\int_{\R^d}\lrbr{\nabla u,\nabla v}\,d\mu_n+c\int_{\R^d}\lrbr{Jx,\nabla u}\,v\,\mu_n(dx)\quad\text{on}\quad H_n:=L^2(\R^d;\mu_n).\]
Their generators are, formally, $L^nu=\Delta u-\lrbr{t_n^{-1}x+cJx,\nabla u}$, and the associated diffusions are Ornstein--Uhlenbeck processes with an additional rotation, given by $dX_s=-(t_n^{-1}\Id+cJ)X_s\,ds+\sqrt2\,dW_s$. Since $\lrbr{Jx,x}=0$, the rotation leaves $\mu_n$ invariant. Hence the rotation does not contribute to the symmetric part $\int\lrbr{\nabla u,\nabla v}\,d\mu_n$ of $\Escr^n$, and it produces the antisymmetric part $\frac12\left(\Escr^n(u,v)-\Escr^n(v,u)\right)=c\int\lrbr{Jx,\nabla u}\,v\,d\mu_n$. As $t_n\to t$, one expects that the resolvents, the semigroups and the laws of these processes converge to the corresponding objects for $t$. Two obstacles appear. First, the forms live on the different Hilbert spaces $L^2(\R^d;\mu_n)$, which do not coincide with each other or with $L^2(\R^d;\mu)$, and in other examples the densities $\sigma_n$ may even degenerate or concentrate. Hence Mosco convergence on a fixed Hilbert space cannot be applied directly. One needs a notion of convergence of vectors $u_n\in H_n$ to $u\in H$, which is provided by the framework of varying Hilbert spaces introduced by Kuwae and Shioya~\cite{KS}. Second, the forms are non-symmetric. The diagonal $\Escr^n(u,u)$ only sees the symmetric part, while the rotation enters only through the antisymmetric part. In the symmetric case, Mosco convergence, a condition on the diagonal only, characterizes the convergence of resolvents. In the non-symmetric case, the off-diagonal values $\Escr^n(w_n,u_n)$ have to be controlled as well, which is the purpose of the conditions (F1) and (F2) below. Forms of this type, with general coefficients and speed measures, are treated rigorously in Theorem~\ref{convergenceaij}. Since the finite dimensional distributions of the associated processes are expressed by their semigroups, the strong convergence of the semigroups along the varying spaces yields the convergence of the finite dimensional distributions. If, in addition, the laws of the processes are tight, which in the example can be checked directly and in general with the help of the Lyons--Zheng decomposition, Theorem~\ref{lawconv} yields the weak convergence of these laws on path space, and this probabilistic consequence is one of the main motivations for the whole framework.

We deal with various concepts of convergence ``along'' a ``sequence of Hilbert spaces'', including
the concept of \textit{generalized convergence} of non-symmetric forms or, more precisely, \textit{generalized Dirichlet forms}.
Our notion of convergence is a generalization of the concept of \textit{Mosco convergence}, a variational convergence of symmetric
quadratic forms introduced by U. Mosco \cite{M69}, \cite{M}. Not only do the considered forms fit the general framework of
generalized Dirichlet forms introduced by W. Stannat~\cite{S2}, \cite{S1}\footnote{For the stochastic calculus associated with generalized Dirichlet forms, see~\cite{Tr}.}, they are also assumed to be defined on different
Hilbert spaces. This idea is due to K. Kuwae and T. Shioya, who developed this framework in~\cite{KS} as a consequence of research on convergence of
metric measure spaces, see also the extension by the same authors~\cite{KS2}. The framework has been extended to dual pairs of Banach spaces by the author of this work in his PhD thesis~\cite{Toe2}. In particular, in~\cite{Toe2}, it was verified that Kuwae--Shioya convergence generates a sequential topology in the sense of convergence classes, see Remark~\ref{seqtop}. A. Kolesnikov proved in~\cite[Proposition 7.2]{K1} that Kuwae--Shioya convergence is metrizable for a countable sequence of separable infinite dimensional Hilbert spaces.

Former research on convergence of forms considers, with a few exceptions, forms on one fixed Hilbert space. To understand the conceptual difference, we
explain what we mean by ``convergence of Hilbert spaces'' in applications. An abstract functional analytic introduction to this
framework and its consequences can be found in Subsection~\ref{ss:spectral}, where many proofs are taken either from~\cite{KS} or from~\cite{K1} and can be
found in Appendix A in order to make this work as self-contained as possible.

Let $E$ be an infinite dimensional real locally convex topological vector space. Let $\mu_n$, $n\in\N$, $\mu$ be fully supported
Borel probability measures such that $\mu_n\to\mu$ weakly. Define $H_n:=L^2(E;\mu_n)$, $n\in\N$, $H:=L^2(E;\mu)$ and let
$C:=\FCb$ be the space of cylindrical test functions. Assume that $E$ is a Souslin space, so that $C$ is dense\footnote{The argument for this can e.g. be found in~\cite[Section II.3.a), (3.2)]{MR}.} in $H$. Now $\{H_n\}$ converges
to $H$ in the following sense: there exists a sequence $\{\Phi_n:C\subset H\to H_n\}$ of injective linear operators with dense linear domain $C$ such that
\[\lim_n\norm{\Phi_n(u)}_{H_n}=\norm{u}_H,\quad\forall u\in C.\]
Clearly, if we take each $\Phi_n$ as the identity operator on $C$,
$H_n$ converges to $H$ in this sense by weak convergence of measures.

If we let $E:=\R^d$ and $\mu_n$, $n\in\N$, $\mu$ be fully supported Radon measures such that
$\mu_n\to\mu$ vaguely, we clearly obtain by setting $C:=C_0^\infty(\R^d)$, $H_n:=L^2(\R^d;\mu_n)$, $n\in\N$, $H:=L^2(\R^d;\mu)$ that
$H_n$ converges to $H$ in this sense. Let us denote the disjoint union of Hilbert spaces by $\Hscr:=\disjoint_n H_n\dot\cup H$.

Now one can define the concept of strong and weak convergence of vectors $\{u_n\}$ to $u$ with $u_n\in H_n$, $n\in\N$, $u\in H$. This is done via
approximation in $C$, since the spaces $H_n$ and $H$ are not subspaces of a common normed space.
\begin{defi}[Strong and weak convergence]
A sequence of vectors $\{u_n\}$, $u_n\in H_n$, $n\in\N$  is said to \textit{strongly converge} to $u\in H$ if
there exists an approximating sequence $\{\phi_m\}\subset C$ with $\phi_m\to u$ in $H$ and
\[\lim_m\varlimsup_n\norm{\Phi_n(\phi_m)-u_n}_{H_n}=0.\]
A sequence of vectors $\{u_n\}$, $u_n\in H_n$, $n\in\N$  is said to \textit{weakly converge} to $u\in H$ if
\[\lim_n (u_n,v_n)_{H_n}=(u,v)_H\]
for every strongly convergent sequence $v_n\to v$.
\end{defi}
Within this framework, we also define weak and strong convergence of bounded linear operators from $H_n$ to $H_n$, $n\in\N$ and the analog
to Mosco convergence and to $\Gamma$-convergence in this framework. This has mainly been done by K. Kuwae and T. Shioya in~\cite{KS}. The analysis of
Dirichlet forms by convergence problems goes back to~\cite{DGS}, \cite{DGF}, where $\Gamma$-convergence introduced by E. De Giorgi was used
to obtain asymptotic properties of Dirichlet forms. A complete treatment of $\Gamma$-convergence on arbitrary topological spaces can be found in~\cite{DaMa}.
In~\cite{M69}, U. Mosco introduced a variant of this convergence for convex sets and convex functionals. It amounts to $\Gamma$-convergence with respect to both the strong and the weak topology. In~\cite{M}, he studied it for closed quadratic forms. His main result
was to identify Mosco convergence as a necessary and sufficient condition on a sequence of symmetric closed forms $\{\Escr^n\}$
such that the associated $C_0$-contraction resolvents
$\{G_\a^n\}_{\a>0}$ converge in the strong operator topology.
For an introduction to symmetric Dirichlet forms we refer to~\cite{FOT} or \cite{MR}. An excellent exposition of Mosco convergence of convex sets and convex functionals can be found in H. Attouch's monograph~\cite{A}.
\begin{defi}[Mosco convergence on one fixed Hilbert space]
We say that a sequence of symmetric closed forms $\{\Escr^n\}$ on some common separable Hilbert space $H_0$ \textit{Mosco converges} to a symmetric closed form $\Escr$
on the same Hilbert space $H_0$ if
the following two conditions are fulfilled, where $\Dscr(\Escr)$ denotes the domain of $\Escr$, and where we extend every $\Escr^n$, $n\in\N$, and $\Escr$ to $H_0$ by setting $\Escr(u):=+\infty$ if $u\in H_0\setminus\Dscr(\Escr)$, and similarly for each $\Escr^n$:
\begin{enumerate}
\item[(M1)] For every weakly convergent sequence of vectors $u_n\to u$ we have
\[\Escr(u)\le\varliminf_n\Escr^n(u_n).\]
\item[(M2)] For every $v\in\Dscr(\Escr)$ there exists a strongly convergent sequence $v_n\to v$ such that
\[\Escr(v)=\lim_n\Escr^n(v_n).\]
\end{enumerate}
\end{defi}
The convergence of the associated $C_0$-contraction semigroups $\{T_t^n\}_{t\ge 0}$ in the strong operator topology is
obtained either by an analysis of the associated spectral measures or by the more general theorem of T. Kato~\cite[Theorem IX.2.16]{Ka}.

Combining the concepts of~\cite{M} and \cite{KS}, analogous results for Mosco convergence can be obtained, see Subsection~\ref{ss:symmetric}.
There has been some research on this approach, since, as can be seen in the examples above, a change of reference measures or speed measures
can be considered here. We refer to the papers of A.V. Kolesnikov~\cite{K1}, \cite{K2} and \cite{K3}, where other useful extensions
of the Kuwae--Shioya framework are also proved. These results are reviewed in Section~\ref{sec:mosco}, since Mosco
convergence of the symmetric parts of non-symmetric forms helps to establish generalized convergence in the sense of Definition~\ref{genint} below. With some important exceptions, the proofs are omitted.

Convergence problems for Dirichlet forms and the associated operators have been examined in
\cite{AHKS}, \cite{AKS}, \cite{CS}, \cite{Can}, \cite{H}, \cite{Kas}, \cite{K1}, \cite{K2}, \cite{K3}, \cite{KS}, \cite{KU1}, \cite{KU2}, \cite{LZa}, \cite{LZb},
\cite{LZ}, \cite{Mat}, \cite{Me}, \cite{M}, \cite{OTT}, \cite{Po}, \cite{PoZ}, \cite{RZ2}, \cite{Str}, \cite{Sun} and \cite{Uem}.

We note that S. Mataloni in~\cite{Mat} and P. Mertens in~\cite{Me}
first considered an abstract convergence of sectorial forms, though non-symmetric cases
can also be found in~\cite{RZ} and other papers, and that M. Hino in~\cite{H} first stated abstract conditions on generalized forms in the sense of~\cite{S1}. In all three
papers the strong convergence of the associated resolvents was proved. More precisely, M. Hino proved in~\cite[Theorem 3.1]{H} that, for generalized forms on one fixed Hilbert space, conditions of the type (F1), (F2) below are equivalent to the strong convergence of the resolvents, whereas~\cite{Mat} and \cite{Me} give conditions for sectorial forms on a fixed Hilbert space. The present paper transfers Hino's characterization to varying Hilbert spaces. The main additional tools are the isometries $\Psi_n$ of Theorem~\ref{Kolmetric} and the diagonal argument of Corollary~\ref{doubletrouble}.
For an introduction to sectorial forms we refer to~\cite{MR}. The theory of generalized forms can be found in~\cite{S1}.

Our notion of generalized convergence below generalizes M. Hino's conditions \cite{H} to the Kuwae--Shioya framework. Although we shall use the case of generalized forms later, we formulate this notion for
sectorial forms here:
\begin{defi}[Generalized convergence]\label{genint}
A sequence of coercive closed forms $\{(\Escr^n,\Dscr(\Escr^n))\}$ defined on $H_n$, $n\in\N$, \textit{converges to} a coercive closed form $(\Escr,\Dscr(\Escr))$ on $H$ \textit{in the generalized sense} if the
sector constants $K_n$ of the $\Escr^n$'s are uniformly bounded and
$C\subset\Dscr(\Escr)$ is dense w.r.t. $\tilde{\Escr}^{1/2}_1$, and if the following two conditions hold:
\begin{enumerate}
\item[(F1)] For every weakly convergent sequence $u_n\to u$ with $\varliminf_n\tilde{\Escr}_1^{n}(u_n)<\infty$ it follows that $u\in\Dscr(\Escr)$.
\item[(F2)] For every $w\in C$ and for every weakly convergent sequence $u_n\to u$ with $u_n\in\Dscr(\Escr^n)$, $u\in\Dscr(\Escr)$ there exists a strongly convergent sequence $w_n\to w$ such that
\[\lim_n\Escr^{n}(w_n,u_n)=\Escr(w,u)\]
\end{enumerate}
\end{defi}
Here, the symmetric part $\tilde{\Escr}$ and the sector constants are defined as in Subsection~\ref{ss:nonsym}. Definition~\ref{genint} is the special case of Definition~\ref{Fs} below for coercive closed forms with uniformly bounded sector constants. In Definition~\ref{Fs}, condition (F1) is formulated with the functional
\[\Theta^n(u):=\sup\left\{\Escr^n_1(v,u)\;\vert\;v\in\Dscr(\Escr^n),\;\tilde{\Escr}^n_1(v,v)=1\right\},\quad u\in\Dscr(\Escr^n),\]
which is the norm of the linear functional $\Escr^n_1(\cdot,u)$ on the form domain $\Dscr(\Escr^n)$, equipped with the norm $(\tilde{\Escr}^n_1)^{1/2}$. In the non-symmetric case, $\Theta^n(u)$ takes over the r\^ole of the energy $\tilde{\Escr}^n_1(u,u)^{1/2}$, since it also controls the off-diagonal values $\Escr^n_1(v,u)$. For coercive closed forms with sector constants $K_n$, we have
\[\tilde{\Escr}^n_1(u,u)^{1/2}\le\Theta^n(u)\le K_n\,\tilde{\Escr}^n_1(u,u)^{1/2},\quad u\in\Dscr(\Escr^n),\]
by Lemma~\ref{PsiEquiv} (ii) and the weak sector condition. Hence, for uniformly bounded sector constants, the above condition (F1) is equivalent to condition (F1) of Definition~\ref{Fs}, and the density hypotheses of both definitions coincide.
We stress that, already for Mosco convergence, the density of $C$ in the domain of the limiting form is essential in order to obtain convergence. In the examples of Section~\ref{sec:mosco}, it is ensured by the assumption that the minimal and the maximal extensions of the limiting form coincide. This is a uniqueness property of the limiting form, which is closely related to Markov uniqueness. For these notions, we refer to~\cite{RZ3}, \cite{RZ4} and \cite{Eb}.
We prove in the Kuwae--Shioya framework, more precisely, in the setting of Subsubsection~\ref{sss:sequence} under Assumption~\ref{infdim}, that this is equivalent to the strong convergence of the associated resolvents and the weak
convergence of the associated co-resolvents. This was first shown in the author's diploma thesis~\cite{Toe1}. The proof relies on the results of M. Hino~\cite{H}. An alternative proof by a different method can be found in the author's doctoral dissertation~\cite[Theorem 7.15]{Toe2}.
The condition that the sector constants are uniformly bounded can be replaced by a weaker version of (F1).
This characterization, which is Theorem~\ref{Hino}, is the central result of the paper. Condition (F1) corresponds to (M1): weak limits of sequences with bounded energy belong to the domain of the limiting form. Condition (F2) corresponds to (M2), i.e., to the existence of recovery sequences. However, it has to be tested against arbitrary weakly convergent sequences $\{u_n\}$, since in the non-symmetric case the off-diagonal values $\Escr^n(w_n,u_n)$ are not controlled by the diagonal. For symmetric forms, (F1) and (F2) are just another characterization of Mosco convergence, see Remark~\ref{symmosco}.

\subsection*{Relation to the author's diploma thesis and to the later literature}

This paper is a revised version of the author's diploma thesis~\cite{Toe1}, which was written at Universit\"at Bielefeld in 2006 and circulated as BiBoS Preprint E06-09-234. The results of Sections~\ref{sec:general}--\ref{sec:stoch} and of Appendix A, including Theorems~\ref{Kato}, \ref{Hino}, \ref{convergenceaij}, \ref{infaij} and \ref{lawconv} and the characterization of generalized convergence by strong graph convergence in Subsection~\ref{ss:graph}, already appear in~\cite{Toe1}, with the same structure and, apart from the corrections listed below, the same proofs. Apart from Subsection~\ref{ss:subMarkov}, which is not contained in~\cite{Toe1}, the present paper does not contain substantial new mathematical results beyond~\cite{Toe1}. Its purpose is a corrected, updated and more accessible presentation. Compared to~\cite{Toe1}, the following changes have been made:
\begin{enumerate}[(1)]
\item The convention that convergence to $H_\infty$ is only tested along indices tending to infinity has been made explicit in Definition~\ref{Hilbertconv} and Remark~\ref{convention}, and the infinite dimensionality of the Hilbert spaces, which is needed for Theorem~\ref{Kolmetric}, has been made the standing Assumption~\ref{infdim}.
\item The constant in Lemma~\ref{PsiEquiv} (i) and the running index in Corollary~\ref{doubletrouble} have been corrected.
\item In~\cite{Toe1}, Theorem~\ref{closability} asserted the existence of a lower bound $\a$ independent of $n$, which does not follow from the hypotheses, as Example~\ref{nounifshift} shows. Theorem~\ref{closability} is now stated for each form separately, a uniform lower bound is assumed explicitly in Assumption~\ref{unifshift}, and Theorem~\ref{convergenceaij} and its proof have been adapted accordingly. In particular, no uniform bound on the sector constants is needed there.
\item Theorem~\ref{lawconv} has been reformulated: the state space is a separable Banach space or $\R^d$, the processes are realized on the canonical path space $C([0,\infty);E)$, and the argument showing that a certain family of functions separates the points of $E^J$ has been corrected, since the family of products of the form $\sin l_1(x_1)\cdots\sin l_k(x_k)$ used in~\cite{Toe1} does not separate the points of $E^J$ for $k\ge 2$. Remark~\ref{cadlag} explains the modifications for c\`adl\`ag processes.
\item The discussion of~\cite{KS2} and \cite{Toe2}, including Remark \ref{seqtop} on the sequential topology, has been added, and the bibliography has been updated. Subsection~\ref{ss:subMarkov} on the preservation of the sub-Markov property, which is not contained in~\cite{Toe1}, has been added.
\item A further careful check of all proofs led to the following additional corrections. The proof of (F1b) in Theorems~\ref{convergenceaij} and \ref{infaij} identified two weak limits without justification. Condition (M1) for the symmetric parts is now assumed, see Remark~\ref{M1remark}. In Theorem~\ref{infaij}, uniform bounds and uniform tail conditions on the coefficients are now assumed in order to justify the interchange of limits and infinite sums, and the hypothesis on the logarithmic derivatives has been corrected, also in Proposition~\ref{weakbeta}. In Subsection~\ref{idsc}, the construction of the orthonormal bases $\{e_n^i\}$ and the domain of the maximal gradient form have been corrected, and the coincidence of minimal and maximal extensions is now stated explicitly in Remark~\ref{weakderi}. In Section~\ref{sec:general} and Appendix A, Theorem~\ref{Kolmetric} is now stated along general index sequences, auxiliary statements on subsequences have been added in Corollary~\ref{exstrongconv} (iii), (iv) and are used where hypotheses on full sequences were applied to subsequences, and several statements and proofs, in particular those of Lemmas~\ref{opnormconv}, \ref{compactB}, \ref{metricdoubletrouble} and Theorem~\ref{KSConv}, have been corrected in detail. Moreover, in Theorem~\ref{infaij} a summability condition of trace class type on the coefficient matrices, which contradicted the uniform ellipticity in infinite dimensions, has been removed, the index convention has been made consistent, and the bounds $\sup_nK_n<\infty$ and $c_0Q_1^n\le\tilde{\Escr}^n_1$ are now derived rather than assumed. In Section~\ref{sec:mosco}, the full support of $\mu_n$, the condition $\lrbr{l,h}_X\neq0$ and the local integrability of $1/\rho$ for the limit measure are now assumed explicitly in Subsection~\ref{idsc}, the hypotheses of Theorem~\ref{gradientmosco} have been completed, the Mosco limit in Theorem~\ref{Kolconvfinite} has been identified as $\Escr_0=\Escr_+$, and Remark~\ref{weakisstrong} has been added. In Section~\ref{sec:stoch}, the notion of association of a right process with a form is now defined explicitly.
\end{enumerate}

Since 2006, the Kuwae--Shioya framework has been developed further. In his PhD thesis~\cite{Toe2}, the author extended it to varying Banach spaces and showed that strong convergence generates a sequential topology~\cite[Theorem 5.32]{Toe2}. In~\cite{Toe2}, versions of Theorem~\ref{Kato}, including the equivalence with strong graph convergence of the generators, and of Theorem~\ref{Hino} are also obtained by a different method, see~\cite[Proposition 7.2, Theorem 7.10 and Theorem 7.15]{Toe2}. We do not use these results here, apart from the remark on the sequential topology, and keep the direct proofs of~\cite{Toe1}. J.-U. L\"obus~\cite{Lob} developed a related Mosco-type convergence of non-symmetric bilinear forms with applications to the weak convergence of $n$-particle systems. K. Suzuki~\cite{Suz} studied the convergence of non-symmetric diffusion processes on RCD spaces, and H. Ito~\cite{Ito} used the framework of~\cite{Toe1} to obtain discrete approximations of semi-Dirichlet forms with non-symmetric perturbations of resistance forms. The diploma thesis~\cite{Toe1} has moreover been cited in works on $N$-particle Langevin dynamics and their scaling limits~\cite{CG10,CG11,GN20}, on averaging principles for diffusion processes~\cite{BvR14}, on Mosco convergence of gradient forms~\cite{GW24,GW25,Wit23}, on the approximation of partial differential equations on resistance spaces~\cite{HM22,Mei20}, on Venttsel' problems in irregular domains~\cite{LRV16}, on Markov chain approximations of pure jump processes~\cite{MSS18} and on the stability of invariant measures of diffusion processes under synthetic lower Ricci curvature bounds~\cite{Suz22}. For the general theory of semi-Dirichlet forms and the associated Markov processes, we refer to~\cite{Osh}, and for the stochastic calculus of generalized Dirichlet forms to~\cite{Tr} and \cite{Tr2}.

\subsection*{Structure of the paper}

Let us now give a brief overview of the structure of this work and the main results.

In Subsection~\ref{ss:spectral} we develop the Kuwae--Shioya framework of convergence of spectral structures along a sequence of Hilbert spaces.
Subsubsection~\ref{sss:general} introduces the original, rather abstract setting, whereas in view of our later applications we will use
the setting of Subsubsection~\ref{sss:sequence} exclusively, which is the setting described above, seen from an abstract point of view. Here we obtain
as a main result in Theorem~\ref{Kolmetric} the complete and separable metric structure of our space $\Hscr$ and a sequence of isometric isomorphisms
$\{\Psi_n:H\to H_n\}$ with the property that $u_n\to u$ strongly if and only if $\lim_n\norm{\Psi_n(u)-u_n}_{H_n}=0$. The proof is essentially
due to A.V. Kolesnikov~\cite[Proposition 7.2]{K1}, but has been completed, simplified and entirely rewritten in~\cite{Toe1}. The remaining lemmas and propositions of Subsection~\ref{ss:spectral} are taken from K. Kuwae and T. Shioya~\cite[Section 2]{KS}, unless stated otherwise in the text, as for Lemma~\ref{OpPhiConv}, Theorem~\ref{Kato} and Corollary~\ref{exstrongconv}. All these proofs, partially completed and rewritten, can be found in Appendix A, as well as some of the proofs of A.V. Kolesnikov in~\cite{K1} and \cite{K2} and one proof of H. Attouch in~\cite{A}. We suggest reading Subsections~\ref{ss:spectral}--\ref{ss:symmetric} and Appendix A in parallel.

In Subsubsection~\ref{sss:semigroups} we prove a result which was new in the diploma thesis~\cite{Toe1} from 2006, namely the generalization of T. Kato's theorem mentioned above to the Kuwae--Shioya framework.
It gives the equivalence of strong convergence of resolvents and semigroups.

Subsection~\ref{ss:symmetric} is taken from~\cite{KS} and \cite{K1}.

In Subsection~\ref{ss:nonsym} we arrive at the main result for generalized non-symmetric forms. The result, Theorem~\ref{Hino}, has been described above.
In this generality, the conditions differ from those of Definition~\ref{genint}, and they are necessary and sufficient for the
convergence of resolvents. After that we examine the relation between Mosco convergence and generalized convergence. These results were also new in the diploma thesis~\cite{Toe1} from 2006.

Subsection~\ref{ss:graph} examines the connection between generalized convergence and strong graph convergence of the infinitesimal generators of our forms. Their equivalence
is proved as well. Additionally, we give a characterization of strong graph convergence for closed operators associated with generalized forms.
All this is done in the Kuwae--Shioya framework. In the non-symmetric case, these results first appeared in~\cite{Toe1}. Note, however, that strong graph convergence may be difficult to verify in applications, since the domains of the generators are often not explicitly known.

Subsection~\ref{ss:subMarkov} shows that the sub-Markov property of semigroups is preserved under strong convergence along a sequence of $L^2$-spaces with weakly converging finite reference measures.

Section~\ref{sec:mosco} reviews the results of~\cite{K1} and \cite{K2} for the reasons mentioned above.

Section~\ref{sec:genconv} contains the main application results, which first appeared in~\cite{Toe1}. We consider finite and infinite dimensional elliptic and sectorial
$a_{i,j}$-forms and give sufficient conditions for the generalized convergence of these forms.

Section~\ref{sec:stoch} gives a short application to stochastics, namely, that the associated processes converge weakly, in the sense of the generating
path measures, provided the associated forms converge in the generalized sense and the collection of path measures is tight. This result was also new in the diploma thesis~\cite{Toe1} from 2006, and we remark that here the concept of changing reference measures is particularly useful.

Appendix A is a collection of proofs from some reference papers, which may have been completed and rewritten.

\subsubsection*{Origin of the results}
To summarize, the results of this paper have the following origin. Further details are given in the preceding subsection.
\begin{enumerate}[(a)]
\item Results recalled from Kuwae and Shioya~\cite{KS} and from Kolesnikov~\cite{K1,K2}: the basic properties of convergence along $\{H_n\}$ in Subsection~\ref{ss:spectral}, the symmetric theory in Subsection~\ref{ss:symmetric}, and Section~\ref{sec:mosco}.
\item Results adapted or completed from the literature: Theorem~\ref{Kolmetric}, which is based on~\cite[Proposition 7.2]{K1}, and the extension of Hino's theorem~\cite[Theorem 3.1]{H} to varying Hilbert spaces in Theorem~\ref{Hino}.
\item Results which were new in the diploma thesis~\cite{Toe1} from 2006: Theorem~\ref{Kato}, the characterization in Theorem~\ref{Hino} and its relation to Mosco convergence, the strong graph convergence theorem in Subsection~\ref{ss:graph}, the applications in Section~\ref{sec:genconv} and Theorem~\ref{lawconv}.
\item Results which are not contained in~\cite{Toe1}: the preservation of the sub-Markov property in Subsection~\ref{ss:subMarkov}.
\end{enumerate}

\subsection*{Acknowledgements}

This work is based on the diploma thesis\footnote{A German diploma thesis is equivalent to a master's thesis.} \cite{Toe1} by the author from 2006, supervised by Michael R\"ockner at Bielefeld University. Large parts of the thesis were written during a stay at Purdue University, West Lafayette, Indiana in 2006.
The reason for publishing it now is the ongoing interest of the scientific community. The original diploma thesis had been cited 15 times by May 2026. The author is grateful to the anonymous referees for a very careful reading and for many valuable suggestions, including Example~\ref{nounifshift}.

\newpage
\subsection*{Notation}\mbox{}\par\nopagebreak
\noindent\begin{minipage}{\textwidth}
The main notation is summarized in the following table.
\begin{center}\small
\begin{tabular}{@{}>{\raggedright\arraybackslash}p{0.27\textwidth}>{\raggedright\arraybackslash}p{0.42\textwidth}l@{}}
Symbol & Meaning & Introduced in\\
\hline
$H_n$, $H=H_\infty$ & Hilbert spaces with $H_n\to H$ & Definition~\ref{Hilbertconv}\\
$C$ & dense subspace of $H$ & Definition~\ref{Hilbertconv}\\
$\Phi_n$ & linear maps $C\to H_n$ & Definition~\ref{Hilbertconv}\\
$\Hscr$ & disjoint union of the spaces $H_n$ and $H$ & Subsubsection~\ref{sss:sequence}\\
$u_n\to u$ & strong and weak convergence in $\Hscr$ & Definitions~\ref{strongconv}, \ref{weakconv}\\
$\Psi_n$ & isometric isomorphisms $H\to H_n$ & Theorem~\ref{Kolmetric}\\
$d_\Hscr$ & metric on $\Hscr$ & Theorem~\ref{Kolmetric}\\
$B_n\to B$ & strong and weak convergence of bounded operators & Definition~\ref{Operatorconvergence}\\
$G^n_\a$, $G_\a$ & resolvents & Subsubsection~\ref{sss:semigroups}\\
$T^n_t$, $T_t$ & semigroups & Subsubsection~\ref{sss:semigroups}\\
$A^n$, $A$ & infinitesimal generators & Subsubsection~\ref{sss:semigroups}\\
$\Escr^n\to\Escr$ & Mosco convergence of symmetric closed forms & Definition~\ref{Mosco}\\
$\Hscr_\Escr$, $\Psi^\Escr_n$ & disjoint union of the form domains and the corresponding isometries & Definition~\ref{Dconv}\\
$(\Ascr^n,\Vscr_n)$ & coercive closed forms & Subsubsection~\ref{sss:genforms}\\
$\tilde{\Ascr}^n$, $\check{\Ascr}^n$ & symmetric and antisymmetric parts of $\Ascr^n$ & Subsubsection~\ref{sss:genforms}\\
$K_n$ & sector constants & Subsubsection~\ref{sss:genforms}\\
$\Lambda^n$ & linear operator defining the generalized form & Subsubsection~\ref{sss:genforms}\\
$\Fscr_n$, $\hat{\Fscr}_n$ & domains of the closure of $\Lambda^n$ and of the adjoint operator $\hat{\Lambda}^n$ in $\Vscr_n$ & Subsubsection~\ref{sss:genforms}\\
$\Escr^n$, $\hat{\Escr}^n$ & generalized forms and co-forms & Subsubsection~\ref{sss:genforms}\\
$\hat{G}^n_\a$, $\hat{T}^n_t$ & co-resolvents and co-semigroups & Proposition~\ref{Resolvent}\\
$\Escr^{(\b)}$ & approximate forms & Proposition~\ref{BetaMR}\\
$\Theta^n$ & dual norm of $\Escr^n_1(\cdot,u)$ on $\Vscr_n$ & Lemma~\ref{PsiEquiv}\\
(F1), (F2), (F2'), (R), (CR) & convergence conditions & Definition~\ref{Fs}\\
$\Gamma_\infty$ & strong graph limit & Subsection~\ref{ss:graph}\\
$\phi_\eps$, $\psi_\eps$ & smooth truncations of the identity & Definition~\ref{smoothtrunc}\\
$\mu_n=\sigma_n\,dx$ & speed measures on $\R^d$ & Assumption~\ref{speed}\\
$a^n_{i,j}$, $b^n_i$, $d^n_i$, $c^n$ & coefficients of the $a_{i,j}$-forms & Subsection~\ref{ss:finaij}\\
$X\subset E$ & Hilbert space densely embedded into a Banach space & Subsection~\ref{idsc}\\
$\FCE$, $\FCb$ & cylindrical test functions & Subsections~\ref{idsc}, \ref{ss:infaij}\\
$\nabla_X$ & gradient along $X$ & Equation~\eqref{nabla}\\
$\beta^{\mu_n}_h$ & logarithmic derivative of $\mu_n$ along $h$ & Subsection~\ref{idsc}\\
$\M^{(n)}$, $\P^{(n)}$ & right processes and their laws on path space & Subsection~\ref{ss:laws}\\
$m_n$, $\mu_n$ & reference measures and initial distributions in Section~\ref{sec:stoch} & Subsection~\ref{ss:laws}
\end{tabular}
\end{center}
\end{minipage}

\section{General functional analytic theory}\label{sec:general}

\subsection{Convergence of spectral structures}\label{ss:spectral}

We start with some notation. Let $\N$, $\Z$, $\Q$, $\R$, $\C$ denote the sets of natural, integer, rational, real and complex numbers, respectively.
We write $\K:=\R$ or $\C$ if we do not want to specify whether we use the real or complex numbers. Let $\ol{\R}:=\R\cup\{-\infty\}\cup\{+\infty\}$ and
$\ol{\C}:=\C\cup\{\infty\}$, i.e., the extended real numbers and the extended complex plane. For two sets $A,B$ we write $A\dot\cup B:=A\cup B$,
if $A$ and $B$ are disjoint, i.e., $A\cap B=\emptyset$. This notation also extends to arbitrary unions of pairwise disjoint sets.
For a Banach space $E$ we denote by $E'$ its topological dual. For Banach spaces $E,F$, we denote the set of all bounded
linear operators from $E$ to $F$ by $\Lscr (E,F)$ with operator norm
$\norm{\;}_{\Lscr (E,F)}$. For convenience we set $\Lscr(E):=\Lscr(E,E)$. We abbreviate $\a:=\a\Id$ for any $\a\in\K$.

We follow a framework developed by K. Kuwae and T. Shioya in~\cite{KS}. It provides a functional
analytic theory introducing convergence along a sequence of Hilbert spaces, and applications to convergence
of spectral structures and forms. The most important proofs
are repeated in Appendix A.

\subsubsection{Convergence of Hilbert spaces - the general case}\label{sss:general}

We first describe the general framework, which we restrict to a countable index set afterwards.
Let $\Nscr$ be an arbitrary index set and $\{H_\nu\;\vert\;\nu\in\Nscr\}$ a family of separable Hilbert spaces over $\K$. Assume
that there exists a family
\[\{\Phi_{\nu,\mu}:\Cscr_\nu\to H_\mu\;\vert\;\nu,\mu\in\Nscr\}\]
of linear maps with dense linear domains
$\Cscr_\nu\subset H_\nu$ such that $\Phi_{\nu,\nu}$ for each $\nu\in\Nscr$ is the identity operator on $\Cscr_\nu$.
Assume that
$\{H_\nu\;\vert\;\nu\in\Nscr\}$, or equivalently $\Nscr$, carries a topology, not necessarily Hausdorff, such that a sequence $\{H_{\nu_n}\}_{n\in\N}$
converges to some $H_\nu$, $\nu\in\Nscr$ if and only if for any $u\in\Cscr_\nu$
\[\lim_n\norm{\Phi_{\nu,\nu_n}u}_{H_{\nu_n}}=\norm{\Phi_{\nu,\nu}u}_{H_\nu}\quad(=\norm{u}_{H_\nu}).\]

Now we can define a topology on the disjoint union $\Hscr:=\disjoint_{\nu\in\Nscr}H_\nu$ as follows.
\begin{defi}\label{abstractstrconv}
Assume that a sequence $\{H_{\nu_n}\}$ converges to an $H_\nu$.
We say that a sequence $\{u_n\}_{n\in\N}$ with $u_n\in H_{\nu_n}$ \textit{converges strongly} to a vector $u\in H_\nu$ if
for one sequence, and \textbf{hence for all} sequences, $\{\tilde{u}_m\}_{m\in\N}\subset\Cscr_\nu$ tending to $u$ in $H_\nu$ we have
\[\lim_m \varlimsup_n \norm{\Phi_{\nu,\nu_n}\tilde{u}_m -u_n}_{H_{\nu_n}}=0.\]
\end{defi}
This topology is called the \textit{strong topology} on $\Hscr$. Note that this notion of convergence depends explicitly on the sequence $\{\nu_n\}$ and can only be
defined in a reasonable way
if $H_{\nu_n}\to H_\nu$. ``Reasonable'' means that ``\textbf{hence all}'' in
the above definition holds,
which is a consequence of the easily provable fact that $H_{\nu_n}\to H_\nu$ implies that for every
$\{\tilde{u}_m\}\subset \Cscr_\nu$ such that $\tilde{u}_m\to 0\in H_\nu$ we have
\[\lim_m\varlimsup_n\norm{\Phi_{\nu,\nu_n}\tilde{u}_m}_{H_{\nu_n}}=0.\]
Moreover, the notion of a subsequence is defined, since, provided $H_{\nu_n}\to H_\nu$, clearly
$H_{\nu_{n_k}}\to H_\nu$ for any sequence of natural numbers $n_k\uparrow \infty$, $n_{k+1}>n_k$. It can easily be
proved that $u_n\to u$, $u_n\in H_{\nu_n}$, $u\in H_\nu$ strongly if and only if every subsequence $u_{n_k}\to u$, $u_{n_k}\in H_{\nu_{n_k}}$ strongly.
On the other hand, it does not make sense to ask for convergence
of a sequence $\{u_n\}\subset\Hscr$ with $u_n\in H_{\nu_n}$ for some $\nu_n$ such that the sequence $\{H_{\nu_n}\}_{n\in\N}$
does not have a limit.

\begin{lem}\label{hausimdorf}
The strong topology on $\Hscr$ is Hausdorff if and only if $\{H_\nu\;\vert\;\nu\in\Nscr\}$
is Hausdorff.\footnote{Here, a topology defined by convergent sequences is called Hausdorff if limits of convergent sequences are unique.}
\end{lem}
\begin{proof}
This is a more detailed proof than the original one found in~\cite[Corollary 2.2]{KS}.
To prove the ``if'' part assume that $\{H_\nu\;\vert\;\nu\in\Nscr\}$ is Hausdorff, i.e.,
a convergent sequence $\{H_{\nu_n}\}$ has at most one limit point $H_{\nu_0}$. Now let
$u_n\in H_{\nu_n}$, $u,v\in H_{\nu_0}$ such that $u_n\to u$ and $u_n\to v$ strongly. If we can
prove that $u=v$, we are done. Let $\{\tilde{u}_m\},\{\tilde{v}_m\}\subset\Cscr_{\nu_0}$ as in
Definition~\ref{abstractstrconv} such that $\lim_m\norm{\tilde{u}_m-u}_{H_{\nu_0}}=0$,
$\lim_m\norm{\tilde{v}_m-v}_{H_{\nu_0}}=0$,
\[\lim_m \varlimsup_n \norm{\Phi_{\nu_0,\nu_n}\tilde{u}_m -u_n}_{H_{\nu_n}}=0\quad\text{and}\quad
\lim_m \varlimsup_n \norm{\Phi_{\nu_0,\nu_n}\tilde{v}_m -u_n}_{H_{\nu_n}}=0.\]
Now,
\begin{align*}
\norm{u-v}_{H_{\nu_0}}&\le\varlimsup_m\left[\norm{u-\tilde{u}_m}_{H_{\nu_0}}+\norm{\tilde{u}_m-\tilde{v}_m}_{H_{\nu_0}}+\norm{\tilde{v}_m-v}_{H_{\nu_0}}\right]\\
&=\varlimsup_m\lim_n\norm{\Phi_{\nu_0,\nu_n}(\tilde{u}_m-\tilde{v}_m)}_{H_{\nu_n}}=\varlimsup_m\lim_n\norm{\Phi_{\nu_0,\nu_n}\tilde{u}_m-\Phi_{\nu_0,\nu_n}\tilde{v}_m}_{H_{\nu_n}}\\
&\le\varlimsup_m\varlimsup_n\norm{\Phi_{\nu_0,\nu_n}\tilde{u}_m-u_n}_{H_{\nu_n}}+\varlimsup_m\varlimsup_n\norm{u_n-\Phi_{\nu_0,\nu_n}\tilde{v}_m}_{H_{\nu_n}}\\
&=0.
\end{align*}
Hence $u=v$. 

To prove the ``only if'' part, assume that $\{H_\nu\;\vert\;\nu\in\Nscr\}$ is
not Hausdorff, i.e., there exists a sequence of Hilbert spaces $\{H_{\nu_n}\}$ with $H_{\nu_n}\to H_{\nu_1}$ and $H_{\nu_n}\to H_{\nu_2}$ where $H_{\nu_1}\not=H_{\nu_2}$. Set $u_n:=0\in H_{\nu_n}$, $n\in\N$. Then
clearly $u_n\to 0\in H_{\nu_1}$ strongly and $u_n\to 0\in H_{\nu_2}$ strongly.
Hence $\Hscr$ is not Hausdorff which proves the ``only if'' part.
\end{proof}
Note that the fact that for every $u\in H_{\nu_0}$ there exists a strongly convergent sequence $\{u_n\}$ along
$H_{\nu_n}$ is verified for a sequence of Hilbert spaces below in Corollary~\ref{exstrongconv}, but it is not clear whether it holds
in the general case. For a sequence of Hilbert spaces and under Assumption~\ref{infdim}, Theorem~\ref{Kolmetric} even shows that $\Hscr$ is \textit{Polish}.

\subsubsection{Convergence of a sequence of Hilbert spaces}\label{sss:sequence}

In the following, we consider a single sequence $H_n$, $n\in\N$, converging to a limit space $H_\infty$. This covers all applications in this paper, whereas the general case of Subsubsection~\ref{sss:general}, with an arbitrary topological index set, is only included for context.

\textbf{From now on}, we shall consider the case that $\Nscr=\N\cup\{\infty\}$. Therefore, let $H_n$, $n\in\N$, and $H_\infty$
be real separable Hilbert spaces with inner products $\deli{\;}{\;}_{H_n}$, $\deli{\;}{\;}_{H_\infty}$ and
norms $\lrnorm{\;}_{H_n}=\deli{\;}{\;}_{H_n}^{1/2}$, $\lrnorm{\;}_{H_\infty}=\deli{\;}{\;}_{H_\infty}^{1/2}$, respectively.
As a special case of the above, we set $\Cscr_n:=H_n$, $n\not=\infty$, and fix some dense linear subspace $\Cscr_\infty\subset H_\infty$. Furthermore, we
fix some injective linear operators $\Phi_{\infty,n}:\Cscr_\infty\to H_n$, and for $m,n\in\N$ set $\Phi_{m,n}=\Id$, if $m=n$ and $\Phi_{m,n}\equiv 0$ if $m\not=n$. Moreover, $\Phi_{\infty,\infty}:=\Id_{\Cscr_\infty}$. We assume that $H_n\neq\{0\}$ for all $n\in\N$. Because of the additional requirement $n_k\to\infty$ in Definition~\ref{Hilbertconv} below, this is not literally a special case of the general setting. However, Lemma~\ref{hausimdorf} remains valid, since its proof only uses that convergence of Hilbert spaces implies \eqref{Phiuconv}.
Let $\Hscr:=\disjoint_{n\in\N}H_n\dot\cup H_\infty$ be the disjoint union of Hilbert spaces. As above
we define:
\begin{defi}[Convergence of Hilbert spaces]\label{Hilbertconv}
Let $\{n_k\}$ be a not necessarily increasing sequence in $\N\cup\{\infty\}$.
A sequence of Hilbert spaces
\[\{H_{n_k}\}\subset\left\{H_n\;\vert\;n\in\N\cup\{\infty\}\right\}\]
is said to \textit{converge} to a Hilbert space
$H_{n_0}$, $n_0\in\N\cup\{\infty\}$ if
\begin{equation}\label{Phiuconv}
\lim\limits_{k\to\infty}\lrnorm{\Phi_{n_0,n_k}u}_{H_{n_k}}=\lrnorm{u}_{H_{n_0}}
\end{equation}
for every $u\in\Cscr_{n_0}$, and if, in the case $n_0=\infty$, in addition $n_k\to\infty$ as $k\to\infty$, that is, for every $N\in\N$ there is $k_0\in\N$ such that $n_k\ge N$ for all $k\ge k_0$, where $n_k=\infty$ is allowed. Here we set $\Phi_{n,\infty}:\equiv 0$ for $n\in\N$.
\end{defi}

\begin{rem}\label{convention}
The additional requirement $n_k\to\infty$ in the case $n_0=\infty$ encodes the distinguished limit $n\to\infty$. Without it, Definition~\ref{Hilbertconv} would allow degenerate situations: if, for example, $H_n=H_\infty$ and $\Phi_{\infty,n}=\Id$ for all $n$, then the constant sequence $n_k=1$ would satisfy \eqref{Phiuconv} for $n_0=\infty$ as well as for $n_0=1$, and the topology on $\{H_n\;\vert\;n\in\N\cup\{\infty\}\}$ would not be Hausdorff. With the above convention, convergence to $H_\infty$ is only tested along indices tending to infinity, whereas convergence to $H_{n_0}$, $n_0\in\N$, is ordinary convergence inside the fiber $H_{n_0}$, see Remark~\ref{2cases} below.
\end{rem}

\begin{rem}\label{seqtop}
In~\cite[Theorem 5.32]{Toe2}, the author has shown that Kuwae--Shioya convergence generates a sequential topology in the sense of $\mathscr{S}^\ast$-convergence classes, see~\cite[Definition A.5]{Toe2}, so that, contrary to the remark of Kuwae and Shioya that nets are needed, it suffices to consider sequences.
\end{rem}

From now on, we make the following assumption.
\begin{assumption}
$\{H_n\}_{n\in\N}$ converges to $H_\infty$.
\end{assumption}
For the motivation of this assumption and the relation to a sequence of $L^2$-spaces and forms defined on them we refer either to the
introduction or to the later sections.
\begin{defi}[Strong convergence]\label{strongconv}
Let $\{n_k\}$ be a not necessarily increasing sequence in $\N\cup\{\infty\}$. Let $n_0\in\N\cup\{\infty\}$. Assume that $H_{n_k}\to H_{n_0}$ in the above sense.
$\{u_{k}\}$, $u_{k}\in H_{n_k}$ is said to \textit{strongly converge} to some $u\in H_{n_0}$ if there exist $\{\tilde{u}_m\}\subset \Cscr_{n_0}$
such that
\begin{gather}\label{STRONG1}
\lim\limits_m\lrnorm{\tilde{u}_m-u}_{H_{n_0}}=0,\\ \label{STRONG2}
\lim\limits_m\varlimsup\limits_k\lrnorm{\Phi_{n_0,n_k}\tilde{u}_m-u_{k}}_{H_{n_k}}=0.
\end{gather}
Then \eqref{STRONG2} even holds for every sequence $\{\tilde{v}_m\}\subset\Cscr_{n_0}$ with
$\tilde{v}_m\to u$ in $H_{n_0}$.
\end{defi}
\begin{defi}[Weak convergence]\label{weakconv}
Let $\{n_k\}$ be a not necessarily increasing sequence in $\N\cup\{\infty\}$. Let $n_0\in\N\cup\{\infty\}$. Assume that $H_{n_k}\to H_{n_0}$ in the above sense.
$\{u_{k}\}$, $u_{k}\in H_{n_k}$ is said to \textit{weakly converge} to $u\in H_{n_0}$ if
\begin{equation}\label{WEAK}
\deli{u_k}{v_k}_{H_{n_k}}\to\deli{u}{v}_{H_{n_0}}
\end{equation}
for every sequence $\{v_k\}$, $v_k\in H_{n_k}$ strongly convergent to $v\in H_{n_0}$.
\end{defi}
Strong convergence
creates a topology on $\Hscr=\disjoint_{n\in\N}H_n\dot\cup H_\infty$, called the \textit{strong topology}, and weak convergence creates the \textit{weak topology}.

\begin{rem}\label{2cases}
We point out that this is a very special case, which, in a way, restricts strong and weak convergence of ``interesting'' sequences
to those converging along the sequence of Hilbert spaces $H_n\to H_\infty$. To see what this means, consider the following
statements, which can easily be proved, taking into account that for $n,m\in\N$, $\Cscr_{n}=H_n$, $\Phi_{m,n}=\Id$ if $m=n$ and $\Phi_{m,n}\equiv 0$,
if $m\not=n$,
\begin{enumerate}[(1)]
\item It can easily be seen that the topology of $\left\{H_n\;\vert\;n\in\N\cup\{\infty\}\right\}$ is Hausdorff, which also follows from Lemma~\ref{hausimdorf} and Theorem~\ref{Kolmetric}. It is also second countable, see \textup{\cite[Lemma 2.13]{KS}}.
\item Let $n_0\in\N$, that is, $n_0\not=\infty$. Then $H_{n_k}\to H_{n_0}$ if and only if $\exists k_0\ge 1$ with $n_k=n_0$ $\forall k\ge k_0$.
\item Let $n_0\in\N$, that is, $n_0\not=\infty$, and assume $H_{n_k}\to H_{n_0}$. Then $u_k\to u$, $u_k\in H_{n_k}$, $u\in H_{n_0}$ strongly if and only if
$\exists k_0\ge 1$ such that $u_k\in H_{n_0}$ $\forall k\ge k_0$ and
\[\lim_{k\to\infty,k\ge k_0}\norm{u_{k}-u}_{H_{n_0}}=0.\]
A similar characterization holds for the weak convergence.
\end{enumerate}
As a consequence, in our particular case, depending on the ``limiting'' Hilbert space,
one should understand strong and weak convergence either in the ``usual'' way, ending up in one ``final'' Hilbert space, or in the
sense ``along $H_n$'' ending in $H_\infty$, which leads to a different notion of strong convergence, where ``$u_n$ converges strongly to $u$'' means that the ``distance'' between some approximating sequence $\{\tilde{u}_m\}\subset\Cscr_\infty$ and $u$ is small, and the ``distance'' between
its images in the $H_n$'s under the $\Phi_{\infty,n}$'s and the $u_n$'s is small.
\end{rem}

\textbf{From now on}, unless stated otherwise, we shall always understand strong and weak convergence along the sequence $H_n\to H_\infty$, since by the
above remark the other cases are more or less trivial. This convention also extends to the Appendix. \textbf{From now on} we shall also set
$H:=H_\infty$, $C:=\Cscr_\infty$, $\Phi_n:=\Phi_{\infty,n}$, unless this leads to confusion. Subsequences are also understood in this way. All statements below which are formulated for sequences indexed by $n\in\N$ remain valid, with the same proofs, along every sequence $\{n_k\}$ with $n_k\to\infty$ in the sense of Definition~\ref{Hilbertconv}, replacing $H_n$ by $H_{n_k}$ and $\Phi_n$ by $\Phi_{n_k}$, since their proofs only use \eqref{Phiuconv} and previously established results along the given indices, in particular the general form of \eqref{Psitwoconv} in Theorem~\ref{Kolmetric}.

Now let us start with some useful facts.

\begin{lem}\label{basicconv}
\begin{enumerate}[(1)]
\item Let $u_n\in H_n$, $n\in\N$. Then $u_n\to 0\in H$ strongly in $\Hscr$ if and only if $\norm{u_n}_{H_n}\to 0$.
\item Let $u_n,v_n\in H_n$, $n\in\N$, $u,v\in H$ such that $u_n\to u$ strongly and $v_n\to v$ strongly.
Then $\alpha u_n+\beta v_n\to \alpha u+\beta v$ strongly in $\Hscr$ for any $\alpha,\beta\in\R$.
\item Let $\{u_n\}$ be a sequence with $u_n\in H_n$, $n\in\N$ and $u_n\to u\in H$ strongly. Then $\norm{u_n}_{H_n}\to\norm{u}_H$.
In particular, the sequence of norms $\{\norm{u_n}_{H_n}\}$ of a strongly convergent sequence is bounded.
\item If $u_n,v_n\in H_n$, $n\in\N$, $u,v\in H$ such that $u_n\to u$ and $v_n\to v$ $\Hscr$-strongly, then
$(u_n,v_n)_{H_n}\to (u,v)_H$. In particular, every strongly convergent sequence converges weakly.
\end{enumerate}
\end{lem}
\begin{proof}
See Appendix A.
\end{proof}

\begin{lem}\label{twoconv}
Let $\{u_n\}$, $\{v_n\}$ be two sequences of vectors in $\Hscr$ with $u_n,v_n\in H_n$, $n\in\N$, and let $u\in H$. Suppose that $u_n\to u$
strongly in $\Hscr$.
Then $v_n\to u$ strongly in $\Hscr$ if and only if $\norm{u_n-v_n}_{H_n}\to 0$.
\end{lem}
\begin{proof}
The proof follows directly from Lemma~\ref{basicconv} (1) and (2).
\end{proof}

\begin{rem}\label{phiphi}
\begin{enumerate}[(1)]
\item One immediately obtains by definition that $\Phi_n(\phi)\to\phi$ $\Hscr$-strongly for every $\phi\in C$.
\item Note that if $v_n\in H_n$, $n\in\N$, with $v_n\to\phi\in C$ $\Hscr$-strongly, it follows from Lemma~\ref{twoconv}, or even
from ``\textbf{hence all}'' in Definition~\ref{abstractstrconv}, that
\[\lim_n\norm{v_n-\Phi_n(\phi)}_{H_n}=0.\]
We will use this later.
\end{enumerate}
\end{rem}

The next theorem is due to A.V. Kolesnikov, see~\cite[Proposition 7.2]{K1}. The proof has
been rewritten and completed for this paper.
\begin{thm}\label{Kolmetric}
Assume that all $H_n$, $n\in\N$ and $H$ are infinite dimensional and separable. Then there exists a complete separable metric $d_\Hscr$ on $\Hscr$ such that the
convergence in $d_\Hscr$ coincides with the strong convergence and there exists a bijective isometry of metric spaces $\Psi:(\Hscr,d_\Hscr)\to I\times \ell^2$,
where $I:=\{0\}\cup\bigcup_{n\in\N}\{\frac{1}{n}\}\subset\R$.

There also exists a sequence of isometric isomorphisms of Hilbert spaces
\[\Psi_n:H_\infty\to H_n\]
such that
$\Psi_\infty=\Id_{H_\infty}$
and $u_n\to u$, $u_n\in H_n$, $u\in H_\infty$ strongly in $\Hscr$ if and only if
\begin{equation}\label{Psitwoconv}
\lim_n\norm{\Psi_n u-u_n}_{H_n}=0.
\end{equation}
More generally, if $\{n_k\}$ is any sequence in $\N\cup\{\infty\}$ with $n_k\to\infty$ in the sense of Definition~\ref{Hilbertconv}, then $u_k\to u$, $u_k\in H_{n_k}$, $u\in H_\infty$, strongly in $\Hscr$ if and only if $\lim_k\norm{\Psi_{n_k}u-u_k}_{H_{n_k}}=0$.

Furthermore, if we fix $n_0\in\N$, the metric has the property that for
$f_m\in H_{n_0}$ for all $m\in\N$, $f\in H_{n_0}$ one has $\norm{f_m-f}_{H_{n_0}}\to 0$ as $m\to\infty$ if and only if $d_\Hscr(f_m,f)\to 0$ as $m\to\infty$.
\end{thm}
\begin{proof}
As in the beginning of the paper we write $H_\infty=H$.
Let us first construct an orthonormal basis $\{e_i\;\vert\;i\in\N\}$ in $H_\infty$ consisting of vectors from $C$ using a standard
orthogonalization procedure, applied to a countable subset of $C$ which is dense in $H_\infty$. Such a subset exists, since $H_\infty$ is separable and $C$ is dense.

\textbf{Step 1.} Assume for the moment that the $\Psi_n$'s have already been constructed. Then for fixed $n_0\in\N\cup\{\infty\}$,
$\{\Psi_{n_0}(e_i)\;\vert\;i\in\N\}$ clearly is an orthonormal basis of $H_{n_0}$.
Let $u\in H_n$, $v\in H_m$, $n,m\in\N\cup\{\infty\}$. We claim that the metric
\[d_\Hscr(u,v)=\sqrt{\abs{\delta_n-\delta_m}^2+\sum_{i=1}^\infty (u_i-v_i)^2},\]
where $u=\sum_{i=1}^\infty u_i \Psi_n(e_i)\in H_n$, $v=\sum_{i=1}^\infty v_i \Psi_m(e_i)\in H_m$, $\delta_n:=\frac{1}{n}$, $\delta_\infty:=0$,
is the desired one.

Let us prove that $d_\Hscr$ is a metric generating strong convergence on $\Hscr$. Obviously, $d_\Hscr$ is symmetric.
We prove that $d_\Hscr(u,v)=0$ if and only if $u=v$. The ``if'' part is trivial. To see the ``only if'' part, assume that both $\abs{\delta_n-\delta_m}^2=0$ and
$\sum_{i=1}^\infty (u_i-v_i)^2=0$. The former implies that $u,v\in H_{n_0}$ for some fixed $n_0\in\N\cup\{\infty\}$,
and the latter shows that $u=v$.
The triangle inequality follows directly from the triangle inequality in the metric space $\ell^2$.
Thus $d_\Hscr$ defines a metric on $\Hscr$. Next we prove that convergence in this metric coincides with strong convergence in $\Hscr$.
\begin{enumerate}
\item[Case 1:] ``$H_n\to H_\infty$''. We treat full sequences. For a general sequence $n_k\to\infty$ as in Definition~\ref{Hilbertconv}, replace $n$ by $n_k$ in the following argument and use the general form of \eqref{Psitwoconv} proved at the end of Step 2. Let $u_n\in H_n$, $u\in H_\infty$, $u_n\to u$ strongly. Then by \eqref{Psitwoconv}
\[\lim_n\norm{\Psi_n(u)-u_n}_{H_n}=0.\]
Let $u_i^{(n)}:=(u_n,\Psi_n(e_i))_{H_n}$, $u_i:=(u,e_i)_{H_\infty}=(\Psi_n(u),\Psi_n(e_i))_{H_n}$. Thus by Parseval's identity, see~\cite[Theorem II.6]{ReSi1},
\begin{equation}\label{Parseval}
\sum_{i=1}^\infty (u_i^{(n)}-u_i)^2=\sum_{i=1}^\infty\lrabs{(u_n-\Psi_n(u),\Psi_n(e_i))_{H_n}}^2=\norm{u_n-\Psi_n(u)}_{H_n}^2\to 0
\end{equation}
as $n\to\infty$. Clearly, $\lim_n\abs{\delta_n-\delta_\infty}^2=0$. Thus $\lim_n d_\Hscr(u_n,u)=0$.

Now assume $\lim_n d_\Hscr(u_n,u)=0$. Let $u_i^{(n)}$, $u_i$ as above. Now we can prove that
$\lim_n\norm{\Psi_n(u)-u_n}_{H_n}=0$ by reading \eqref{Parseval} backwards.
\item[Case 2:] ``$H_{n_k}\to H_{n_0}$''. Let $\{n_k\}$ be a not necessarily increasing sequence of natural numbers and $n_0\in\N$. Assume that
$H_{n_k}$ converges to $H_{n_0}$. Then there exists a $k_0\in\N$ such that
\begin{equation}\label{kgek0}
n_k=n_0\quad\forall k\ge k_0.
\end{equation}
Let $u_k\in H_{n_k}$, $u\in H_{n_0}$. We want to prove that $u_k\to u$ strongly if and only if $\lim_k d_\Hscr(u_k,u)=0$.

Assume that $u_k\to u$ strongly in $\Hscr$. By \eqref{kgek0} $\delta_{n_k}=\delta_{n_0}$ $\forall k\ge k_0$. Thus
$\lim_k\abs{\delta_{n_k}-\delta_{n_0}}^2=0$. Let $u_i^{(k)}:=(u_k,\Psi_{n_k}(e_i))_{H_{n_k}}$, $u_i:=(u,\Psi_{n_0}(e_i))_{H_{n_0}}$.
Since $u_k\in H_{n_0}$ $\forall k\ge k_0$,
we have by strong convergence
\begin{equation}\label{kk0conv}
\begin{split}
\lim_{k\to\infty,k\ge k_0}d_\Hscr(u_k,u)^2=&\lim_{k\to\infty,k\ge k_0}\sum_{i=1}^\infty (u_i^{(k)}-u_i)^2\\
=&\lim_{k\to\infty,k\ge k_0}\sum_{i=1}^\infty\lrabs{(u_k-u,\Psi_{n_0}(e_i))_{H_{n_0}}}^2\\
=&\lim_{k\to\infty,k\ge k_0}\norm{u_k-u}^2_{H_{n_0}}=0.
\end{split}
\end{equation}
Thus $\lim_k d_\Hscr(u_k,u)=0$.

Suppose now $\lim_k d_\Hscr(u_k,u)=0$. Then $\delta_{n_k}=\delta_{n_0}$, i.e., $n_k=n_0$, for all sufficiently large $k$, and by reading \eqref{kk0conv} backwards we get $\lim_{k}\norm{u_k-u}_{H_{n_0}}=0$, which gives the desired result.
\end{enumerate}
By Definition~\ref{Hilbertconv} and Remark~\ref{2cases} all types of strongly convergent sequences are contained in these two cases, with general sequences $n_k\to\infty$ in Case 1. Conversely, if
$\lim_k d(u_k,u)=0$ for some $u_k\in H_{n_k}$, $u\in H_{n_0}$, then $\delta_{n_k}\to\delta_{n_0}$, so that either $n_0=\infty$ and $n_k\to\infty$, or $n_0\in\N$ and $n_k=n_0$ eventually. Thus we are in one of the two cases above.

\textbf{Step 2.} It remains to prove the existence of the isometric isomorphisms $\Psi_n:H_\infty\to H_n$ with the property \eqref{Psitwoconv}.
To this end, for $k\in\N$ let $\Lscr^k:=\operatorname{lin}\{e_1,\ldots,e_k\}$, $\Lscr_n^k:=\Phi_n(\Lscr^k)$, which
is $k$-dimensional as well since the $\Phi_n$'s are injective. Now by \eqref{Phiuconv} the sequence of symmetric
matrices $\{M_n^k\}\subset\R^k\otimes\R^k$ given by
\[(M_n^k)_{1\le i,j\le k}:=\left((\Phi_n(e_i),\Phi_n(e_j))_{H_n}\right)_{1\le i,j\le k}\]
tends to $\mathbbm{1}_k:=\operatorname{diag}\{1,\ldots,1\}\in\R^k\otimes\R^k$ in the usual matrix norm
\[\norm{(a_{i,j})_{1\le i,j\le k}}_{\R^k\otimes\R^k}:=\max_{1\le i,j\le k}\abs{a_{i,j}},\]
and thus in the standard operator norm $\norm{\;}_{\Lscr_n^k\to\Lscr_n^k}$, i.e.,
\[\norm{M_n^k-\mathbbm{1}_k}_{\Lscr_n^k\to\Lscr_n^k}\le\text{const.}\norm{M_n^k-\mathbbm{1}_k}_{\R^k\otimes\R^k}<\frac{1}{k}\quad\forall n\ge n(k),\]
for some strictly increasing sequence of natural numbers $n(k)$. Here $M_n^k$ is regarded as an operator on $\Lscr_n^k$ via the identification of $\R^k$ with $\Lscr_n^k$ by an orthonormal basis $\{e_1^n,\ldots,e_k^n\}$ of $\Lscr_n^k$ as below.
Here we have used Remark~\ref{phiphi} (1) and Lemma~\ref{basicconv} (4).
Let $n\ge n(k)$. For each such $n$ we
want to define
\[\Psi_n:\Lscr^k\subset C\subset H_\infty\to\Lscr_n^k\subset H_n\]
such that for the operator
\[(\Psi_n-\Phi_n)\r_{\Lscr^k}:\Lscr^k\to\Lscr_n^k\subset H_n\]
we have
\[\norm{(\Psi_n-\Phi_n)\r_{\Lscr^k}}_{\Lscr^k\to H_n}<\frac{1}{k},\]
provided $n\ge n(k)$. This is needed to establish \eqref{Psitwoconv}.

There exists an orthonormal basis $\{e_1^n,\ldots,e_k^n\}$ of $\Lscr_n^k$ and a bijective linear operator $V_n^k:\Lscr_n^k\to\Lscr_n^k$ with
\[\Phi_n(e_i)=\sum_{j=1}^k (V_n^k e_i^n,e_j^n)_{H_n} e_j^n,\]
such that $(V_n^k)^\ast V_n^k=M_n^k$.
More precisely, $V_n^k e_i^n:=\Phi_n(e_i)$. Here $\{e_1^n,\ldots,e_k^n\}$ can be any orthonormal basis of $\Lscr_n^k$ with respect to the inner product of $H_n$. Then the matrix of $(V_n^k)^\ast V_n^k$ with respect to this basis is $M_n^k$.

By the polar decomposition theorem~\cite[Theorem 5.7.1]{LM} we can
represent $V_n^k=B_n^k U_n^k$ as a composition of a positive definite, self-adjoint
operator $B_n^k=\sqrt{V_n^k(V_n^k)^\ast}$ on $\Lscr_n^k$ and an isometric operator $U_n^k$ on $\Lscr_n^k$.
$B_n^k$ is even strictly positive definite, since $V_n^k$ is invertible, hence all eigenvalues $\lambda_{n,i}^k$, $1\le i\le k$, of $B_n^k$
are strictly positive.

For a bounded self-adjoint linear operator $A$ and a polynomial $P$ we have
\[\norm{P(A)}=\sup_{\lambda\in\sigma(A)}\abs{P(\lambda)},\]
where $\sigma(A)$ denotes the spectrum of $A$, see~\cite[Section VII.1 Lemma 2]{ReSi1}.

Furthermore, we have $\mathbbm{1}_k=U_n^k(U_n^k)^\ast$ and
\[(B_n^k)^2=(B_n^k)^\ast B_n^k=(V_n^k(U_n^k)^\ast)^\ast V_n^k(U_n^k)^\ast=U_n^k (V_n^k)^\ast V_n^k (U_n^k)^\ast=U_n^k M_n^k (U_n^k)^\ast.\]
Hence, using also the fact that $U_n^k$ is isometric, for each $n\ge n(k)$ we have
\[\max_{1\le i\le k}\{\abs{(\lambda_{n,i}^k)^2-1}\}=\norm{(B_n^k)^2-\mathbbm{1}_k}_{\Lscr_n^k\to\Lscr_n^k}\le \norm{M_n^k-\mathbbm{1}_k}_{\Lscr_n^k\to\Lscr_n^k}<\frac{1}{k}.\]
Hence $\sqrt{1-\frac{1}{k}}<\lambda_{n,i}^k<\sqrt{1+\frac{1}{k}}$ for all $1\le i\le k$.
Define $\Psi_n^k:\Lscr^k\to\Lscr_n^k$ by
\[\Psi_n^k(e_i):=U_n^k e_i^n\]
and linear extension.
Then $(\Psi_n^k(e_i),\Psi_n^k(e_j))_{H_n}=\delta_{i,j}$, where $\delta_{i,j}$ denotes the Kronecker-delta. Moreover, for $1\le i\le k$
\[(\Phi_n-\Psi_n^k)(e_i)=(V_n^k-U_n^k)(e_i^n)\]
and consequently for $n\ge n(k)$ again by the above remark about spectra and the fact that $U_n^k$ is isometric
\begin{multline*}
\norm{\Phi_n-\Psi_n^k}_{\Lscr^k\to\Lscr_n^k}=\norm{V_n^k-U_n^k}_{\Lscr_n^k\to\Lscr_n^k}
\le\norm{B_n^k-\mathbbm{1}_k}_{\Lscr_n^k\to\Lscr_n^k}\\=\max_{1\le i\le k}\{\abs{\lambda_{n,i}^k-1}\}\le\max_{1\le i\le k}\{\abs{(\lambda_{n,i}^k)^2-1}\}<\frac{1}{k}.
\end{multline*}
To extend $\Psi_n^k$ to $H_\infty$, pick for every $n\in\N$ some arbitrary isometric isomorphisms $\tilde{\Psi}^k_n:(\Lscr^k)^\perp\to(\Lscr_n^k)^\perp$
and extend $\Psi_n^k$ by setting $\Psi_n^k(v):=\tilde{\Psi}_n^k(v)$ for every $v\in(\Lscr^k)^\perp$. Set $\Psi_\infty:=\Id_{H_\infty}:H_\infty\to H_\infty$ and
pick some arbitrary isometric isomorphisms $\Psi_n:H_\infty\to H_n$ for $1\le n< n(1)$. For $n(k)\le n<n(k+1)$, set $\Psi_n:=\Psi_n^k$. Then each $\Psi_n:H_\infty\to H_n$ is an isometric isomorphism and
\[\norm{(\Phi_n-\Psi_n)\r_{\Lscr^{k}}}_{\Lscr^{k}\to H_n}<\frac{1}{k}\]
for all $n\ge n(k)$: indeed, if $n\ge n(k)$, then $n(j)\le n<n(j+1)$ for some $j\ge k$, and since $\Lscr^k\subset\Lscr^j$,
\[\norm{(\Phi_n-\Psi_n^j)\r_{\Lscr^{k}}}_{\Lscr^{k}\to H_n}\le\norm{\Phi_n-\Psi_n^j}_{\Lscr^{j}\to\Lscr_n^{j}}<\frac{1}{j}\le\frac{1}{k}.\]
Evidently, for every fixed $k_0\in\N$ and every $v\in\Lscr^{k_0}$ we have
\[\lim_n\norm{(\Phi_n-\Psi_n)v}_{H_n}=0.\]
Finally, if we can prove \eqref{Psitwoconv}, we are done. So let
$u=\sum_{i=1}^\infty\a_i e_i\in H_\infty$. Let us first prove that $\Psi_n u\to u$ strongly in $\Hscr$. Therefore,
set $\tilde{u}_m:=\sum_{i=1}^m \a_i e_i\in C$. Then $\tilde{u}_m\to u$ in $H_\infty$ as $m\to\infty$ and
\begin{align*}
&\lim_m\varlimsup_n\norm{\Phi_n\tilde{u}_m-\Psi_n u}_{H_n}\\
=&\lim_m\varlimsup_n\lrnorm{(\Phi_n-\Psi_n)\sum_{i=1}^m \a_i e_i-\Psi_n\sum_{i=m+1}^\infty \a_i e_i}_{H_n}\\
\le&\lim_m\varlimsup_n\lrnorm{(\Phi_n-\Psi_n)\sum_{i=1}^m \a_i e_i}_{H_n}+\lim_m\varlimsup_n\lrnorm{\Psi_n\sum_{i=m+1}^\infty \a_i e_i}_{H_n}\\
=&\lim_m\lrnorm{\sum_{i=m+1}^\infty \a_i e_i}_{H_\infty}\\
=&0.
\end{align*}
Now by Lemma~\ref{twoconv}, since $\Psi_n u\to u$ strongly for any $u\in H_\infty$, for $u_n\in H_n$, $n\in\N$, $u\in H_\infty$, we have
$u_n\to u$ $\Hscr$-strongly if and only if
\[\lim_n\norm{\Psi_n u-u_n}_{H_n}=0.\]
Let now $\{n_k\}$ be a general sequence with $n_k\to\infty$ in the sense of Definition~\ref{Hilbertconv}, with $\Psi_\infty=\Id$ and $\Phi_{\infty,\infty}=\Id$. Since $\norm{\Phi_{n}\tilde u_m-\Psi_{n}u}_{H_{n}}\le\norm{(\Phi_n-\Psi_n)\tilde u_m}_{H_n}+\norm{\tilde u_m-u}_{H_\infty}$ for all $n\in\N\cup\{\infty\}$, the above computation gives $\lim_m\varlimsup_k\norm{\Phi_{n_k}\tilde{u}_m-\Psi_{n_k}u}_{H_{n_k}}=0$, i.e., $\Psi_{n_k}u\to u$ strongly along $\{n_k\}$. Since Lemmas~\ref{basicconv} and \ref{twoconv} hold verbatim along $\{n_k\}$, we obtain that $u_k\to u$ strongly if and only if $\lim_k\norm{\Psi_{n_k}u-u_k}_{H_{n_k}}=0$.
In addition, it is clear that the following mapping $\Psi:\Hscr\to I\times \ell^2$, $I:=\{0\}\cup\bigcup_{n\in\N}\{\frac{1}{n}\}$ is the desired
isometry of metric spaces:
\[\Psi(v):=(\delta_n,(v_i)_{i\in\N})\in I\times\ell^2,\]
where $v\in H_n$ and $v_i:=(v,\Psi_n(e_i))_{H_n}$, and where $I\times\ell^2$ carries the metric $((\delta,x),(\delta',x'))\mapsto(\abs{\delta-\delta'}^2+\norm{x-x'}_{\ell^2}^2)^{1/2}$.
It is clearly bijective. Note that $I\times\ell^2$ is complete as a product of complete metric spaces, since $I$ is compact in $\R$ and hence complete. Therefore, $(\Hscr,d_\Hscr)$ is complete as well. Furthermore, $\bigcup_{n\in\N}\operatorname{lin}_\Q\bigcup_{i\in\N}\{\Psi_n e_i\}$ is a countable set which is
dense in $\Hscr$ w.r.t. $d_\Hscr$, where $\operatorname{lin}_\Q$ denotes the \textit{linear span} w.r.t. $\Q$.

The proof is complete.
\end{proof}

\begin{assumption}[Standing assumption]\label{infdim}
From now on, all Hilbert spaces $H_n$, $n\in\N$, and $H=H_\infty$ are assumed to be infinite dimensional and separable.
\end{assumption}
Assumption~\ref{infdim} guarantees that the isometries $\Psi_n$ and the metric $d_\Hscr$ of Theorem~\ref{Kolmetric} exist. They are used in Corollary~\ref{exstrongconv}, Lemma~\ref{dualconv} (1), and hence in Lemmas~\ref{Bdual}, \ref{opnormconv}, \ref{compactB} and Corollary~\ref{Badjoint}, in Lemma~\ref{OpPhiConv}, Theorem~\ref{Kato}, Corollary~\ref{doubletrouble} and in all results whose proofs rely on these statements. The preceding Lemmas~\ref{basicconv} and \ref{twoconv} and Remark~\ref{phiphi} remain valid without Assumption~\ref{infdim}, and so do Lemma~\ref{normmosco} and the ``if'' part of Lemma~\ref{weakmosco} below. All concrete spaces considered in Sections~\ref{sec:mosco}--\ref{sec:stoch} are of the form $L^2(E;\mu)$ with a fully supported measure $\mu$ on $E=\R^d$ or on an infinite dimensional space $E$. Since the support of $\mu$ contains infinitely many disjoint open sets of positive measure, these spaces are infinite dimensional, so that Assumption~\ref{infdim} is satisfied in all applications. The same applies to the domains $\Dscr(\Escr^n)$, $\Dscr(\Escr)$ equipped with the norms $\Escr^n_1$, $\Escr_1$, since in all applications $\Phi_n$ is the identity and these domains contain the set $C$, which is dense in the infinite dimensional spaces $H_n$ and $H$, respectively.

\begin{cor}\label{exstrongconv}
\begin{enumerate}[(i)]
\item For any $u\in H$ there exists a sequence $\{u_n\}$, $u_n\in H_n$, $n\in\N$, such that $u_n\to u$ strongly.
\item If $n_k\uparrow\infty$ strictly as $k\to\infty$ and $v_{k}\in H_{n_k}$, $k\in\N$,
such that $v_{k}\to u\in H$ in $d_\Hscr$, then there exist $u_n\in H_n$, $n\in\N$,
such that $u_n\to u$ in $d_\Hscr$ and $u_{n_k}=v_k$ for every $k\in\N$.
\item If $n_k\uparrow\infty$ strictly as $k\to\infty$ and $v_{k}\in H_{n_k}$, $k\in\N$, converge weakly to $u\in H$ along $\{n_k\}$, then there exist $u_n\in H_n$, $n\in\N$, such that $u_n\to u$ weakly and $u_{n_k}=v_k$ for every $k\in\N$.
\item A sequence $\{u_n\}$, $u_n\in H_n$, $n\in\N$, converges weakly to $u\in H$ if and only if every subsequence of $\{u_n\}$ has a further subsequence which converges weakly to $u$. The same holds for strong convergence.
\end{enumerate}
\end{cor}
\begin{proof}
\begin{enumerate}[(i):]
\item Obvious.
\item For $n\not\in\{n_k\;\vert\;k\in\N\}$ define $u_n:=\Psi_n(u)$ and
$u_{n_k}:=v_k$ for $k\in\N$. Then $\{u_n\}$ is as desired.
\item Define $u_n$ as in (ii). Let $w_n\to w$ strongly. Along $\{n_k\}$, $w_{n_k}\to w$ strongly, hence $(u_{n_k},w_{n_k})_{H_{n_k}}=(v_k,w_{n_k})_{H_{n_k}}\to(u,w)_H$. Along the complementary indices, if there are infinitely many, $(\Psi_n u,w_n)_{H_n}\to(u,w)_H$ by Lemma~\ref{basicconv} (4). Hence $(u_n,w_n)_{H_n}\to(u,w)_H$.
\item We prove the weak case. The strong case is similar and easier. If $u_n\to u$ weakly and $\{n_k\}$ is a subsequence, let $w_k\to w$ strongly along $\{n_k\}$. By (ii), $\{w_k\}$ can be extended to a strongly convergent sequence $\{\tilde w_n\}$ with $\tilde w_{n_k}=w_k$, hence $(u_{n_k},w_k)_{H_{n_k}}\to(u,w)_H$, i.e., $u_{n_k}\to u$ weakly along $\{n_k\}$. Conversely, if $u_n\not\to u$ weakly, there exist $w_n\to w$ strongly, $\eps>0$ and a subsequence $\{n_k\}$ with $\abs{(u_{n_k},w_{n_k})_{H_{n_k}}-(u,w)_H}\ge\eps$ for all $k$. Then no further subsequence of $\{u_{n_k}\}$ converges weakly to $u$, since $w_{n_k}\to w$ strongly along every subsequence.
\end{enumerate}
\end{proof}

The above theorem clarifies the geometric structure of $\Hscr$ as follows. We start with a rather weak limiting structure
of $\Hscr$ based on a ``uniform approximation'' of $H$ by the $H_n$'s via our linear maps $\Phi_n$, which are neither unitary nor even bounded,
on a set of vectors $C$, which may have, for instance, ``nice'' or
``controllable'' properties. The definition of strong and weak convergence seems natural in this setting, but is nevertheless hard to handle and does not
even provide existence of strongly convergent sequences along the $H_n$'s for every possible limit $u\in H$, unless much stronger assumptions on the
$\Phi_n$'s, such as uniform bounds on the operator norms, are made. Finally, by the above theorem, using only basic properties of our concepts and a little
linear algebra, we construct a natural metric on $\Hscr$, and moreover, a sequence of isometric isomorphisms $\{\Psi_n\}$ which in a way carry
over the geometry of $H$ to each $H_n$ and contain the ``asymptotics'' of strong convergence, namely, $u_n\to u$ strongly if and only if $\norm{u_n-\Psi_n(u)}_{H_n}\to 0$. From now on, strong convergence should always be thought
of as given by this characterization: along a ``limit of orthonormal bases'' or ``geometric structures'' via the $\Psi_n$'s. \textbf{From now on}, $\Psi_n$ always denotes the isometries of Theorem~\ref{Kolmetric}.

We note that in applications the $\Phi_n$'s are usually identity operators. Keeping this in mind,
the next two lemmas and Lemma~\ref{OpPhiConv} as well as the later convergence theorems for forms turn out to be very useful, since they reduce the verification of various properties to sequences of the form $\{\Phi_n(\phi)\}$, $\phi\in C$.

\begin{lem}\label{normmosco}
A sequence $\{u_n\}$, $u_n\in H_n$, $n\in\N$ converges to $u\in H$ $\Hscr$-strongly if and only if $\norm{u_n}_{H_n}\to\norm{u}_H$ and
$\deli{u_n}{\Phi_n(\phi)}_{H_n}\to\deli{u}{\phi}_H$ for every $\phi\in C$.
\end{lem}
\begin{proof}
See Appendix A.
\end{proof}

\begin{lem}\label{weakmosco}
Let $\{u_n\}$, $u_n\in H_n$, $n\in\N$ be a sequence in $\Hscr$ and let $u\in H$. Then $u_n\to u$ weakly if and only if $\sup_n\norm{u_n}_{H_n}<\infty$ and
$(u_n,\Phi_n(\phi))_{H_n}\to(u,\phi)_H$ for every $\phi\in C$.
\end{lem}
\begin{proof}
See Appendix A.
\end{proof}

\begin{lem}\label{dualconv}
\begin{enumerate}[(1)]
\item Let $\{u_n\}$ be a sequence with $u_n\in H_n$, $n\in\N$. If the sequence of norms $\{\norm{u_n}_{H_n}\}$ is bounded, there
exists a weakly convergent subsequence of $\{u_n\}$.
\item Let $\{u_n\}$, $u_n\in H_n$, $n\in\N$ be a sequence which weakly converges to $u\in H$. Then
\[\supl_n \norm{u_n}_{H_n}<\infty,\quad \norm{u}_H\le\varliminf_n\norm{u_n}_{H_n}.\]
Moreover, $u_n\to u$ strongly if and only if
\[\norm{u}_H=\lim_n\norm{u_n}_{H_n}.\]
\item A sequence $\{u_n\}$, $u_n\in H_n$, $n\in\N$ tends to $u\in H$ $\Hscr$-strongly if and only if
\[\deli{u_n}{v_n}_{H_n}\to\deli{u}{v}_H\]
for every $\{v_n\}$, $v_n\in H_n$, $n\in\N$ $\Hscr$-weakly tending to $v\in H$.
\end{enumerate}
\end{lem}
\begin{proof}
See Appendix A.
\end{proof}

\subsubsection{Convergence of bounded operators}\label{sss:operators}

\begin{defi}[Convergence of bounded operators]\label{Operatorconvergence}
$\{B_n\}$, $B_n\in \Lscr(H_n)$ are said to \textit{converge strongly} to $B\in \Lscr(H)$ if
for every sequence $\{u_n\}$, $u_n\in H_n$, converging strongly to $u\in H$,
$\{B_n u_n\}$ converges strongly to $Bu$. They are said to \textit{converge weakly} to $B$ if the same holds with strong convergence replaced by weak convergence.

$\{B_n\}$, $B_n\in \Lscr(H_n)$ is said to \textit{compactly converge} to $B\in \Lscr(H)$ if
for every sequence $\{u_n\}$, $u_n\in H_n$ weakly converging to $u\in H$, $\{B_n u_n\}$ strongly converges
to $Bu$.
\end{defi}

Clearly, compact convergence of a sequence of bounded operators implies both weak and strong convergence of this sequence.

\begin{lem}\label{Bdual}
Let $\{B_n\}$ be a sequence of bounded operators, $B_n\in \Lscr(H_n)$, $B\in \Lscr(H)$. Then we have:
\begin{enumerate}[(1)]
\item $B_n\to B$ strongly if and only if
\begin{equation}\label{Bconv}
\lim_n\deli{B_n u_n}{v_n}_{H_n}=\deli{Bu}{v}_H
\end{equation}
for any $\{u_n\}$, $\{v_n\}$, $u$, $v$ such that $u_n\to u$ strongly and $v_n\to v$ weakly.
\item $B_n\to B$ weakly if and only if \eqref{Bconv} holds for any $\{u_n\}$, $\{v_n\}$, $u$, $v$ such
that $u_n\to u$ weakly and $v_n\to v$ strongly.
\item $B_n\to B$ compactly if and only if \eqref{Bconv} holds for any $\{u_n\}$, $\{v_n\}$, $u$, $v$
such that $u_n\to u$ weakly and $v_n\to v$ weakly.
\end{enumerate}
\end{lem}
\begin{proof}
The lemma follows from the definitions of the convergences and Lemma~\ref{dualconv} (3).
\end{proof}

\begin{lem}\label{opnormconv}
\begin{enumerate}[(1)]
\item If $B_n\to B$ strongly, then
\[\varliminf_n\norm{B_n}_{\Lscr(H_n)}\ge\norm{B}_{\Lscr(H)}.\]
\item If $B_n\to B$ compactly, then
\[\lim_n\norm{B_n}_{\Lscr(H_n)}=\norm{B}_{\Lscr(H)}.\]
\end{enumerate}
\end{lem}
\begin{proof}
See Appendix A.
\end{proof}

Denote by $\hat{A}$ the adjoint of an operator $A$. The following is a direct consequence of Lemma~\ref{Bdual}.
\begin{cor}\label{Badjoint}
\begin{enumerate}[(1)]
\item $B_n\to B$ strongly if and only if $\hat{B}_n\to\hat{B}$ weakly. In particular, strong convergence is
equivalent to weak convergence for symmetric operators.
\item $B_n\to B$ compactly if and only if $\hat{B}_n\to\hat{B}$ compactly.
\end{enumerate}
\end{cor}

Note that strong convergence in the sense of Definition~\ref{Operatorconvergence} does not imply weak convergence in this sense, not even if $H_n=H$ and $\Phi_n=\Id$ for all $n$, since weak convergence of operators is tested along weakly convergent sequences of vectors. Thus these notions should not be confused with the strong and weak operator topologies. Also, uniform operator convergence does
not make sense in this framework.

\begin{lem}\label{compactB}
If $B_n\to B$ compactly, then $B$ and $\hat{B}$ are both compact operators.
\end{lem}
\begin{proof}
See Appendix A.
\end{proof}

The next lemma shows that $C$ contains enough information to verify strong convergence of a sequence of operators, if
strong convergence along the sequence $\{\Phi_n(\phi)\}$, $\phi\in C$ is assumed. This observation goes back to~\cite{Toe1}.

\begin{lem}\label{OpPhiConv}
Let $B_n\in\Lscr(H_n)$, $n\in\N$, $B\in\Lscr(H)$, such that $\sup_n\norm{B_n}_{\Lscr(H_n)}<\infty$.
\begin{enumerate}[(1)]
\item If $B_n\Phi_n\phi\to B\phi$
$\Hscr$-strongly for every $\phi\in C$,
\item or if there exists a dense linear subspace $C'\subset H$, which might be taken equal to $H$, such that
$B_n\Psi_n\psi\to B\psi$ $\Hscr$-strongly for every $\psi\in C'$,
\end{enumerate}
then $B_n\to B$ strongly in the sense of Definition~\ref{Operatorconvergence}.
\end{lem}
\begin{proof}
We prove only (2), since the proof of (1) is similar, except that $B_n\Psi_n\psi_m$ in $\eqref{BBB}$ is replaced by $B_n\Phi_n\psi_m$.
Let $u\in H$ and $u_n\in H_n$, $n\in\N$ with $u_n\to u$ $\Hscr$-strongly. Pick $\psi_m\to u$ strongly in $H$ with $\psi_m\in C'$ for every $m$.
Clearly,
\begin{equation}\label{BBB}
\begin{split}
&\norm{B_n u_n-\Psi_n B u}_{H_n}\\
\le&\norm{B_n u_n-B_n \Psi_n \psi_m}_{H_n}+\norm{B_n\Psi_n\psi_m-\Psi_n B\psi_m}_{H_n}+\norm{\Psi_n B\psi_m-\Psi_n Bu}_{H_n}\\
\le&\sup_n\norm{B_n}_{\Lscr(H_n)}\norm{u_n-\Psi_n \psi_m}_{H_n}+\norm{B_n\Psi_n\psi_m-\Psi_n B\psi_m}_{H_n}+\norm{B}_{\Lscr(H)}\norm{\psi_m-u}_{H}.
\end{split}
\end{equation}
After taking $\varlimsup_n$ and then letting $m\to\infty$, the first term tends to $0$ by the uniform boundedness of $\{\norm{B_n}_{\Lscr(H_n)}\}$ and
\[\norm{u_n-\Psi_n\psi_m}_{H_n}\le\norm{u_n-\Psi_n u}_{H_n}+\norm{\Psi_n u-\Psi_n \psi_m}_{H_n}=\norm{u_n-\Psi_n u}_{H_n}+\norm{u-\psi_m}_H.\]
Similarly, the second term tends to $0$ by assumption, and so does the last term.
This proves the assertion.
\end{proof}

\subsubsection{Convergence of semigroups}\label{sss:semigroups}

In this subsection we prove a generalization of Kato's theorem~\cite[Theorem IX.2.16]{Ka} for strong convergence of semigroups, i.e., that
strong convergence of not necessarily symmetric contraction semigroups is equivalent to strong convergence of the associated contraction
resolvents. This result first appeared in the diploma thesis~\cite{Toe1} from 2006. An alternative proof by a different method, which also includes strong graph convergence of the generators, is given in~\cite[Theorem 7.10]{Toe2}. The proof below uses the isometries $\Psi_n$ of Theorem~\ref{Kolmetric}, since the standard proof cannot be transferred one-to-one to the present framework.

To this end, let $(A^n,D(A^n))$, $n\in\N$, and $(A,D(A))$ be the infinitesimal generators of not necessarily symmetric
$C_0$-contraction semigroups $(T_t^n)_{t\ge 0}$, $n\in\N$, and $(T_t)_{t\ge 0}$, respectively, defined on subspaces $D(A^n)\subset H_n$, $n\in\N$, $D(A)\subset H$.
Let $(G^n_\a)_{\a>0}$, $n\in\N$, and $(G_\a)_{\a>0}$ be the associated $C_0$-contraction resolvents.

\begin{thm}\label{Kato}
Let $(T_t^n)_{t\ge 0}$, $n\in\N$, $(T_t)_{t\ge 0}$, $(G^n_\a)_{\a>0}$, $n\in\N$, $(G_\a)_{\a>0}$ be as above.
Then $G_\a^n\to G_\a$ strongly for any $\a>0$ if and only if $T_t^n\to T_t$ strongly for any $t\ge 0$.
\end{thm}
\begin{proof}
Let us first prove the ``if'' part: Let $\a>0$, $u_n\in H_n$, $n\in\N$, $u\in H$, $u_n\to u$ $\Hscr$-strongly.
Then $T_t^n u_n\to T_t u$ strongly for every $t\ge 0$.
We express $G_\a^n u_n$ in terms of
the following $H_n$-valued Bochner integral
\[G_\a^n u_n=\int_0^\infty e^{-\a s}T_s^n u_n ds,\]
see~\cite[Section I.1]{MR} and, for Bochner integrals,~\cite[Section V.5]{Y}.
Let $v_n\in H_n$, $n\in\N$, $v\in H$ such that $v_n\to v$ weakly in $\Hscr$. It suffices to prove that
\[(G_\a^n u_n,v_n)_{H_n}\to(G_\a u,v)_H.\]
We note that
\[\sup_n \abs{e^{-\a s}(T_s^n u_n,v_n)_{H_n}}\le e^{-\a s}\sup_n\norm{u_n}_{H_n}\norm{v_n}_{H_n}\le Ce^{-\a s}\]
for some constant $C>0$ by the
contraction property of the $T_s^n$'s, the Cauchy--Schwarz inequality and strong and weak
convergence of $\{u_n\}$ and $\{v_n\}$, respectively, together with Lemma~\ref{dualconv} (2).
The right-hand side is integrable, so we can use Lebesgue's dominated convergence theorem and 
the well-known fact that Bochner integrals interchange with continuous linear functionals, see~\cite[Corollary V.5.2]{Y}, to obtain
\begin{align*}
\lim_{n\to\infty}(G_\a^n u_n,v_n)_{H_n}&=\lim_{n\to\infty}\lrdel{\int_0^\infty e^{-\a s}T_s^n u_n\;ds}{v_n}_{H_n}\\
&=\lim_{n\to\infty}\int_0^\infty e^{-\a s}(T_s^n u_n,v_n)_{H_n}ds=\int_0^\infty e^{-\a s}(T_s u,v)_H ds\\
&=\lrdel{\int_0^\infty e^{-\a s} T_s u\;ds}{v}_H=(G_\a u,v)_H.
\end{align*}
The ``if'' part is proved.

To prove the ``only if'' part, let $G_\a^n\to G_\a$ strongly for $\a>0$. Recall that for each $\a>0$, $n\in\N$, we have
$G_\a^n(H_n)=D(A^n)$, $G_\a(H)=D(A)$, see~\cite[Section I.1.a)]{MR}. We show for each $\psi\in H$, $t\ge 0$, $\a>0$ that
\begin{equation}\label{TPsiGpsi}
T^n_t\Psi_n G_\a \psi\to T_t G_\a \psi
\end{equation}
strongly in $\Hscr$. Since $D(A)\subset H$ densely, see e.g.~\cite[Proposition I.1.10]{MR}, and since $\sup_n\norm{T_t^n}_{\Lscr(H_n)}\le 1<\infty$, we can
apply Lemma~\ref{OpPhiConv} (2) and we are done.

In the following, fix $\a>0$. Now let $u\in H$ and $s\ge 0$. We have
\[\frac{d}{ds}T_s G_\a u=T_s A G_\a u=-T_s(1-\a G_\a)u\]
and the same with $T_s^n$, $A^n$ and $G_\a^n$ replacing $T_s$, $A$ and
$G_\a$, respectively.
Now for $t\ge s\ge 0$
\begin{align*}
&\frac{d}{ds}T_{t-s}^n G_\a^n\Psi_n T_s G_\a u\\
&=T_{t-s}^n G_\a^n\Psi_n\left(\frac{d}{ds}T_s G_\a\right)+\frac{d}{ds'}T_{t-s'}^n G_\a^n \Psi_n T_s G_\a u\Bigg\vert_{s'=s}\\
&=T_{t-s}^n G_\a^n\Psi_n (-T_s(1-\a G_\a)u)+T_{t-s}^n(1-\a G_\a^n)\Psi_n T_s G_\a u\\
&=-T_{t-s}^n(G_\a^n\Psi_n T_s -\Psi_n T_s G_\a)u,
\end{align*}
where the ``product rule'' is proved the same way as in finite dimensions, with one-sided derivatives at $s=0$ and $s=t$. Here $s\mapsto T_{t-s}^nG_\a^n\Psi_n$ is uniformly bounded and strongly continuous, and $s\mapsto T_{t-s}^nG_\a^n\Psi_nx$ is continuously differentiable for every $x\in H$, since $G_\a^n\Psi_nx\in D(A^n)$. Note that we have used that $\Psi_n$ is continuous and the fact that
$T_s G_\a u=G_\a T_s u$ for $u\in H$, see~\cite[Chapter IX]{Ka}. Applying the
fundamental theorem of calculus we obtain
\begin{multline*}
-\int_0^t T_{t-s}^n(G_\a^n\Psi_n T_s-\Psi_n G_\a T_s)u\;ds=\int_0^t \frac{d}{ds} T_{t-s}^n G_\a^n \Psi_n G_\a T_s u\;ds\\
=T_{t-s}^n G_\a^n \Psi_n G_\a T_s u\Bigg\vert_{s=0}^{s=t}=G_\a^n\Psi_n G_\a T_t u- T_t^n G_\a^n \Psi_n G_\a u.
\end{multline*}
Hence by Bochner's inequality, see~\cite[Corollary V.5.1]{Y},
\begin{equation}\label{gaga}
\norm{G_\a^n\Psi_n G_\a T_t u- T_t^n G_\a^n \Psi_n G_\a u}_{H_n}\le\int_0^t \norm{(G_\a^n\Psi_n T_s-\Psi_n G_\a T_s)u}_{H_n} ds.
\end{equation}
First note that $\norm{G_\a^n\Psi_n T_su-\Psi_n G_\a T_s u}_{H_n}\to 0$ as $n\to\infty$ by Lemma~\ref{twoconv} and the
fact that both $G_\a^n \Psi_n T_s u\to G_\a T_s u$ and $\Psi_n G_\a T_s u\to G_\a T_s u$ $\Hscr$-strongly, which follows
from the strong convergence of resolvents. It is easy to see that
$\norm{(G_\a^n\Psi_n T_s-\Psi_n T_s G_\a)u}_{H_n}\le\frac{2}{\a}\norm{u}_{H}$ for every $n$. We conclude that by Lebesgue's
dominated convergence theorem the right-hand side of \eqref{gaga} tends to zero as $n\to\infty$. Altogether, for $u\in H$,
\[\lim_n\norm{G_\a^n\Psi_n T_t G_\a u- T_t^n G_\a^n \Psi_n G_\a u}_{H_n}=0\]
and since $G_\a(H)=D(A)$, and $D(A)\subset H$ densely, we get by a $3$-$\eps$-argument, taking into account that
$\norm{G_\a^n}_{\Lscr(H_n)}\le\a^{-1}$ for every $n$, that for every $\psi\in H$
\begin{equation}\label{gaga2}
\lim_n\norm{G_\a^n\Psi_n T_t \psi- T_t^n G_\a^n \Psi_n \psi}_{H_n}=0.
\end{equation}
As another result of strong convergence of resolvents and Lemma~\ref{twoconv}, using the contraction property of the
$T^n_t$'s, we get
\begin{equation}\label{gaga1}
\varlimsup_n \norm{T_t^n \Psi_n G_\a \psi-T_t^n G_\a^n \Psi_n \psi}_{H_n}\le\lim_n\norm{\Psi_n G_\a \psi-G_\a^n \Psi_n \psi}_{H_n}=0,
\end{equation}
and
\begin{equation}\label{gaga3}
\lim_n \norm{G_\a^n \Psi_n T_t \psi-\Psi_n G_\a T_t \psi}_{H_n}=0.
\end{equation}
To prove \eqref{TPsiGpsi}, let $\psi\in H$, $t\ge 0$, $\a >0$. One easily observes
\begin{align*}
&\norm{T_t^n\Psi_n G_\a \psi-\Psi_n T_t G_\a \psi}_{H_n}\\
\le&\norm{T_t^n\Psi_n G_\a \psi- T_t^n G_\a^n \Psi_n\psi}_{H_n}+\norm{T_t^n G_\a^n \Psi_n\psi-G_\a^n\Psi_n T_t \psi}_{H_n}\\
&+\norm{G_\a^n\Psi_n T_t\psi-\Psi_n G_\a T_t \psi}_{H_n}
\end{align*}
Combining \eqref{gaga1}, \eqref{gaga2} and \eqref{gaga3} and the last equation above, we get
the desired result as $n\to\infty$. The proof is complete.
\end{proof}

\subsection{Convergence of symmetric forms}\label{ss:symmetric}

Now we consider three different types of convergence of a sequence of symmetric forms $\{\Escr_n\}$ along $\Hscr$, where each form is
defined on $H_n$.
Recall that a bilinear form is a bilinear mapping $\Escr:\Dscr(\Escr)\times\Dscr(\Escr)\to\R$ defined
on some subspace $\Dscr(\Escr)\subset H$. In this subsection we consider only \textit{non-negative} and \textit{symmetric} quadratic forms, that
is, $\Escr(u,u)\ge 0$ for every $u\in\Dscr(\Escr)$ and $\Escr(u,v)=\Escr(v,u)$ for every $u,v\in\Dscr(\Escr)$.
Define for $\alpha>0$ and a form $(\Escr,\Dscr(\Escr))$ the inner product
\[\Escr_\a(u,v):=\Escr(u,v)+\a\deli{u}{v}_H,\quad u,v\in\Dscr(\Escr),\]
which makes $\Dscr(\Escr)$ a pre-Hilbert space.
Recall that a form $\Escr$ is \textit{closed}
if $\Dscr(\Escr)$ equipped with the norm $\Escr_1^{1/2}$ is
complete. We identify a quadratic form $\Escr$ with the functional
\begin{equation}\label{DEinfty}
\Escr(u):u\mapsto\begin{cases}\Escr(u,u)&,\; u\in\Dscr(\Escr)\\\infty&,\; u\not\in\Dscr(\Escr).\end{cases}
\end{equation}
It is well known that $\Escr$ is closed if and only if $\Escr:H\to\ol{\R}$ is lower semicontinuous, see for instance~\cite[p. 372]{M}.
We shall use the notions quadratic form and bilinear form interchangeably if the form is symmetric and non-negative, which
is justified by the polarization identity:
\[\Escr(u,v)=\frac{1}{4}\left[\Escr(u+v,u+v)-\Escr(u-v,u-v)\right],\quad u,v\in\Dscr(\Escr).\]
Also note that $\Escr(u):=\Escr(u,u)$, $u\in\Dscr(\Escr)$ is used for the diagonal even for non-symmetric forms.

\begin{defi}\label{Mosco}
A sequence $\{\Escr^n:H_n\to\ol{\R}\}$ of symmetric, non-negative, closed forms is said to \textit{Mosco converge}
to a quadratic form $\Escr$ on $H$ if the following two conditions hold:
\begin{enumerate}
\item[(M1)] If $\{u_n\}$, $u_n\in H_n$, $n\in\N$ weakly converges to $u\in H$ then
\[\Escr(u)\le\varliminf\limits_n \Escr^n(u_n).\]
\item[(M2)] For every $u\in H$ there exists a strongly convergent sequence $u_n\to u$, $u_n\in H_n$, $n\in\N$
such that
\[\Escr(u)=\lim\limits_n \Escr^n(u_n).\]
\end{enumerate}
\end{defi}
\begin{defi}
A sequence $\{\Escr^n:H_n\to\ol{\R}\}$ of symmetric, non-negative, closed forms is said to \textit{$\Gamma$-converge}
to a quadratic form $\Escr$ on $H$ if the following conditions are fulfilled:
\begin{enumerate}
\item[(G1)] If $\{u_n\}$, $u_n\in H_n$, $n\in\N$ strongly converges to $u\in H$ then
\[\Escr(u)\le\varliminf\limits_n \Escr^n(u_n).\]
\item[(G2)] For every $u\in H$ there exists a strongly convergent sequence $u_n\to u$ with $u_n\in H_n$, $n\in\N$ such that
\[\Escr(u)=\lim\limits_n \Escr^n(u_n).\]
\end{enumerate}
\end{defi}
The above conditions, especially (M2) and (G2), make sense by Corollary~\ref{exstrongconv} and \eqref{DEinfty}.
It is clear that Mosco convergence implies $\Gamma$-convergence. To determine when they are equivalent, consider the following:
\begin{defi}
A sequence $\{\Escr^n\}$ is called \textit{asymptotically compact} if for every strictly increasing sequence $\{n_k\}$ of natural numbers and every
$\{u_k\}$, $u_k\in H_{n_k}$, $k\in\N$ such that
\[\varlimsup_k(\Escr^{n_k}(u_k)+\lrnorm{u_k}_{H_{n_k}}^2)<\infty,\]
there exists a strongly convergent subsequence of $\{u_k\}$.
\end{defi}
\begin{lem}\label{gammamosco}
Assume that $\{\Escr^n\}$ is asymptotically compact. Then $\{\Escr^n\}$ $\Gamma$-converges to $\Escr$ if and only if
$\{\Escr^n\}$ Mosco converges to $\Escr$.
\end{lem}
\begin{proof}
The proof of~\cite[Lemma 2.3.2]{M} can easily be extended to our framework. Note that (G1) and (M1) also hold along subsequences. Indeed, if $u_k\to u$ strongly along $n_k\uparrow\infty$, extend $\{u_k\}$ by Corollary~\ref{exstrongconv} (ii) to a strongly convergent sequence $\{w_n\}$ with $w_{n_k}=u_k$. Then $\Escr(u)\le\varliminf_n\Escr^n(w_n)\le\varliminf_k\Escr^{n_k}(u_k)$. For weakly convergent sequences, the same argument applies with Corollary~\ref{exstrongconv} (iii).
\end{proof}
Note that $\Gamma$-convergence can be defined for arbitrary functionals on a topological space with values in $\ol{\R}$, see for instance~\cite{DaMa},
and every $\Gamma$-limit is lower semicontinuous. In particular, if the form $\Escr$ is a $\Gamma$-limit, then it is closed.

It is a well-known fact that every sequence of functionals on a second-countable space with values in $\ol{\R}$ has a $\Gamma$-convergent
subsequence, see~\cite[Theorem 8.5]{DaMa}. So in our case we have the following result from~\cite[Theorem 2.3]{KS}:
\begin{thm}\label{subgamma}
Every sequence $\{\Escr^n\}$ of symmetric, non-negative, closed quadratic forms with values in $\ol{\R}$
has a $\Gamma$-convergent subsequence whose $\Gamma$-limit is a symmetric, non-negative, closed quadratic form on $H$.
\end{thm}
This just means that the space of symmetric, non-negative forms is relatively sequentially compact
w.r.t. the $\Gamma$-topology.
\begin{defi}
We say that $\Escr^n\to\Escr$ \textit{compactly} if $\Escr^n\to\Escr$ Mosco and if $\{\Escr^n\}$ is asymptotically compact.
\end{defi}
Lemma~\ref{gammamosco} and Theorem~\ref{subgamma} together imply:
\begin{cor}
If $\{\Escr^n\}$ is asymptotically compact, it has a compact convergent subsequence.
\end{cor}

With every non-negative symmetric closed form $\Escr$ we associate a non-negative self-adjoint operator
$-A$ with $\Dscr(A)\subset\Dscr(\sqrt{-A})=\Dscr(\Escr)$ and $\Escr(u,v)=\deli{\sqrt{-A}u}{\sqrt{-A}v}$,
$u,v\in\Dscr(\Escr)$. Let $T_t:=e^{tA}$, $t\ge 0$ be the associated semigroup and
$G_\a:=(\a-A)^{-1}$, $\a>0$ the associated resolvent, see~\cite{FOT} for details. If $\Dscr(\Escr)$ is not dense in $H$, which may happen for a Mosco limit by Theorem~\ref{subgamma}, then $-A$ is a self-adjoint operator on the closed subspace $\ol{\Dscr(\Escr)}$, and $T_t$ and $G_\a$ are extended to $H$ by setting them equal to $0$ on the orthogonal complement of $\ol{\Dscr(\Escr)}$, cf.~\cite{KS}.

For further references on convergence problems of spectral structures in the symmetric case, we again refer to~\cite[Section 2]{KS}.
The next theorem, whose proof has been reformulated for this paper and can be found in the Appendix, shows the strength of Mosco
convergence, which is necessary and sufficient for strong convergence of resolvents and semigroups.

\begin{thm}[Mosco, Kuwae, Shioya]\label{KSConv}
Let $\{\Escr^n:H_n\to\ol{\R}\}$ be a sequence of non-negative, symmetric, closed forms and let $\Escr$ be
a non-negative, symmetric, closed form on $H$. The following statements are equivalent:
\begin{enumerate}[(1)]
\item $\{\Escr^n\}$ Mosco converges to $\Escr$,
\item $\{G_{\a}^n\}$ strongly converges to $G_\a$ for all $\a>0$,
\item $\{T_{t}^n\}$ strongly converges to $T_t$ for all $t\ge 0$.
\end{enumerate}
\end{thm}
\begin{proof}
See Appendix A.
\end{proof}

Now we shall state some useful conditions on a sequence of forms $\{\Escr^n:H_n\to\ol{\R}\}$, which are easy to check in many applications and
provide a convenient criterion for Mosco convergence. The idea is to restrict strong and weak convergence of vectors
along the sequence of domains $\{\Dscr(\Escr^n)\}$, which contains enough information to verify Mosco convergence.
\begin{defi}\label{Dconv}
Suppose that we are given a convergent sequence of Hilbert spaces $H_n\to H$ and a sequence of non-negative symmetric closed forms
$\{\Escr^n:H_n\to\ol{\R}\}$.
We say that \textit{a sequence of pairs $\{(H_n,\Escr^n)\}$ converges to $(H,\Escr)$} if the following conditions hold:
\begin{enumerate}[(1)]
\item $\Phi_n(C)\subset\Dscr(\Escr^n)$ for every $n\in\N$.
\item $C\subset\Dscr(\Escr)$ is dense in $(\Dscr(\Escr),\Escr_1^{1/2})$.
\item $\lim_n \Escr^n (\Phi_n(\phi))=\Escr(\phi)$ for every $\phi\in C$.
\end{enumerate}
\end{defi}
If we have a convergent sequence $\{(H_n,\Escr^n)\}$ of pairs in the above sense, we have a convergent sequence of Hilbert spaces
$\{(\Dscr(\Escr^n),\Escr_1^n(\;,\;))\}$ which converges to $(\Dscr(\Escr),\Escr_1(\;,\;))$ in the sense of Definition~\ref{Hilbertconv}.
Let us denote the space $\disjoint_n \Dscr(\Escr^n)\dot\cup\Dscr(\Escr)$
by $\Hscr_\Escr$. Hence we can construct the corresponding isometries $\Psi_n^\Escr$, $n\in\N$
as in Theorem~\ref{Kolmetric}. Note that $(\Dscr(\Escr^n),\Escr^n_1)$ is separable, since $u\mapsto(\sqrt{-A_n}u,u)$ is an isometry from $(\Dscr(\Escr^n),\Escr^n_1)$ into $H_n\times H_n$, where $A_n$ denotes the generator associated with $\Escr^n$ on $\ol{\Dscr(\Escr^n)}$. Then for every $u\in\Dscr(\Escr)$, $\Psi_n^\Escr(u)\to u$ strongly in $\Hscr_\Escr$.
\begin{lem}\label{strongstrong}
Assume that $\{(H_n,\Escr^n)\}$ converges to $(H,\Escr)$. If $u_n\to u$ strongly in $\Hscr_\Escr$, then $u_n\to u$ strongly
in $\Hscr$.
\end{lem}
\begin{proof}
Let $u_n\to u$, $u_n\in\Dscr(\Escr^n)$, $n\in\N$, $u\in\Dscr(\Escr)$ strongly in $\Hscr_\Escr$. Then, by the definition of strong convergence in $\Hscr_\Escr$, there exists a sequence
$\{\phi_m\}\subset C$ with $\lim_m\varlimsup_n\Escr^n_1(u_n-\Phi_n(\phi_m))=0$ and
\[\norm{u-\phi_m}_H^2\le\textup{const.}\Escr_1(u-\phi_m)\to 0\]
as $m\to\infty$. Clearly,
\[\lim_m\varlimsup_n\norm{u_n-\Phi_n(\phi_m)}_{H_n}^2
\le\textup{const.}\lim_m\varlimsup_n\Escr^n_1(u_n-\Phi_n(\phi_m))=0,\]
hence $u_n\to u$ strongly in $\Hscr$. The case $\Dscr(\Escr^{n_k})\to\Dscr(\Escr^{n_0})$ for some $n_0\in\N$, i.e., $n_k=n_0$ eventually, is trivial.
\end{proof}
For weak convergence, we have the following result.
\begin{prop}\label{weakweakmosco}
Assume that $\{(H_n,\Escr^n)\}$ converges to $(H,\Escr)$.
If for any sequence
$u_n\to u$ weakly in $\Hscr_\Escr$, $u_n\in\Dscr(\Escr^n)$, $n\in\N$,
$u\in\Dscr(\Escr)$ we have that $u_n\to u$ weakly in $\Hscr$, too, then
$\Escr^n\to\Escr$ Mosco.
\end{prop}
For the proof, we need the following lemma.
\begin{lem}\label{densemosco}
Assume that $\{(H_n,\Escr^n)\}$ converges to $(H,\Escr)$ and (M1) holds. Then $\Escr^n\to\Escr$ Mosco, i.e., (M1) implies (M2) in this case.
\end{lem}
\begin{proof}
Let us take $u\in \Dscr(\Escr)$. Clearly, $\Psi_n^\Escr(u)\to u$ $\Hscr_\Escr$-strongly and also $\Hscr$-strongly by Lemma~\ref{strongstrong}. Since $\Psi_n^\Escr$ is an isometry, we have
\[\Escr_1^n(\Psi_n^\Escr(u))=\Escr_1(u)\quad\text{for every }n\in\N,\]
and, by Lemma~\ref{basicconv} (3),
\[\lim_n\norm{\Psi_n^\Escr(u)}_{H_n}^2=\norm{u}_H^2.\]
Therefore
\[\lim_n\Escr^n(\Psi_n^\Escr(u))=\Escr(u)\]
which gives us (M2) for $u\in\Dscr(\Escr)$. If $u\in H\setminus\Dscr(\Escr)$,
then (M2) holds for any $\Hscr$-strongly convergent sequence $u_n\to u$ as a consequence of (M1).
\end{proof}

\begin{proof}[Proof of Proposition~\ref{weakweakmosco}]
We only have to prove (M1) by the preceding lemma.
Let $u_n\in H_n$, $n\in\N$, $u\in H$ with $u_n\to u$ $\Hscr$-weakly.
If $\varliminf_n\Escr^n(u_n,u_n)=+\infty$ there is nothing to prove. So assume
that $\varliminf_n\Escr^n(u_n,u_n)<+\infty$, which implies that for a subsequence
$\{u_{n_k}\}$ we have $\sup_k\Escr^{n_k}(u_{n_k})<+\infty$.
But by Lemma~\ref{dualconv} (2) $\sup_n\norm{u_n}_{H_n}<+\infty$. By extracting
another subsequence if necessary we can use Lemma~\ref{dualconv} (1) to ensure
that for some $\tilde{u}\in\Dscr(\Escr)$ we have $u_{n_k}\to\tilde{u}$ $\Hscr_\Escr$-weakly and
$\lim_k\Escr^{n_k}(u_{n_k})=\varliminf_n\Escr^n(u_n)$.
By Corollary~\ref{exstrongconv} (iii), applied in $\Hscr_\Escr$, the sequence $\{u_{n_k}\}$ can be extended to a sequence converging weakly in $\Hscr_\Escr$ to $\tilde u$. By assumption, this sequence converges $\Hscr$-weakly to $\tilde u$, so that, by Corollary~\ref{exstrongconv} (iv), $u_{n_k}\to\tilde{u}$ $\Hscr$-weakly, which shows
$\tilde{u}=u$. If $\Escr(u,u)=0$, there is nothing to prove. Otherwise,
\[\lim_k\Escr^{n_k}_1(u_{n_k},\Psi_{n_k}^\Escr(u))=\Escr_1(u,u)\]
and by Lemma~\ref{strongstrong}
\[\lim_k (u_{n_k},\Psi_{n_k}^\Escr(u))_{H_{n_k}}=(u,u)_H.\]
Combining the last two equalities, we get
\[\lim_k\Escr^{n_k}(u_{n_k},\Psi_{n_k}^\Escr(u))=\Escr(u,u).\]
Therefore,
\begin{multline*}\Escr(u,u)^2=\lim_k\abs{\Escr^{n_k}(u_{n_k},\Psi_{n_k}^\Escr(u))}^2\\
\le\lim_k\Escr^{n_k}(u_{n_k},u_{n_k})\Escr^{n_k}(\Psi_{n_k}^\Escr(u),\Psi_{n_k}^\Escr(u))=
\varliminf_n\Escr^n(u_n,u_n)\Escr(u,u),
\end{multline*}
which gives us (M1) after division by $\Escr(u,u)>0$.
\end{proof}

\subsection{Convergence of non-symmetric forms}\label{ss:nonsym}

\subsubsection{Generalized forms}\label{sss:genforms}

To analyze convergence of a sequence of non-symmetric forms defined on different Hilbert spaces, we
will assume our forms to be generalized Dirichlet forms, following the framework of~\cite[Section I]{S1}.
In other words, we assume our form $\Escr^n$ to be associated with some coercive closed form $(\Ascr^n,\Vscr_n)$ and
some suitable linear operator
$(\Lambda^n,\Dscr(\Lambda^n,H_n))$. We prove necessary and sufficient conditions
on the forms $\{\Escr^n\}$ for strong convergence of the associated resolvents and semigroups, and we discuss the relation to the Mosco convergence of the symmetric parts $\{\tilde{\Ascr}^n\}$ of $\{\Escr^n\}$ in Proposition~\ref{moscotohino}.

Now let $\Ascr$ be a bilinear form on $H$ with a domain $\Vscr\subset H$. The \textit{symmetric part} $\tilde{\Ascr}$ of $\Ascr$ is defined by
\[\tilde{\Ascr}(u,v):=\frac{1}{2}\left[\Ascr(u,v)+\Ascr(v,u)\right],\quad u,v\in\Vscr.\]
The \textit{antisymmetric part} $\check{\Ascr}$ of $\Ascr$ is defined by
\[\check{\Ascr}(u,v):=\frac{1}{2}\left[\Ascr(u,v)-\Ascr(v,u)\right],\quad u,v\in\Vscr.\]
It is clear that $\Ascr=\tilde{\Ascr}+\check{\Ascr}$.
For $\a >0$, set
\[\Ascr_\a(u,v):=\Ascr(u,v)+\a\deli{u}{v}_H,\quad u,v\in \Vscr.\]
$\tilde{\Ascr}_\a$ is defined similarly. We suppose that $(\Ascr,\Vscr)$ is a \textit{coercive closed form} with sector constant $K\ge 1$, that is,
\begin{enumerate}[(1)]
\item $(\tilde{\Ascr},\Vscr)$ is a non-negative definite, symmetric, closed form with $\Vscr$ dense in $H$,
\item $(\Ascr,\Vscr)$ satisfies the \textit{weak sector condition}, i.e., there exists a sector constant $K\ge 1$ such that
\[\lrabs{\Ascr_1(u,v)}\le K\Ascr_1(u,u)^{1/2}\Ascr_1(v,v)^{1/2}\quad\text{for all}\;u,v\in\Vscr.\]
\end{enumerate}
Equipped with the norm $\norm{\;}_\Vscr:=\tilde{\Ascr}_1^{1/2}(\;)$, $\Vscr$ becomes a Hilbert space. Identifying $H$ with its
dual $H'$ we obtain a dense and continuous embedding $\Vscr\subset H\equiv H'\subset \Vscr'$. The pairing between $\Vscr$ and $\Vscr'$ is
expressed by $\dualdel{\Vscr}{\;}{\;}$.

Let $\Lambda$ be a linear operator on $H$ with a linear domain $\Dscr(\Lambda,H)$. We assume the following:
\begin{enumerate}[(1)]
\item $\Lambda$ generates a $C_0$-semigroup of contractions $(U_t)_{t\ge 0}$ on $H$.
\item $(U_t)_{t\ge 0}$ can be restricted to a $C_0$-semigroup of contractions on $\Vscr$.
\end{enumerate}
Denote the infinitesimal generator of the restricted semigroup by
$(\Lambda,\Dscr(\Lambda,\Vscr))$. Note that the adjoint
operator $(\hat{\Lambda},\Dscr(\hat{\Lambda},\Vscr'))$ of $(\Lambda,\Dscr(\Lambda,\Vscr))$ also satisfies the conditions above.
Since $\Dscr(\Lambda,\Vscr)\subset\Dscr(\Lambda,H)\cap\Vscr$ and $\Dscr(\Lambda,\Vscr)$ is dense in $\Vscr$ as the domain of the generator of a $C_0$-semigroup on $\Vscr$, $\Dscr(\Lambda,H)\cap\Vscr$ is dense in $\Vscr$. It follows from~\cite[Lemma I.2.3]{S1} that $\Lambda:\Dscr(\Lambda,H)\cap\Vscr\to\Vscr'$
is closable. Let us denote its closure by $(\Lambda,\Fscr)$. Then $\Fscr$ is a Hilbert space with norm
\[\norm{\cdot}_\Fscr:=\Bigg(\norm{\cdot}_\Vscr^2+\norm{\Lambda\cdot}_{\Vscr'}^2\Bigg)^{1/2}.\]
Furthermore, define $\hat{\Fscr}:=\Dscr(\hat{\Lambda},\Vscr')\cap\Vscr$ with norm
\[\norm{\cdot}_{\hat{\Fscr}}:=\Bigg(\norm{\cdot}_\Vscr^2+\norm{\hat{\Lambda}\cdot}_{\Vscr'}^2\Bigg)^{1/2}.\]
$\Fscr$ and $\hat{\Fscr}$ are dense in $\Vscr$, $\dualdel{\Vscr}{\Lambda u}{u}\le 0$ for $u\in\Fscr$, $\dualdel{\Vscr}{\hat{\Lambda}u}{u}\le 0$
for $u\in\hat{\Fscr}$ and $\Dscr(\Lambda,\Vscr)$ is dense in $\Fscr$, see~\cite[Lemma I.2.5 and Lemma I.2.6]{S1}.

Now for given $\Ascr$ and $\Lambda$, define the \textit{bilinear form} $\Escr$ \textit{associated with} $(\Ascr,\Vscr)$ \textit{and}
$(\Lambda,\Dscr(\Lambda,H))$ on $H$ by
\[\Escr(u,v):=\begin{cases}\Ascr(u,v)-\dualdel{\Vscr}{\Lambda u}{v},&\;\text{if}\;u\in\Fscr,v\in\Vscr,\\
\Ascr(u,v)-\dualdel{\Vscr}{\hat{\Lambda}v}{u},&\;\text{if}\;u\in\Vscr,v\in\hat{\Fscr}.\end{cases}\]
We extend $\Escr$ to a form defined on $H$ and taking values in $\ol{\R}$
by setting $\Escr(u,v)=+\infty$ for every other case, even if $u\in H\setminus\Vscr$ and $v=0$.

We also define the \textit{co-form} $\hat{\Escr}$ by
\[\hat{\Escr}(u,v):=\begin{cases}\Ascr(v,u)-\dualdel{\Vscr}{\hat{\Lambda} u}{v},&\;\text{if}\;u\in\hat{\Fscr},v\in\Vscr,\\
\Ascr(v,u)-\dualdel{\Vscr}{\Lambda v}{u},&\;\text{if}\;u\in\Vscr,v\in \Fscr.\end{cases}\]

\begin{rem}
Let $(\Ascr,\Vscr)$ be a coercive closed form and $\Lambda=0$. Clearly $\Fscr=\Vscr=\hat{\Fscr}$ and $\Escr=\Ascr$ is a generalized form
by \textup{\cite[Example I.4.9 (i)]{S1}}. This is the case
of coercive closed forms, in particular of Dirichlet forms as described in \textup{\cite{MR}}.
\end{rem}

Let us recall some useful facts.

As usual, we define for $\a>0$ 
\[\Escr_\a(u,v):=\Escr(u,v)+\a\deli{u}{v}_H.\]
\begin{prop}\label{Resolvent}
For all $\a>0$ there exist continuous, linear bijections $W_\a:\Vscr'\to\Fscr$
and $\hat{W}_\a:\Vscr'\to\hat{\Fscr}$ such
that
\[\Escr_\a(W_\a f,v)=\dualdel{\Vscr}{f}{v}=\Escr_\a(v,\hat{W}_\a f)\]
for all $f\in\Vscr'$, $v\in\Vscr$. $(W_\a)_{\a>0}$ and $(\hat{W}_\a)_{\a>0}$
satisfy the resolvent equation, see~\cite[Proposition I.3.4]{S1}.

Furthermore, there exists a unique $C_0$-resolvent $(G_\a)_{\a>0}$ and a unique $C_0$-co-resolvent $(\hat{G}_\a)_{\a>0}$ on $H$, which are the restrictions
of $W_\a$ and $\hat{W}_\a$, respectively, to $H$,
such that for all $\a>0$, $f\in H$ and $u\in\Vscr$
\begin{equation}\label{GML}\begin{split}
&G_\a(H)\subset\Fscr,\quad\hat{G}_\a(H)\subset\hat{\Fscr},\\
&\Escr_\a(G_\a f,u)=\Escr_\a(u,\hat{G}_\a f)=\deli{f}{u}_H.
\end{split}\end{equation}
$\hat{G}_\a$ is the adjoint of $G_\a$ and $\a G_\a$, $\a\hat{G}_\a$ are contraction operators. Also, we have for $u\in\Vscr$ that
\[\lim_{\a\to\infty}\a G_\a u=u\]
strongly in $\Vscr$ and thus in $H$.
\end{prop}
\begin{proof}
See~\cite[Section I.3]{S1}.
\end{proof}

Note that the second line of \eqref{GML} is equivalent to
\[\dualdel{\Vscr}{(M_\a-\Lambda)G_\a f}{g}=(f,g)_H=\dualdel{\Vscr}{(\hat{M}_\a-\hat{\Lambda})\hat{G}_\a f}{g},\quad f\in H,\;g\in\Vscr,
\]
where for $\a>0$ we set $M_\a:\Vscr\to\Vscr'$, $\dualdel{\Vscr}{M_\a u}{\cdot}:=\Ascr_\a(u,\cdot)$ and
$\hat{M}_\a:\Vscr\to\Vscr'$, $\dualdel{\Vscr}{\hat{M}_\a u}{\cdot}:=\Ascr_\a(\cdot,u)$.

Define approximate forms $\Escr^{(\b)}$, $\b>0$ of $\Escr$ by
\[\Escr^{(\b)}(u,v)=\b\deli{u-\b G_\b u}{v}_H,\quad u,v\in H\]
and set $\Escr_\a^{(\b)}(u,v)=\Escr^{(\b)}(u,v)+\a\deli{u}{v}_H$.
\begin{prop}\label{BetaMR}
\begin{enumerate}[(i)]
\item $\Escr^{(\b)}(u,v)=\Escr(\b G_\b u,v)$ for $u\in H$, $v\in\Vscr$.
\item $\Escr^{(\b)}(u,u)=\Escr(\b G_\b u,\b G_\b u)+\b\norm{u-\b G_\b u}_H^2$ for $u\in H$.
\item $\lim_{\b\to\infty}\Escr^{(\b)}(u,v)=\Escr(u,v)$ for $u\in\Fscr$, $v\in\Vscr$.
\item If $\sup_{\b >0}\Escr_1^{(\b)}(u,u)<\infty$, then $u\in\Vscr$.
\end{enumerate}
\end{prop}
\begin{proof}
For (i)--(iii), see~\cite[Lemma I.2.11]{MR} and \cite[Proposition 2.7 (iii)]{S2}.

(iv): Since $\Escr(v,v)\ge\Ascr(v,v)$ for $v\in\Fscr$ and $\b G_\b$ is contractive, we have by (ii)
\begin{multline*}
\Escr_1^{(\b)}(u,u)=\Escr(\b G_\b u,\b G_\b u)+\b\norm{u-\b G_\b u}_H^2+\norm{u}_H^2\\
\ge\Ascr_1(\b G_\b u,\b G_\b u)+\b\norm{u-\b G_\b u}_H^2.
\end{multline*}
Hence the assumption $\sup_{\b > 0}\Escr_1^{(\b)}(u,u)<\infty$ implies that
\begin{align}
\label{beta1}&\sup_{\b >0}\Ascr_1(\b G_\b u,\b G_\b u)<\infty,\\
\label{beta2}&\sup_{\b >0}\b\norm{u-\b G_\b u}_H^2 <\infty.
\end{align}
From \eqref{beta2}, $\b G_\b u\to u$ in $H$ as $\b\to\infty$. Combining this and \eqref{beta1}, we
have that $u\in\Vscr$ by~\cite[Lemma I.2.12]{MR}.
\end{proof}

Let $(T_t)_{t\ge 0}$ and $(\hat{T}_t)_{t\ge 0}$ be the $C_0$-semigroup of contractions and the $C_0$-co-semigroup of contractions
associated with $(G_\a)_{\a >0}$ and $(\hat{G}_\a)_{\a>0}$, respectively.

\subsubsection{Criteria of convergence}\label{sss:criteria}

In this subsection we shall give necessary and sufficient conditions on a sequence of generalized forms $\{\Escr^n\}$ for the strong convergence
of the associated resolvents $\{G_\a^n\}$ and equivalently the weak convergence of co-resolvents $\{\hat{G}_\a^n\}$, see Corollary~\ref{Badjoint}.
We point out that, in a natural
way, we have to introduce a condition on the asymmetry of the $\{\Escr^n\}$, which differs from (M1) and (M2), since clearly for symmetric forms
all information is contained in the diagonal. On the other hand, much stronger assumptions are needed to obtain Mosco convergence of
the symmetric parts $\{\tilde{\Ascr}^n\}$.

In order to handle double indexed sequences we need the following lemma, which is a convenient formulation of the standard diagonal argument.
\begin{lem}\label{metricdoubletrouble}
Let $X$ be a metrizable space with some metric $d$ and $\{x_{n,m}\;\vert\;n,m\in\N\}$ a double indexed sequence in $X$, $\{x_m\}_{m\in\N}\subset X$, $x\in X$
such that
\[d-\liml_{n\to+\infty}x_{n,m}=x_m\]
and
\[d-\liml_{m\to+\infty}x_m=x.\]
Then there exists a nondecreasing mapping $n\mapsto m(n)$ tending to $+\infty$ so that
\[d-\liml_{n\to\infty}x_{n,m(n)}=x.\]
\end{lem}
\begin{proof}
See Appendix A.
\end{proof}

\begin{cor}\label{doubletrouble}
Suppose that double sequences $\{u_{i,j}\}_{i,j\in\N}\subset\Hscr$, $\{a_{i,j}\}_{i,j\in\N}\subset\R$ and $u\in H$, $a\in\R$ satisfy
\begin{align*}
&d_\Hscr-\liml_{j\to\infty}d_\Hscr-\liml_{i\to\infty}u_{i,j}=u,\\
&\liml_{j\to\infty}\liml_{i\to\infty}a_{i,j}=a.
\end{align*}
Then there exists a nondecreasing mapping $i\mapsto j(i)$, tending to $+\infty$, so that
\begin{align*}
&d_\Hscr-\liml_{i\to\infty}u_{i,j(i)}=u,\\
&\liml_{i\to\infty}a_{i,j(i)}=a.
\end{align*}
\end{cor}
\begin{proof}
Apply Lemma~\ref{metricdoubletrouble} to the product space $\Hscr\times\R$ with
the \textit{$1$-product metric} $d_1((h_1,a_1),(h_2,a_2)):=d_\Hscr(h_1,h_2)+\abs{a_1-a_2}$, $h_1,h_2\in\Hscr$, $a_1,a_2\in\R$.
\end{proof}
\noindent\emph{Remark.} Corollary~\ref{doubletrouble} remains valid if $a_{i,j}\in\ol{\R}$ and, for every $j$, $a_{i,j}\in\R$ for all $i\ge N_j$ with some $N_j\in\N$. Indeed, in the proof of Lemma~\ref{metricdoubletrouble} in Appendix A, the thresholds $n_p$ can be chosen such that $n_p\ge N_{m_p}$. Then $a_{i,j(i)}\in\R$ for all $i\ge n_1$, and the assertion follows by applying Corollary~\ref{doubletrouble} to the real numbers $a_{i,j}\mathbbm{1}_{\{i\ge N_j\}}$.

We now turn to the criteria of convergence. In the symmetric case, compactness arguments are based on bounds of the energy $\Escr(u,u)$. In the non-symmetric case, the energy only controls the symmetric part, whereas the convergence of $\Escr^n(w_n,u_n)$ for $w_n$ in a dense set requires a bound on $\Escr^n_1(\cdot,u_n)$ as a functional on $\Vscr_n$. This is the r\^ole of the functional $\Theta$ defined below, which replaces the energy bound of the symmetric case. The functional $\Theta$ measures the size of $\Escr_1(\cdot,u)$ as a functional on $\Vscr$ and is, indeed, an equivalent norm to $\norm{\;}_{\hat{\Fscr}}$ by Lemma~\ref{PsiEquiv} below. So let
\[\Theta(u):=\sup_{\norm{v}_\Vscr =1}\Escr_1(v,u)=\norm{\Escr_1(\cdot,u)}_{\Vscr'},\quad\text{for}\;\;u\in\hat{\Fscr},\]
which is finite. If $u\in H\setminus\hat{\Fscr}$, we extend $\Theta$ to a functional on $H$ with values in $\ol{\R}$ by setting $\Theta(u):=+\infty$.
\begin{lem}\label{PsiEquiv}
For $u\in\hat{\Fscr}$, we have
\begin{enumerate}[(i)]
\item $\Theta(u)\le \sqrt{K^2+1}\norm{u}_{\hat{\Fscr}}$,
\item $\norm{u}_\Vscr\le\Theta(u)$,
\item $\norm{\hat{\Lambda}u}_{\Vscr'}\le (K+1)\Theta(u)$,
\item In particular, $\norm{\;}_{\hat{\Fscr}}\sim\Theta(\;)$.
\end{enumerate}
\end{lem}
\begin{proof}
\begin{enumerate}[(i):]
\item \begin{align*}\Theta(u)&\le\sup_{\norm{v}_\Vscr =1}\lrabs{\Ascr_1(v,u)}
+\sup_{\norm{v}_\Vscr =1}\lrabs{\dualdel{\Vscr}{\hat{\Lambda}u}{v}}\\
&\le K\norm{u}_\Vscr+\norm{\hat{\Lambda}u}_{\Vscr'}\le \sqrt{K^2+1}\norm{u}_{\hat{\Fscr}},
\end{align*}
where the last step follows from the Cauchy--Schwarz inequality in $\R^2$ and the definition of $\norm{\;}_{\hat{\Fscr}}$.
\item
\[
\norm{u}_\Vscr^2\le\Escr_1(u,u)=\norm{u}_\Vscr \Escr_1\left(\frac{u}{\norm{u}_\Vscr},u\right)
\le\norm{u}_\Vscr \supl_{\norm{v}_\Vscr=1}\Escr_1 (v,u)=\norm{u}_\Vscr \Theta(u).
\]
\item For $v\in\Vscr$ we have by (ii)
\begin{align*}
\dualdel{\Vscr}{\hat{\Lambda}u}{v}&=\Ascr_1(v,u)-\Escr_1(v,u)\\
&\le K\norm{v}_\Vscr \norm{u}_\Vscr +\norm{v}_\Vscr\Theta(u)\\
&\le (K+1)\norm{v}_\Vscr\Theta(u).
\end{align*}
Hence $\norm{\hat{\Lambda}u}_{\Vscr'}\le(K+1)\Theta(u)$.
\item Obvious by (i)--(iii).
\end{enumerate}
\end{proof}

We arrive at the main convergence theorem of this section. It is a generalization of~\cite[Theorem 3.1]{H} by M. Hino. From now on,
we consider given forms $\{\Escr^n\}$ and $\Escr$ on $H_n$ and $H$, respectively. The operators, spaces and norms related to
$\Escr^n$ are denoted by adding the index $n$, such as $G_\a^n$, $\Theta^n$ and $\Vscr_n$. Note that the
sector constants $K_n$ of the $\Ascr^n$'s are not necessarily assumed to be uniformly bounded.

\begin{defi}\label{Fs}
Suppose $C\subset\Fscr$ densely w.r.t. $\norm{\;}_\Fscr$. Consider the following conditions:
\begin{enumerate}
\item[\textup{(F1)}] If a sequence $\{u_n\}$ weakly convergent to $u$ in $\Hscr$ satisfies
$\varliminf_n \Theta^n(u_n)<\infty$, then $u\in\Vscr$.
\item[\textup{(F2)}] For any $w\in C$, any $u\in\Vscr$ and any sequence $\{u_n\}$ weakly convergent to $u$ in $\Hscr$, $u_n\in\Vscr_n$, $n\in\N$, there exists $\{w_n\}$, $w_n\in H_n$, $n\in\N$ converging to $w$ strongly in $\Hscr$ such that $\lim_n \Escr^n(w_n,u_n)=\Escr(w,u)$.
\item[\textup{(F2')}] There exists a linear subspace $\tilde{C}$ of $C$ such that
$\tilde{C}\subset\Fscr$ densely w.r.t. $\norm{\;}_\Fscr$ and for any sequence $n_k\uparrow\infty$ and every $w\in\tilde{C}$, $u\in\Vscr$ and any sequence $\{u_k\}$, $u_k\in H_{n_k}$, $k\in\N$ converging weakly to $u$ and satisfying
$\sup_k\Theta^{n_k}(u_k)<\infty$ one has a
sequence $\{w_k\}$, $w_k\in H_{n_k}$, $k\in\N$ converging $\Hscr$-strongly to $w$
with $\varliminf_k \Escr^{n_k}(w_k,u_k)\le\Escr(w,u)$.
\item[\textup{(R)}] $\{G_\a^n\}$ converges to $G_\a$ strongly for $\a>0$.
\item[\textup{(CR)}] $\{\hat{G}_\a^n\}$ converges to $\hat{G}_\a$ weakly for $\a >0$.
\end{enumerate}

Define also \textup{(F1a)} and \textup{(F1b)} by replacing $\Theta^n(u_n)$ in \textup{(F1)} by $\norm{u_n}_{\hat{\Fscr}_n}$ and by $\norm{u_n}_{\Vscr_n}$, respectively. Similarly, define \textup{(F2'a)} and \textup{(F2'b)} by replacing $\Theta^{n_k}(u_{k})$ in \textup{(F2')} by $\norm{u_{k}}_{\hat{\Fscr}_{n_k}}$ and by $\norm{u_{k}}_{\Vscr_{n_k}}$, respectively. Here we set $\norm{u}_{\hat{\Fscr}_n}:=+\infty$ for $u\in H_n\setminus\hat{\Fscr}_n$ and $\norm{u}_{\Vscr_n}:=+\infty$ for $u\in H_n\setminus\Vscr_n$.
\end{defi}
We have (F1b) $\Rightarrow$ (F1) and (F2'b) $\Rightarrow$ (F2') by Lemma~\ref{PsiEquiv} (ii).
Condition (F1) ensures that weak limits of sequences for which $\Theta^n$ is bounded belong to the limiting form domain $\Vscr$, whereas the conditions (F2) and (F2') describe the limit behavior of $\Escr^n(w_n,u_n)$ along recovery sequences $\{w_n\}$ and thereby handle the non-symmetric part. The following theorem characterizes (R) and (CR) by these conditions. For symmetric forms, combined with the theorem of Mosco, Kuwae and Shioya, Theorem~\ref{KSConv}, it shows that (F1) and (F2) are equivalent to Mosco convergence, see Remark~\ref{symmosco}.

\begin{thm}\label{Hino}
Suppose that $C\subset\Fscr$ densely w.r.t. $\norm{\;}_\Fscr$.
Then
\[\textup{(F2)} \;\;\Rightarrow\;\; \textup{(F2')},\]
\[\textup{(F1)(F2')}\;\;\Leftrightarrow\;\;\textup{(F1)(F2)}\;\;\Leftrightarrow\;\;\textup{(R)}
\;\;\Leftrightarrow\;\;\textup{(CR)}.\]
\end{thm}
\begin{proof}
The equivalence between (R) and (CR) follows immediately from Corollary~\ref{Badjoint} (1).

\textit{(F2) $\Rightarrow$ (F2').} By choosing $u_n=0$ for every $n$ in (F2), we see that for each $w\in C$ there exists $\{w_n\}$ converging to $w$ strongly in $\Hscr$ such
that $\Escr^n(w_n,0)$ is finite for all sufficiently large $n$, hence $w_n\in\Vscr_n$ for all large $n$ by the convention that $\Escr^n(v,0)=+\infty$ for $v\notin\Vscr_n$. Changing finitely many terms, we may assume $w_n\in\Vscr_n$ for every $n$. Since $C$ is dense in $\Fscr$ and therefore also dense in $\Vscr$, we can find, by Corollary
\ref{doubletrouble}, for any $w\in\Vscr$ a sequence $\{w'_n\}$ converging strongly to $w$ with $w'_n\in\Vscr_n$, $n\in\N$. Take $\tilde C:=C$, an arbitrary strictly increasing sequence $n_k\uparrow\infty$, $u\in\Vscr$ and $\{u_k\}$, $u_k\in H_{n_k}$, $k\in\N$, weakly convergent to $u$ in $\Hscr$ with $\sup_k\Theta^{n_k}(u_k)<\infty$.
Then $u_k\in\hat{\Fscr}_{n_k}\subset\Vscr_{n_k}$. Let $\{w'_n\}$ be as above for $w:=u$, and define $u'_{n_k}:=u_k$, $k\in\N$, and $u'_n:=w'_n$ for $n\notin\{n_k\;\vert\;k\in\N\}$. As in the proof of Corollary~\ref{exstrongconv} (iii), $u'_n\to u$ weakly, and $u'_n\in\Vscr_n$ for every $n$. Take an arbitrary $w\in\tilde{C}$. By (F2) there exists a sequence $\{w_n\}$ converging strongly to $w$ with $\lim_n\Escr^n(w_n,u'_n)=\Escr(w,u)$. Setting $w_k:=w_{n_k}$, we obtain a sequence converging strongly to $w$ along $\{n_k\}$ with $w_k\in\Vscr_{n_k}$ for large $k$ and
\[
\lim_k\Escr^{n_k}(w_k,u_k)=\Escr(w,u),
\]
which is even stronger than required in (F2').

\textit{(F1),(F2') $\Rightarrow$ (CR).} We follow the argument of R\"ockner and Zhang in~\cite{RZ2}. Choose a sequence $\{f_n\}$ $\Hscr$-weakly
convergent to some $f\in H$. Let $\a>0$. We shall prove that $\hat{G}_\a^n f_n$ converges $\Hscr$-weakly to $\hat{G}_\a f$.
It suffices to prove that for any sequence $n_k\uparrow\infty$ we can extract a subsequence $\{n_{k_l}\}$ such
that $\hat{G}_\a^{n_{k_l}} f_{n_{k_l}}$ converges to $\hat{G}_\a f$ weakly in $\Hscr$.
Set $u_n:=\hat{G}_\a^n f_n$. Since $\norm{\hat{G}_\a^n}_{\Lscr(H_n)}\le\a^{-1}$ and
$\{f_n\}$ is uniformly bounded by Lemma~\ref{dualconv} (2), one can extract a subsequence $\{u_{n_{k_l}}\}$ converging $\Hscr$-weakly to some $u\in H$ using Lemma
\ref{dualconv} (1). For notational convenience, we denote this subsequence again by $\{u_{n_k}\}$.

Note that for some constant $C_\a>0$ depending only on $\a$ by the Cauchy--Schwarz inequality and the uniform boundedness of $\{f_n\}$
and since $u_{n_k}=\hat{G}_\a^{n_k} f_{n_k}\in\hat{\Fscr}_{n_k}$
\begin{multline*}
\supl_k\Theta^{n_k}(u_{n_k})=\supl_k\supl_{\norm{w}_{\Vscr_{n_k}}=1}\Escr_1^{n_k}(w,u_{n_k})\\\le 
\supl_k\supl_{\norm{w}_{\Vscr_{n_k}}=1}\left[\lrabs{\Escr_\a^{n_k}(w,u_{n_k})}
+\abs{1-\a}\abs{(w,u_{n_k})_{H_{n_k}}}\right]\\
\le\supl_k\supl_{\norm{w}_{\Vscr_{n_k}}=1}\left[\abs{\deli{f_{n_k}}{w}_{H_{n_k}}}+\frac{\abs{1-\a}}{\a}\norm{w}_{H_{n_k}}\norm{f_{n_k}}_{H_{n_k}}\right]\\
\le C_\a\supl_k\supl_{\norm{w}_{\Vscr_{n_k}}=1}\norm{f_{n_k}}_{H_{n_k}}\norm{w}_{H_{n_k}}\le C_\a \supl_k \norm{f_{n_k}}_{H_{n_k}}<\infty
\end{multline*}
which yields $u\in\Vscr$ by (F1), applied to the extension of $\{u_{n_k}\}$ to a weakly convergent full sequence according to Corollary~\ref{exstrongconv} (iii). Note that the $\varliminf$ in (F1) along the full sequence is bounded by the $\varliminf$ along $\{n_k\}$.

For any $w\in \tilde{C}$, by extracting a subsequence if necessary, we can choose $\{w_k\}$ strongly convergent to $w$ in $\Hscr$ such that
$w_k\in\Vscr_{n_k}$ and
\[\varliminf_k \Escr^{n_k}(w_k,u_{n_k})\le\Escr(w,u)\]
by (F2'). Since $\Escr_\a^{n_k}(w_k,u_{n_k})=\deli{w_k}{f_{n_k}}_{H_{n_k}}$, it follows that
\[0=\varliminf_k\left[\Escr_\a^{n_k}(w_k,u_{n_k})-\deli{w_k}{f_{n_k}}_{H_{n_k}}\right]\le\Escr_\a (w,u)-\deli{w}{f}_H.\]
Hence $\Escr_\a(w,u)\ge \deli{w}{f}_H$. By substituting $-w$ for $w$, this becomes an equality. Now we have by Proposition~\ref{Resolvent}, using that $\tilde{C}\subset C\subset\Fscr$ densely and
continuously, that for every $w\in\Fscr$
\[\Escr_\a(w,u)=(w,f)_H=\Escr_\a(w,\hat{G}_\a f).\]
Hence $\dualdel{\Vscr}{(M_\a-\Lambda)w}{u-\hat{G}_\a f}=0$, but $(M_\a-\Lambda)(\Fscr)=\Vscr'$ by
\cite[Section I.3]{S1}, hence $u=\hat{G}_\a f$.

This yields $\hat{G}_\a^n f_n\to\hat{G}_\a f$ weakly in $\Hscr$. Since $f_n$ was an arbitrary weakly convergent sequence, we
have $\hat{G}_\a^n \to\hat{G}_\a$ weakly in the sense of convergence of bounded operators.

\textit{(R) $\Rightarrow$ (F1).} Let $u_n\to u$ weakly in $\Hscr$ and $M:=\varliminf_n \Theta^n(u_n)<\infty$. Choose a sequence $\{v_n\}$ converging strongly to $u$ in $\Hscr$.
From Proposition~\ref{BetaMR} (i),(ii), for $\a>0$,
\begin{align*}
&\a\deli{v_n-\a G_\a^n v_n}{u_n}_{H_n}+\deli{\a G_\a^n v_n}{u_n}_{H_n}\\
=&\Escr_1^n(\a G_\a^n v_n,u_n)\\
\le&\Theta^n(u_n)\norm{\a G_\a^n v_n}_{\Vscr_n}\\
\le&\Theta^n(u_n)\left[\Escr_1^n(\a G_\a^n v_n,\a G_\a^n v_n)\right]^{1/2}\\
\le&\Theta^n(u_n)\Escr_1^{n,(\a)}(v_n,v_n)^{1/2}\\
=&\Theta^n(u_n)\left[\a\deli{v_n-\a G_\a^n v_n}{v_n}_{H_n}+\norm{v_n}_{H_n}^2\right]^{1/2}
\end{align*}
If $u_n\notin\hat{\Fscr}_n$, then $\Theta^n(u_n)=+\infty$ and the inequality is trivial. Since $\varliminf_n\Theta^n(u_n)<\infty$, we may pass to a subsequence along which $\Theta^n(u_n)$ converges to $M$, which is sufficient for the following argument. Taking $\varliminf_n$ on both sides, we have
\[\Escr_1^{(\a)}(u,u)-\norm{u}_H^2+\deli{\a G_\a u}{u}_H\le M\Escr_1^{(\a)}(u,u)^{1/2}.\]
Hence
\[\Escr_1^{(\a)}(u,u)^{1/2}\le\frac{1}{2}\left[M+\sqrt{M^2+4\left(\norm{u}_H^2-\deli{\a G_\a u}{u}_H\right)}\right],\]
which implies that
\[\sup_{\a>0}\Escr_1^{(\a)}(u,u)^{1/2}\le\frac{1}{2}\left[M+\sqrt{M^2+8\norm{u}^2_H}\right]<\infty.\]
From Proposition~\ref{BetaMR} (iv), we obtain that $u\in\Vscr$.

\textit{(R) $\Rightarrow$ (F2).} Let $u_n\to u$ weakly in $\Hscr$, $u_n\in\Vscr_n$, $u\in\Vscr$ and $w\in C$. Since $\Phi_n(w)\to w$ $\Hscr$-strongly, we have
\[d_\Hscr-\liml_{\a\to\infty}d_\Hscr-\liml_{n\to\infty}\a G_\a^n \Phi_n(w)=w,\]
and
\begin{align*}
\liml_{\a\to\infty}\liml_{n\to\infty} \Escr^{n,(\a)}(\Phi_n(w),u_n)&=\liml_{\a\to\infty}\liml_{n\to\infty}\a\deli{\Phi_n(w)-\a G_\a^n\Phi_n(w)}{u_n}_{H_n}\\
&=\liml_{\a\to\infty}\Escr^{(\a)}(w,u)=\Escr(w,u).
\end{align*}
Due to Corollary~\ref{doubletrouble}, applied with $\a=m\in\N$, we can take a nondecreasing sequence $\{\a_n\}$$\subset\N$, $\a_n\to\infty$ such that
\[d_\Hscr-\liml_{n\to\infty}\a_n G_{\a_n}^n\Phi_n(w)=w,\quad\liml_{n\to\infty}\Escr^{n,(\a_n)}(\Phi_n(w),u_n)=\Escr(w,u).\]
Recall that by Proposition~\ref{Resolvent} $G_\a^n(H)\subset\Fscr_n$ for any $\a>0$ and $n\in\N$. Setting $w_n:=\a_n G_{\a_n}^n\Phi_n(w)$, we hence have $w_n\to w$ $\Hscr$-strongly and by
Proposition~\ref{BetaMR} (i) that
\[\Escr^n(w_n,u_n)=\Escr^{n,(\a_n)}(\Phi_n(w),u_n)\to\Escr(w,u)\]
as $n\to\infty$ and (F2) is proved.
\end{proof}

\begin{cor}\label{Knbdd}
Suppose that $C\subset\Fscr$ densely.
\begin{enumerate}[(i)]
\item (F1b),(F2'b) $\Rightarrow$ (R),
\item If the sector constants $K_n$ of the $\Ascr^n$'s are uniformly bounded, then (F1a),(F2'a) $\Leftrightarrow$ (R).
\end{enumerate}
\end{cor}
\begin{proof}
\begin{enumerate}[(i):]
\item This is clear, since (F1b) $\Rightarrow$ (F1) and (F2'b) $\Rightarrow$ (F2') by Lemma~\ref{PsiEquiv} (ii).
\item This is a consequence of Theorem~\ref{Hino} and Lemma~\ref{PsiEquiv}.
\end{enumerate}
\end{proof}

\begin{defi}
We say that a sequence of generalized forms $\{\Escr^n\}$ \textit{converges in the generalized sense} to a generalized form $\Escr$, if
$C\subset\Fscr$ densely and
(F1) and (F2), or equivalently (F1) and (F2'), hold. Since (F1b) $\Rightarrow$ (F1) and (F2'b) $\Rightarrow$ (F2') by Lemma~\ref{PsiEquiv} (ii), this is in particular the case if (F1b) and either (F2) or (F2'b) hold. No bound on the sector constants is needed for this. We shall also use this notion if (F1a) and either (F2) or (F2'a) hold,
provided the sector constants $K_n$
of the $\Escr^n$'s are uniformly bounded, see Corollary~\ref{Knbdd} (ii).
\end{defi}

\begin{rem}
According to Theorem~\ref{Kato} and Corollary~\ref{Badjoint} (1), statements corresponding to (R) and (CR) also hold for the associated semigroups $(T_t^n)_{t\ge 0}$, $n\in\N$, and $(T_t)_{t\ge 0}$
and the associated co-semigroups $(\hat{T}_t^n)_{t\ge 0}$, $n\in\N$, and $(\hat{T}_t)_{t\ge 0}$.
\end{rem}

We now prove some properties of convergence associated with the above theorem. First we point out that, if $C\subset\Fscr$ densely,
$C$ contains enough elements to determine the behavior of a sequence $\Escr^n$ along particular strongly convergent sequences, so that (F2) is less
restrictive in comparison with (M2) than it might seem. This is illustrated by the following proposition.

\begin{prop}\label{mataloni}
Let $\Escr^n$, $n\in\N$, $\Escr$ be as in Theorem~\ref{Hino} and $C\subset\Fscr$ densely.
Assume that (F2) holds. Then the following stronger version of (F2) holds:

For any $w\in\Fscr$ and any $u_n\in\Vscr_n$, $n\in\N$, $u\in\Vscr$ with $u_n\to u$ $\Hscr$-weakly there exists a sequence $\{w_n\}$, $w_n\in H_n$, $n\in\N$, $w_n\to w$
$\Hscr$-strongly such that
\[\lim_n\Escr^n(w_n,u_n)=\Escr(w,u).\]
\end{prop}
\begin{proof}
Let $w\in\Fscr$ and $u_n\in\Vscr_n$, $n\in\N$, $u\in\Vscr$ with $u_n\to u$ $\Hscr$-weakly. Then there exists a sequence $\{w^m\}\subset C$ such that
$w^m\to w$ strongly in $\norm{\;}_\Fscr$, and therefore in $\norm{\;}_H$.
Now pick by (F2) for every $w^m$ a sequence $\{w^m_n\}\subset\Hscr$, $w_n^m\in H_n$, $n\in\N$ such that $w_n^m\to w^m$ $\Hscr$-strongly and
\begin{equation}\label{Escrwnm}
\lim_n\Escr^n(w_n^m,u_n)=\Escr(w^m,u).
\end{equation}
Clearly,
\[\lim_m\Escr(w^m,u)=\Escr(w,u).\]
By \eqref{Escrwnm}, $\Escr^n(w_n^m,u_n)$ is finite for all sufficiently large $n$. Hence, by Theorem~\ref{Kolmetric}, Corollary~\ref{doubletrouble} and the remark after it, there exists a nondecreasing sequence
of natural numbers $m(n)\uparrow\infty$ with
\[w^{m(n)}_n\to w\quad\text{$\Hscr$-strongly}\]
and
\[\lim_n\Escr^n(w_n^{m(n)},u_n)=\Escr(w,u).\]
The proof is complete.
\end{proof}

In the case of coercive closed forms, generalized convergence
implies an even stronger version of (F2) than in the above proposition. More
precisely, we have:

\begin{prop}\label{matalonistronger}
Let $\Escr^n=\Ascr^n$, $n\in\N$, $\Escr=\Ascr$ be coercive closed forms. 
Assume that $C\subset\Vscr$ densely and that $C\subset D(A)$, where $A$ denotes the
infinitesimal generator.
Then (F1),(F2)
are equivalent to (F1) and the following condition:
\begin{enumerate}
\item[\textup{(FF2)}]
For every $w\in D(A)\subset\Vscr$ there exists a sequence $\{w_n\}$,
$w_n\in H_n$, $n\in\N$, converging $\Hscr$-strongly to $w$, such that
\[\lim_n\Escr^n(w_n,u_n)=\Escr(w,u)\]
for any $\{u_n\}$, $u_n\in \Vscr_n$, $n\in\N$, $u\in\Vscr$ with
$u_n\to u$ $\Hscr$-weakly.
\end{enumerate}
\end{prop}
The assumption $C\subset D(A)$ is needed for the implication (F1),(FF2) $\Rightarrow$ (F1),(F2), since (F2) is required for $w\in C$.
\begin{proof}
It is clear that (F1),(FF2) $\Rightarrow$ (F1),(F2). Let us prove the converse.
Assuming that (F1),(F2) hold, Theorem~\ref{Hino} tells us that (R) holds.
Take $w\in D(A)$ and define $v:=(1-A)w\in H$. We have $G_1^n\Psi_n v\to G_1 v$ $\Hscr$-strongly by (R).
Set $w_n:=G_1^n\Psi_n v\in D(A^n)$. We have that
\[(1-A^n)w_n=\Psi_n v\to v=(1-A)w\]
$\Hscr$-strongly. Hence $A^n w_n\to Lw$ $\Hscr$-strongly. But for every $n\in\N$
\[\Escr^n(w_n,u_n)=(-A^n w_n,u_n)_{H_n}\]
making sense for every $u_n\in\Vscr_n$.

Now take any $u_n\in\Vscr_n$, $n\in\N$, $u\in\Vscr$ with $u_n\to u$ $\Hscr$-weakly.
Then,
\[\lim_n\Escr^n(w_n,u_n)=\lim_n(-A^n w_n,u_n)_{H_n}=(-Lw,u)_H=\Escr(w,u)\]
by strong convergence of $A^n w_n\to Lw$.
\end{proof}

Concerning generalized convergence, we are interested in the question
whether the convergence $\Escr^n\to\Escr$ in the sense of Theorem~\ref{Hino} is sufficient for
the Mosco- or $\Gamma$-convergence of the symmetric parts $\tilde{\Ascr}^n\to\tilde{\Ascr}$. In general, this requires much stronger assumptions, such as uniformly bounded sector constants and conditions similar to (F1) and (F2) for the co-forms $\hat{\Escr}^n$, $\hat{\Escr}$. Conversely, if the symmetric parts satisfy (M1), then (F1) holds by Proposition~\ref{moscotohino} below, so that only (F2) remains to be verified.

\begin{rem}\label{symmosco}
By Theorems~\ref{KSConv} and \ref{Hino}, in the case of symmetric forms and provided that $C\subset\Vscr$ densely, (F1) and (F2) are just another characterization of Mosco convergence.
\end{rem}

If we compare Mosco convergence and generalized convergence, the main question is: can we benefit from the fact that the symmetric parts $\tilde{\Ascr}^n$ of
given forms $\Escr^n$ converge Mosco to a symmetric form $\tilde{\Ascr}$ being the symmetric part of a form $\Escr$? As we will see later,
the difficult part is to prove (M1) since (M2) can be obtained easily in most applications using Lemma~\ref{densemosco}. Accordingly,
assuming Mosco convergence of the symmetric parts, only (F2) remains to be verified, as the following proposition shows.

\begin{prop}\label{moscotohino}
Let $\Escr^n$, $n\in\N$, $\Escr$ be as in Theorem~\ref{Hino}.
If the symmetric parts $\tilde{\Ascr}^n$, $n\in\N$, and $\tilde{\Ascr}$ of $\Escr^n$, $n\in\N$, and $\Escr$, respectively, fulfill
(M1), then (F1) holds.
\end{prop}
\begin{proof}
Let $u_n\to u$ $\Hscr$-weakly with $\varliminf_n\Theta^n(u_n)<\infty$. Using Lemma~\ref{PsiEquiv} (ii) and (M1)
\[\tilde{\Ascr}(u,u)\le\varliminf_n\tilde{\Ascr}^n(u_n,u_n)\le\varliminf_n\norm{u_n}_{\Vscr_n}^2\le\varliminf_n(\Theta^n(u_n))^2<\infty,\]
which yields $u\in\Vscr$.
\end{proof}

\begin{rem}\label{LL}
Obviously, if $\Ascr^n\equiv 0$, $\Vscr_n=H_n$, $n\in\N$, $\Ascr\equiv 0$, $\Vscr=H$, i.e., our generalized forms $\{\Escr^n\}$
depend only on operators $(\Lambda^n,D(\Lambda^n,H_n))=(L^n,D(L^n))$, $n\in\N$, see \textup{\cite[Remark I.4.10]{S1}} for details,
condition (F1) can be omitted.
\end{rem}

\begin{rem}\label{lambdalealpha}
All of the above steps can be carried out if $(\Ascr^n,\Vscr_n)$, $n\in\N$ and $(\Ascr,\Vscr)$ are coercive closed forms up to a common shift,
that is, for some constant $\lambda\in\R$ independent of $n$,
\begin{enumerate}[(1)]
\item $(\tilde{\Ascr}_\lambda,\Vscr)$ is a non-negative, symmetric, closed form.
\item $(\Ascr_\lambda,\Vscr)$ satisfies the \textup{weak sector condition}, i.e., there exists a sector constant $K\ge 1$ such that
\[\lrabs{\Ascr_{\lambda+1}(u,v)}\le K\Ascr_{\lambda+1}(u,u)^{1/2}\Ascr_{\lambda+1}(v,v)^{1/2}\quad\text{for all}\;u,v\in\Vscr.\]
\end{enumerate}
The same is assumed for the $\Ascr^n$'s, where the sector constants can depend on $n$, but $\lambda$ cannot.
In this case the resolvents $G_\a^n$, $n\in\N$ and $G_\a$ are only defined for $\a>\lambda$. This case can be reduced to the case $\lambda=0$ by passing to the forms $\Escr^n+\lambda\deli{\cdot}{\cdot}_{H_n}$ and $\Escr+\lambda\deli{\cdot}{\cdot}_H$, whose resolvents are $G^n_{\a+\lambda}$ and $G_{\a+\lambda}$, respectively. The conditions (F1), (F2) and (F2') are invariant under such a common shift, since $\lambda\deli{w_n}{u_n}_{H_n}\to\lambda\deli{w}{u}_H$ whenever $w_n\to w$ strongly and $u_n\to u$ weakly. This also explains why $\lambda$ must not depend on $n$. In this setting, $\Theta^n$, $\Theta$ and the norms $\norm{\cdot}_{\Vscr_n}$, $\norm{\cdot}_\Vscr$ are defined from the shifted forms, i.e., with $\Escr^n_{\lambda+1}$, $\tilde{\Ascr}^n_{\lambda+1}$, $\Escr_{\lambda+1}$ and $\tilde{\Ascr}_{\lambda+1}$ in place of $\Escr^n_1$, $\tilde{\Ascr}^n_1$, $\Escr_1$ and $\tilde{\Ascr}_1$, respectively.
\end{rem}

\subsection{Strong graph limits and generalized forms}\label{ss:graph}

\begin{defi}
For each $n\in\N$ let $A_n$ be a closed linear operator on $H_n$ with dense linear domain $D(A_n)$.
$\{A_n\}$ is said to be \textit{convergent in the strong graph sense} if for each sequence $\{u_n\}$,
$u_n\in D(A_n)$, such that $u_n\to 0\in H$ $\Hscr$-strongly and the $\Hscr$-strong limit of $\{A_n u_n\}$
exists, we have that $A_n u_n\to 0\in H$ strongly in $\Hscr$.
\end{defi}

If $\{A_n\}$ converges in the strong graph sense, the following linear operator $(A,D(A))$ is well-defined:
\begin{align*}
D(A):=\big\{u\in H\;&\left\vert\right.\;\exists\{u_n\},\;u_n\in D(A_n),\;u_n\to u\in H\;\Hscr\text{-strongly},\\
&\quad A_n u_n\;\text{converges $\Hscr$-strongly}\big\},
\end{align*}
and for $u\in D(A)$
\begin{align*}Au:=\lim_n A_n u_n&\quad\text{with}\;\{u_n\}\;\text{such that}\;u_n\in D(A_n),\;u_n\to u\;\text{$\Hscr$-strongly}\\
&\quad\text{and}\;A_n u_n\;\text{converges $\Hscr$-strongly}.
\end{align*}
$(A,D(A))$ is called the \textit{strong graph limit} of $\{A_n\}$ in the sense of operators and we say that \textit{$\{A_n\}$ converges to $A$ in
the strong graph sense}.

\begin{defi}
For each $n\in\N$ let $A_n$ be a closed linear operator on $H_n$ with dense linear domain $D(A_n)$. The
\textit{strong graph limit} $\Gamma_\infty$ of $\{A_n\}$ in the sense of linear spaces is defined to be the set of
pairs $(u,v)\in H\times H$ such that there exists a sequence of vectors $\{u_n\}$, $u_n\in D(A_n)$
with $u_n\to u$ and $A_n u_n\to v$ strongly in $\Hscr$.
\end{defi}

We immediately obtain:
$\{A_n\}$ converges to some densely defined closed linear operator $A$ on $H$
in the strong graph sense if and only if
$\Gamma_\infty$ coincides with the \textit{graph} $\Gamma(A):=\{(u,Au)\in H\times H\;\vert\; u\in D(A)\}$ of $A$, that is, $\Gamma_\infty=\Gamma(A)$.

We remark that the conditions ``closed'' and ``densely defined'' can be relaxed, but turn out to be reasonable in the following
considerations.

The following theorem, which first appeared in~\cite{Toe1}, ensures that strong graph convergence is powerful enough to
characterize convergence of forms and resolvents, even in the non-symmetric case. Unfortunately, this notion of convergence is still hard
to handle, since in many cases the domains of the generators are difficult to specify explicitly.

\begin{thm}
Let $\Escr^n$, $n\in\N$, and $\Escr$ be generalized forms on $H_n$, $n\in\N$, and $H$, respectively. 
Suppose $C\subset\Fscr$ densely.
Denote by $(G_\a^n)_{\a>0}$, $n\in\N$, and $(G_\a)_{\a>0}$
the associated $C_0$-contraction resolvents, by $(T_t^n)_{t\ge 0}$, $n\in\N$, and $(T_t)_{t\ge 0}$ the associated $C_0$-contraction semigroups
and by $A^n$, $n\in\N$, and $A$ the associated infinitesimal generators. Then the following statements are equivalent:
\begin{enumerate}[(1)]
\item $G_\a^n\to G_\a$ strongly for $\a>0$.
\item $T_t^n\to T_t$ strongly for $t\ge 0$.
\item $A^n\to A$ in the strong graph sense.
\item $\Escr^n\to\Escr$ in the generalized sense.
\end{enumerate}
\end{thm}
\begin{proof}
(1) $\Leftrightarrow$ (2) follows from Theorem~\ref{Kato}.

(1) $\Leftrightarrow$ (4) follows from Theorem~\ref{Hino}.

Note that in our setting $(0,\infty)\subset\rho(A)$, where $\rho(\cdot)$ denotes the resolvent set,
and $A$ is closed and densely defined. The same holds for $A^n$, $n\in\N$.

\textit{(1) $\Rightarrow$ (3).} Assume that $G_\a^n\to G_\a$ strongly for any $\a>0$. We prove $\Gamma(A)=\Gamma_\infty$. Let
$u\in D(A)$ and set $u_n:=G_1^n\Psi_n(1-A)u$. Then $u_n\in D(A^n)$ for every $n$, and, since $\Psi_n(1-A)u\to(1-A)u$ $\Hscr$-strongly, the assumption gives $u_n\to G_1(1-A)u=u$ $\Hscr$-strongly, i.e., $\norm{u_n-\Psi_n u}_{H_n}\to 0$ by Theorem~\ref{Kolmetric}. Now,
\begin{align*}
&\norm{A^n u_n-\Psi_n Au}_{H_n}\\
=&\norm{A^n G_1^n\Psi_n(1-A)u-\Psi_n Au}_{H_n}\\
=&\norm{-(1-A^n)G_1^n\Psi_n(1-A)u+G_1^n\Psi_n(1-A)u-\Psi_n Au}_{H_n}\\
=&\norm{-\Psi_n(1-A)u+G_1^n\Psi_n(1-A)u-\Psi_n Au}_{H_n}=\norm{G_1^n\Psi_n(1-A)u-\Psi_n u}_{H_n},
\end{align*}
which clearly tends to $0$ as $n\to\infty$. Therefore $(u,Au)=\lim_n(u_n,A^n u_n)$ in $\Hscr\times\Hscr$, hence $\Gamma(A)\subset\Gamma_\infty$.

Now let $(u,v)\in\Gamma_\infty$. Then there exist $u_n\in D(A^n)$, $n\in\N$, $u_n\to u$ $\Hscr$-strongly and $A^n u_n\to v$ $\Hscr$-strongly.
Furthermore $(1-A^n)u_n\to u-v$ $\Hscr$-strongly and
\[u_n=G_1^n(1-A^n)u_n\to G_1 (u-v)=:w\in D(A)\]
by strong convergence of resolvents. But since $\Hscr$ is Hausdorff, e.g. by Theorem~\ref{Kolmetric}, we must have $u=w\in D(A)$. Furthermore,
\[
Au=Aw=AG_1(u-v)=-(1-A)G_1(u-v)+G_1(u-v)=v-u+\underbrace{G_1(u-v)}_{=w=u}=v.
\]
Hence $\Gamma_\infty\subset\Gamma(A)$. The assertion is proved.

\textit{(3) $\Rightarrow$ (1).} Let $A^n\to A$ in the strong graph sense. Let $\a>0$ and $v\in H$. We have that $v=(\a-A)u$ with $u:=G_\a v\in D(A)$.
By assumption there exist $u_n\in D(A^n)$, $u_n\to u$ $\Hscr$-strongly and $A^n u_n\to Au$. We apply Lemma~\ref{OpPhiConv}:
\begin{align*}
&\norm{G_\a^n\Psi_n v-\Psi_n G_\a v}_{H_n}\\
=&\norm{G_\a^n\Psi_n(\a-A)u-\Psi_n G_\a(\a-A)u}_{H_n}\\
\le&\norm{G_\a^n\Psi_n(\a-A)u-u_n}_{H_n}+\norm{u_n-\Psi_n u}_{H_n}\\
=&\norm{G_\a^n\Psi_n(\a-A)u-G_\a^n(\a-A^n)u_n}_{H_n}+\norm{u_n-\Psi_n u}_{H_n}\\
\le&\frac{1}{\a}\norm{\Psi_n(\a-A)u-(\a-A^n)u_n}_{H_n}+\norm{u_n-\Psi_n u}_{H_n},
\end{align*}
which clearly tends to $0$ as $n\to\infty$, since $(\a-A^n)u_n\to(\a-A)u$ strongly and $\Psi_n u\to u$ strongly.
The proof is complete.
\end{proof}

\subsection{Preservation of the sub-Markov property}\label{ss:subMarkov}

In this subsection, $E$ denotes a real locally convex topological vector space with $E\neq\{0\}$ and Borel $\sigma$-algebra $\Bscr(E)$. Let $m_n$, $n\in\N$, and $m$ be finite Borel measures on $E$ with full support such that $m_n\to m$ weakly. Here a measure has full support if every non-empty open subset of $E$ has positive measure. Let $C\subset C_b(E)$ be a linear subspace which is dense in $L^2(E;m)$, and let $\Phi_n$ be the identity on $C$. Since the functions in $C$ are continuous and the measures have full support, the maps $\Phi_n$ are injective, and $H_n:=L^2(E;m_n)\to H:=L^2(E;m)$ in the sense of Definition~\ref{Hilbertconv} by the weak convergence of the measures. These spaces are infinite dimensional, so that Assumption~\ref{infdim} holds. As usual, we set $\Hscr:=\disjoint_n H_n\dot\cup H$. For classes $f,g$ of functions on $E$, we write $f\le g$ if the inequality holds almost everywhere with respect to the measure under consideration. A version of Theorem~\ref{limitsubMarkov} below for finite dimensional state spaces and measures which are not necessarily finite can be found in~\cite[Proof of Theorem 3.5]{K1}.

\begin{defi}\label{subdefi}
A bounded linear operator $G$ on $L^2(E;m)$ is called \textit{sub-Markovian} if $0\le Gf\le 1$ for every $f\in L^2(E;m)$ with $0\le f\le 1$. A $C_0$-contraction semigroup $(T_t)_{t\ge0}$ is called \textit{sub-Markovian} if every $T_t$, $t\ge 0$, is sub-Markovian, and a $C_0$-contraction resolvent $(G_\a)_{\a>0}$ is called \textit{sub-Markovian} if every $\a G_\a$, $\a>0$, is sub-Markovian. The same definitions apply on $L^2(E;m_n)$.
\end{defi}

By~\cite[Proposition I.4.3]{MR}, a $C_0$-contraction resolvent is sub-Markovian if and only if the associated semigroup is sub-Markovian. By~\cite[Theorem I.4.4]{MR}, the resolvent of a coercive closed form $(\Escr,\Dscr(\Escr))$ on $L^2(E;m)$ is sub-Markovian if and only if $u^+\wedge 1\in\Dscr(\Escr)$ and $\Escr(u^+\wedge 1,u-u^+\wedge 1)\ge 0$ for every $u\in\Dscr(\Escr)$.

In applications, $C$ consists of smooth functions and is in general not stable under the unit contraction $u\mapsto u^+\wedge 1$. The following condition allows us to approximate the unit contraction within $C$.

\begin{defi}\label{smoothtrunc}
We say that $C$ \textit{admits smooth truncations} if $1\in C$ and if for every $\eps>0$ there exist nondecreasing functions $\phi_\eps:\R\to[-\eps,1+\eps]$ and $\psi_\eps:\R\to[-\eps,\infty)$ with Lipschitz constant $1$ such that $\phi_\eps(t)=t$ for $t\in[0,1]$, $\psi_\eps(t)=t$ for $t\ge0$, and $\phi_\eps\circ\phi\in C$ and $\psi_\eps\circ\phi\in C$ for every $\phi\in C$.
\end{defi}

For instance, if $E$ is a separable real Banach space, then $\FCb$ admits smooth truncations. Indeed, choose infinitely differentiable functions $k_\eps,h_\eps:\R\to[0,1]$ such that $k_\eps=1$ on $[0,1]$, $k_\eps=0$ outside $[-\eps,1+\eps]$, $h_\eps=1$ on $[0,\infty)$ and $h_\eps=0$ on $(-\infty,-\eps]$, and set $\phi_\eps(t):=\int_0^tk_\eps(s)\,ds$ and $\psi_\eps(t):=\int_0^th_\eps(s)\,ds$. If $\phi=f(l_1,\ldots,l_j)$ with $f\in C_b^\infty(\R^j)$ and $l_1,\ldots,l_j\in E'$, then $\phi_\eps\circ\phi=(\phi_\eps\circ f)(l_1,\ldots,l_j)$ and $\phi_\eps\circ f\in C_b^\infty(\R^j)$, since $f$ is bounded and all derivatives of $\phi_\eps$ are bounded. The same holds for $\psi_\eps$.

\begin{lem}\label{truncstrong}
Assume that $C$ admits smooth truncations, and let $\eps>0$. If $u_n\in L^2(E;m_n)$, $n\in\N$, converge $\Hscr$-strongly to $u\in L^2(E;m)$, then $\phi_\eps\circ u_n\to\phi_\eps\circ u$ $\Hscr$-strongly.
\end{lem}
\begin{proof}
Note that $\phi_\eps\circ u_n\in L^2(E;m_n)$ and $\phi_\eps\circ u\in L^2(E;m)$, since $\phi_\eps$ is bounded and the measures are finite. Let $\{\tilde u_j\}\subset C$ with $\tilde u_j\to u$ in $L^2(E;m)$ and $\lim_j\varlimsup_n\norm{\tilde u_j-u_n}_{L^2(E;m_n)}=0$. Then $\phi_\eps\circ\tilde u_j\in C$. Since $\abs{\phi_\eps(s)-\phi_\eps(t)}\le\abs{s-t}$ for $s,t\in\R$, we have $\norm{\phi_\eps\circ\tilde u_j-\phi_\eps\circ u}_{L^2(E;m)}\le\norm{\tilde u_j-u}_{L^2(E;m)}\to0$ and
\[\lim_j\varlimsup_n\norm{\phi_\eps\circ\tilde u_j-\phi_\eps\circ u_n}_{L^2(E;m_n)}\le\lim_j\varlimsup_n\norm{\tilde u_j-u_n}_{L^2(E;m_n)}=0.\]
Hence $\phi_\eps\circ u_n\to\phi_\eps\circ u$ $\Hscr$-strongly.
\end{proof}

\begin{lem}\label{constconv}
Assume that $1\in C$. Let $a_n,a\in\R$, $n\in\N$, be regarded as constant functions in $L^2(E;m_n)$ and in $L^2(E;m)$, respectively. Then $a_n\to a$ $\Hscr$-strongly if and only if $a_n\to a$ in $\R$.
\end{lem}
\begin{proof}
Since $1\in C$, the constant function $a$ belongs to $C$, and $\Phi_n(a)=a\to a$ $\Hscr$-strongly by Remark~\ref{phiphi} (1). By Lemma~\ref{basicconv} (1) and (2), $a_n\to a$ $\Hscr$-strongly if and only if $\norm{a_n-a}_{L^2(E;m_n)}=\abs{a_n-a}\,m_n(E)^{1/2}\to0$. Since $m_n(E)\to m(E)>0$ by the weak convergence of the measures, this is equivalent to $a_n\to a$ in $\R$.
\end{proof}

\begin{lem}\label{posweak}
Assume that $C$ admits smooth truncations. If $f_n\in L^2(E;m_n)$, $f_n\ge0$, $n\in\N$, converge $\Hscr$-weakly to $f\in L^2(E;m)$, then $f\ge0$.
\end{lem}
\begin{proof}
Let $g:=\mathbbm{1}_{\{f<0\}}$ and $\eps>0$. Choose $\phi_j\in C$, $j\in\N$, with $\phi_j\to g$ in $L^2(E;m)$, and set $\chi_j:=\psi_\eps\circ\phi_j+\eps$. Then $\chi_j\in C$, since $1\in C$, and $\chi_j\ge0$, since $\psi_\eps\ge-\eps$. For fixed $j$, we have $\Phi_n(\chi_j)\to\chi_j$ $\Hscr$-strongly by Remark~\ref{phiphi} (1), and hence
\[0\le\int_E f_n\chi_j\,dm_n\underset{n\to\infty}{\longrightarrow}\int_E f\chi_j\,dm\]
by the definition of weak convergence. Since $\psi_\eps$ has Lipschitz constant $1$ and $\psi_\eps\circ g=g$, we have $\chi_j\to g+\eps$ in $L^2(E;m)$ as $j\to\infty$. Hence $\int_E f(g+\eps)\,dm\ge0$. Letting $\eps\downarrow0$ and using that $f\in L^1(E;m)$, since $m$ is finite, we obtain $\int_{\{f<0\}}f\,dm\ge0$. Therefore $m(\{f<0\})=0$.
\end{proof}

\begin{thm}\label{limitsubMarkov}
Assume that $C$ admits smooth truncations. Let $(T^n_t)_{t\ge0}$ be sub-Markovian $C_0$-contraction semigroups on $L^2(E;m_n)$, $n\in\N$, and let $(T_t)_{t\ge0}$ be a $C_0$-contraction semigroup on $L^2(E;m)$ such that $T^n_t\to T_t$ strongly for every $t\ge0$. Then $(T_t)_{t\ge0}$ and the associated resolvent are sub-Markovian.
\end{thm}
\begin{proof}
Fix $t\ge0$ and $f\in L^2(E;m)$ with $0\le f\le1$. Let $u_n:=\Psi_nf$, so that $u_n\to f$ $\Hscr$-strongly by Theorem~\ref{Kolmetric}, and let $\eps_k\downarrow0$. Since $\phi_{\eps_k}\circ f=f$, Lemma~\ref{truncstrong} yields $\phi_{\eps_k}\circ u_n\to f$ $\Hscr$-strongly as $n\to\infty$ for every $k$. By Lemma~\ref{metricdoubletrouble}, applied in the metric space $(\Hscr,d_\Hscr)$, there exists a nondecreasing mapping $n\mapsto k(n)$ tending to $+\infty$ such that $f_n:=\phi_{\eps_{k(n)}}\circ u_n\to f$ $\Hscr$-strongly. Set $\delta_n:=\eps_{k(n)}$. Then $\delta_n\to0$ and $-\delta_n\le f_n\le1+\delta_n$, i.e., $0\le f_n^+\le1+\delta_n$ and $0\le f_n^-\le\delta_n$. Since $T^n_t$ is sub-Markovian, it maps functions with values in $[0,b]$, $b\ge0$, to functions with values in $[0,b]$. Hence
\[-\delta_n\le-T^n_tf_n^-\le T^n_tf_n\le T^n_tf_n^+\le1+\delta_n.\]
By the strong convergence of the semigroups, $T^n_tf_n\to T_tf$ $\Hscr$-strongly. By Lemma~\ref{constconv}, the constants $\delta_n$ converge $\Hscr$-strongly to $0$, and the constants $1+\delta_n$ converge $\Hscr$-strongly to $1$. By Lemma~\ref{basicconv} (2) and (4), the non-negative functions $T^n_tf_n+\delta_n$ and $1+\delta_n-T^n_tf_n$ converge $\Hscr$-weakly to $T_tf$ and to $1-T_tf$, respectively. Lemma~\ref{posweak} yields $0\le T_tf\le1$. Hence $T_t$ is sub-Markovian, and so is the associated resolvent by~\cite[Proposition I.4.3]{MR}.
\end{proof}

\begin{rem}\label{subMarkovforms}
Combined with Theorems~\ref{Kato} and \ref{Hino}, Theorem~\ref{limitsubMarkov} yields the following. If generalized forms $\Escr^n$ on $L^2(E;m_n)$ with sub-Markovian resolvents converge in the generalized sense to a generalized form $\Escr$ on $L^2(E;m)$, then the resolvent and the semigroup of $\Escr$ are sub-Markovian. If, in addition, $\Escr$ is a coercive closed form, then $u^+\wedge1\in\Dscr(\Escr)$ and $\Escr(u^+\wedge1,u-u^+\wedge1)\ge0$ for every $u\in\Dscr(\Escr)$ by~\cite[Theorem I.4.4]{MR}. Since $\FCb$ admits smooth truncations, this applies in the setting of Subsection~\ref{ss:infaij}, where the reference measures are probability measures with full support.
\end{rem}

\section{Examples of Mosco convergence}\label{sec:mosco}

\subsection{Finite dimensional symmetric case}\label{ss:finsym}

We briefly recall the results of A.V. Kolesnikov from~\cite{K1}, since the convergence of elliptic non-symmetric
forms can be proved using the Mosco convergence of symmetric $a_{i,j}$-forms. In particular, if the considered forms have
Mosco convergent symmetric parts, we can apply Proposition~\ref{moscotohino}.

Let $d\ge 1$ and let $dx$ be the Lebesgue measure on $\R^d$ and $\Bscr(\R^d)$ the Borel $\sigma$-algebra of $\R^d$. Let $\abs{\cdot}$ denote
the $d$-dimensional Euclidean norm and $\lrbr{\cdot,\cdot}$ the $d$-dimensional Euclidean inner product. Denote by $C_0^\infty(\R^d)$ the set
of all real-valued infinitely differentiable functions with compact support.

\begin{assumption}[Convergence of speed measures]\label{speed}
$\sigma_n>0$ $dx$-a.e., $\sigma_n\in L^1_\loc (dx)$, $n\in\N$ and there exists a function $\sigma$ such that $\sigma>0$ $dx$-a.e., $\sigma\in L^1_\loc(\R^d;dx)$
and $\{\mu_n=\sigma_n dx\}$ tends to $\mu=\sigma dx$ vaguely, i.e.,
\[\lim_n\int_{\R^d}\psi\;\sigma_n dx=\int_{\R^d}\psi\;\sigma dx\]
for every continuous function $\psi$ with compact support.
\end{assumption}

Let $H_n:=L^2(\R^d;\mu_n)$, $n\in\N$, $H:=L^2(\R^d;\mu)$, $C=C_0^\infty(\R^d)$ and $\Phi_n$ the identity operators on $C$. Then $H_n\to H$ in the sense of Definition \ref{Hilbertconv}.
As above $\Hscr=\disjoint_n H_n\dot\cup H$. Note that $\Phi_n$ is well-defined since the measure $\mu$ has full support, and injective since the measures $\mu_n$ have full support. Both hold, since $\sigma,\sigma_n>0$ $dx$-a.e.

\begin{lem}\label{sigmastrongweak}
Consider the following statements:
\begin{enumerate}[(1)]
\item $f_n\to f$ in $L^2(\R^d;dx)$.
\item $\frac{f_n}{\sqrt{\sigma_n}}\to\frac{f}{\sqrt{\sigma}}$ in $\Hscr$.
\end{enumerate}
In both statements, convergence is understood either in the strong or in the weak sense.
Assume that $\sqrt{\sigma_n}\to\sqrt{\sigma}$ weakly in $L^2_\loc(\R^d;dx)$. Then (1) implies (2) for the strong convergence and (2) implies (1)
for the weak convergence. If, in addition, $\sqrt{\sigma_n}\to\sqrt{\sigma}$ strongly in $L^2_\loc(\R^d;dx)$, then (1) and (2) are equivalent.
\end{lem}
\begin{proof}
The proof is taken from~\cite[Lemma 3.1]{K1}. Suppose first that $\sqrt{\sigma_n}\to\sqrt{\sigma}$ weakly in $L^2_\loc (\R^d;dx)$. 
Let us prove (1) $\Rightarrow$ (2) for the strong convergence. To this end, let $f_n\to f$
in $L^2(\R^d;dx)$. Choose a sequence $\phi_m\in C_0^\infty(\R^d)$ such that $\phi_m\sqrt{\sigma}\to f$ in $L^2(\R^d;dx)$. Then $\phi_m\to\frac{f}{\sqrt{\sigma}}$
in $L^2(\R^d;\sigma dx)$. Hence
\begin{multline*}
\lim_m\varlimsup_n\int_{\R^d}\left(\phi_m-\frac{f_n}{\sqrt{\sigma_n}}\right)^2\sigma_n dx\\
=\lim_m\varlimsup_n\int_{\R^d}(\phi_m\sqrt{\sigma_n}-f_n)^2 dx=\lim_m\int_{\R^d} (\phi_m\sqrt{\sigma}-f)^2 dx=0,
\end{multline*}
by weak convergence of $\sqrt{\sigma_n}\to\sqrt{\sigma}$ in $L^2_\loc(\R^d;dx)$ and Assumption~\ref{speed}.

Let us prove (2) $\Rightarrow$ (1) for the weak convergence. To this end, let $\frac{f_n}{\sqrt{\sigma_n}}\to\frac{f}{\sqrt{\sigma}}$ weakly in $\Hscr$.
Then for every $\Hscr$-strongly convergent sequence $u_n\to u$ we have that
\[\lim_n\int_{\R^d}\frac{f_n}{\sqrt{\sigma_n}} u_n\sigma_n dx=\int_{\R^d}\frac{f}{\sqrt{\sigma}}u\;\sigma dx.\]
Now let $v_n\to v$ strongly in $L^2(\R^d;dx)$. Then by the assertion proved before we have that $\frac{v_n}{\sqrt{\sigma_n}}\to\frac{v}{\sqrt{\sigma}}$
$\Hscr$-strongly. Clearly,
\[\lim_n\int_{\R^d}f_n v_n dx=\lim_n\int_{\R^d}\frac{f_n}{\sqrt{\sigma_n}}\frac{v_n}{\sqrt{\sigma_n}}\sigma_n dx=
\int_{\R^d}\frac{f}{\sqrt{\sigma}}\frac{v}{\sqrt{\sigma}}\sigma dx=\int_{\R^d}fv\;dx.\]

Now suppose that $\sqrt{\sigma_n}\to\sqrt{\sigma}$ strongly in $L^2_\loc(\R^d;dx)$. We prove (2) $\Rightarrow$ (1) for
the strong convergence. To this end, let $\frac{f_n}{\sqrt{\sigma_n}}\to\frac{f}{\sqrt{\sigma}}$
strongly in $\Hscr$. Then there exists a sequence of $C_0^\infty(\R^d)$-functions $\{\phi_m\}$ such that $\phi_m\sqrt{\sigma}\to f$ in
$L^2(\R^d;dx)$ and
\[\lim_m\varlimsup_n\int_{\R^d}\left(\phi_m-\frac{f_n}{\sqrt{\sigma_n}}\right)^2\sigma_n dx=\lim_m\varlimsup_n\int_{\R^d}\left(\phi_m\sqrt{\sigma_n}-f_n\right)^2 dx=0.\]
Since $\lim_m\lim_n\phi_m\sqrt{\sigma_n}=f$, where all limits are taken in $L^2(\R^d;dx)$, one can find a subsequence $n_k$ such that
$\phi_k\sqrt{\sigma_{n_k}}$ tends to $f$ and
\[\lim_k\norm{\phi_k\sqrt{\sigma_{n_k}}-f_{n_k}}_{L^2(\R^d;dx)}=0.\]
Hence $f_{n_k}\to f$ in $L^2(\R^d;dx)$. Since we can do the same with every subsequence of $\{f_n\}$, we get that $f_n\to f$ in $L^2(\R^d;dx)$.

Let us prove (1) $\Rightarrow$ (2) for the weak convergence. To this end, let $f_n\to f$ weakly in $L^2(\R^d;dx)$. Then for every $u_n\to u$
strongly in $L^2(\R^d;dx)$ we have
\[\lim_n\int_{\R^d}f_n u_n dx=\int_{\R^d}fu\;dx.\]
Now let $v_n\to v$ strongly in $\Hscr$. By the assertion proved before we have that $v_n\sqrt{\sigma_n}\to v\sqrt{\sigma}$
strongly in $L^2(\R^d;dx)$. Then
\[\lim_n\int_{\R^d}\frac{f_n}{\sqrt{\sigma_n}}v_n\sigma_n dx=\lim_n\int_{\R^d}f_n v_n\sqrt{\sigma_n}dx
=\int_{\R^d}fv\sqrt{\sigma}dx=\int_{\R^d}\frac{f}{\sqrt{\sigma}}v\;\sigma dx.\]
The proof is complete.
\end{proof}

\begin{rem}\label{weakisstrong}
Under Assumption~\ref{speed}, weak convergence $\sqrt{\sigma_n}\to\sqrt{\sigma}$ in $L^2_\loc(\R^d;dx)$ already implies strong convergence in $L^2_\loc(\R^d;dx)$: if $\phi$ is continuous with compact support, then $\phi\sqrt{\sigma_n}\to\phi\sqrt{\sigma}$ weakly in $L^2(\R^d;dx)$, and $\norm{\phi\sqrt{\sigma_n}}_{L^2(dx)}^2=\int\phi^2\sigma_n\,dx\to\int\phi^2\sigma\,dx=\norm{\phi\sqrt{\sigma}}_{L^2(dx)}^2$ by vague convergence, which gives strong convergence of $\phi\sqrt{\sigma_n}$ in $L^2(\R^d;dx)$. Hence, under Assumption~\ref{speed}, the statements (1) and (2) of Lemma~\ref{sigmastrongweak} are equivalent both for strong and for weak convergence, as soon as $\sqrt{\sigma_n}\to\sqrt{\sigma}$ weakly in $L^2_\loc(\R^d;dx)$.
\end{rem}

For each $n\in\N$ consider the following symmetric form:
\[\Escr^n(f,g)=\sum_{i,j=1}^d\int_{\R^d}a_{i,j}^n(x)\frac{\partial f(x)}{\partial x_i}\pard{g(x)}{x_j}dx,\quad f,g\in C_0^\infty(\R^d).\]
Here $a_{i,j}^n$, $1\le i,j\le d$, are Borel locally integrable functions and $a_{i,j}^n=a_{j,i}^n$ for $1\le i,j\le d$.
We will denote by $A^n$ the symmetric $d\times d$-matrix
$(A^n(x))_{i,j}:=a_{i,j}^n(x)$, $x\in\R^d$. Then our forms can be written as
\[\Escr^n(f,g)=\int_{\R^d}\lrbr{A^n\nabla f,\nabla g}dx.\]
We suppose that the $A^n$'s are $dx$-a.e. positive definite.
We denote the elements of the inverse matrix $(A^n)^{-1}$ by $(a^{-1})_{i,j}^{n}$, $((A^n)^{-1}(x))_{i,j}=(a^{-1})_{i,j}^n(x)$, $x\in\R^d$.

For an arbitrary Borel function $f:\R^d\to\R$ let us define the set
\[R(f):=\left\{x\;\bigg\vert\;\exists\eps>0,\;\int_{\abs{x-y}\le\eps}\frac{dy}{\abs{f(y)}}<\infty\right\},\]
where we adopt the convention $\frac{1}{0}:=+\infty$. 
Evidently, $R(f)$ is the largest open set $V$ such that $\frac{1}{f}\in L^1_\loc(V;dx)$.
Let $R(A)$ be the largest open set $V$ such that $(a^{-1})_{i,j}\in L^1_\loc(V;dx)$ for $1\le i,j\le d$, where $A=(a_{i,j})$ denotes the limiting matrix as in Theorem~\ref{Kolconvfinite} below, and define $R(A^n)$ similarly.

We say that an arbitrary Borel function $f$ satisfies \textit{Hamza's condition} if for $dx$-a.e. $x\in\R^d$
$f(x)>0$ implies $x\in R(f)$. If $f>0$ $dx$-a.e., this is equivalent to $dx(\R^d\setminus R(f))=0$.

The following assumption ensures the closability of $(\Escr^n, C_0^\infty(\R^d))$. It is weaker than the standard assumption made
in~\cite[Section II.2.b)]{MR}.
\begin{assumption}\label{closass}
$R(A)\subset R(\sigma)$, $R(A^n)\subset R(\sigma_n)$, $dx(\R^d\setminus R(A^n))=dx(\R^d\setminus R(A))=0$.
\end{assumption}
Under Assumption~\ref{closass}, the minimal extensions $(\Escr^n_0,\Dscr(\Escr^n_0))$ and $(\Escr_0,\Dscr(\Escr_0))$ considered below are well-defined closed forms.

\begin{lem}
The form $(\Escr_{+},\Dscr(\Escr_{+}))$ defined by
\begin{align*}
\Dscr(\Escr_{+})&=\Bigg\{f\in L^2(\R^d;\sigma dx)\;\bigg\vert\;f\;\text{admits weak derivatives $\partial_i f$ in $R(A)$}\\
&\qquad\qquad\text{for every}\;i\in\{1,\ldots,d\},
\;\sum_{i,j=1}^d\int_{\R^d}a_{i,j}(x)\frac{\partial f(x)}{\partial x_i}\frac{\partial f(x)}{\partial x_j} dx<\infty\Bigg\},\\
\Escr_{+}(f,g)&=\sum_{i,j=1}^d\int_{\R^d}a_{i,j}(x)\frac{\partial f(x)}{\partial x_i}\frac{\partial g(x)}{\partial x_j} dx,\quad f,g\in\Dscr(\Escr_{+}),
\end{align*}
is closed.
\end{lem}
\begin{proof}
See~\cite[Lemma 3.4]{K1}.
\end{proof}
$(\Escr_{+},\Dscr(\Escr_{+}))$ is the ``maximal'' extension of $(\Escr,C_0^\infty(\R^d))$.
Let us denote by $\Dscr(\Escr_0)$ the completion of $C_0^\infty(\R^d)$ w.r.t. $\Escr_1^{1/2}(u)=\left(\Escr(u)+\norm{u}_{L^2(\R^d;\sigma dx)}^{2}\right)^{1/2}$.
The extension of $(\Escr,C_0^\infty(\R^d))$ to this space, denoted by $(\Escr_0,\Dscr(\Escr_0))$, is a closed form on $L^2(\R^d;\sigma dx)$ and is called
the ``minimal'' extension. Let us extend our notation to $(\Escr^n_0,\Dscr(\Escr^n_0))$ and $(\Escr_{+}^n,\Dscr(\Escr_{+}^n))$
in the obvious way. A.V. Kolesnikov has shown in~\cite{K1} that the convergence of the coefficients does in general not determine the Mosco limit, not even
in the one dimensional case. The limit may be the minimal extension, the maximal extension or a form in between. Hence additional geometric assumptions have to be stated to identify the Mosco limit. In particular,
we need to assume that the minimal and maximal extensions coincide.

\begin{thm}\label{Kolconvfinite}
Let $\sqrt{\sigma_n}\to\sqrt{\sigma}$ weakly in $L^2_\loc(\R^d;dx)$,
\[\sup_n\int_\Omega\frac{dx}{\sigma_n}<\infty,\quad\sup_n\int_\Omega\abs{(a^{-1})^n_{i,j}}dx<\infty,\quad 1\le i,j\le d\]
for every bounded domain $\Omega\in\Bscr(\R^d)$ and there exists a $d\times d$ matrix $A$ with Borel locally integrable coefficients, which is symmetric and $dx$-a.e.
positive definite. Assume that $(a^{-1})_{i,j}\in L^1_\loc(\R^d;dx)$ and
\[a_{i,j}^n dx\to a_{i,j} dx,\quad (a^{-1})^n_{i,j} dx\to (a^{-1})_{i,j} dx\]
vaguely.

Let
\[\Escr(\phi,\psi)=\int_{\R^d}\lrbr{A\nabla\phi,\nabla\psi}dx,\quad\phi,\psi\in C_0^\infty(\R^d).\]
Suppose that the minimal and maximal extensions of $(\Escr,C_0^\infty(\R^d))$ in $L^2(\R^d;\sigma dx)$ coincide.
Then $\Escr_0^n\to\Escr_0$ and $\Escr_{+}^n\to\Escr_0$ Mosco, where $\Escr_0=\Escr_+$ denotes the closure of $(\Escr,C_0^\infty(\R^d))$ in $L^2(\R^d;\sigma dx)$.
\end{thm}
\begin{proof}
See~\cite[Theorem 1.1]{K1}.
\end{proof}

\subsection{Infinite dimensional symmetric case}\label{idsc}

The following results are mainly based on the work of A.V. Kolesnikov in~\cite{K2}.
We reduce the general case to the case of a Gelfand triple, which is also used in the results on
generalized convergence in Subsection~\ref{ss:infaij}.

Throughout this subsection we fix a separable real Hilbert space $X$ with inner product $\lrbr{\;,\;}_X$ and norm
$\lrnorm{\;}_X=\lrbr{\;,\;}_X^{1/2}$ which
is densely and continuously embedded into some separable real Banach space $E$. We denote by $E'$ the
topological dual of $E$ and by $\lrdual{E}{\;,\;}:E'\times E\to\R$ the corresponding duality pairing. Denote
by $\Bscr(E)$ the Borel $\sigma$-field of $E$. Identifying
$X$ with its dual $X'$ we have that
\[E'\subset X'\equiv X\subset E\quad\text{densely and continuously}\]
and $\lrdual{E}{\;,\;}$ restricted to $E'\times X$ coincides with $\lrbr{\;,\;}_X$. Here $E'$ is
endowed with the operator norm $\lrnorm{\;}_{E'}:=\norm{\;}_{\Lscr(E,\R)}$.

Define the linear space of cylindrical test functions on $E$ by
\[\FCE:=\left\{f(l_1,\ldots,l_m)\;\left\vert\right.\; m\in\N,\;f\in C_0^\infty(\R^m),\;l_1,\ldots,l_m\in E'\right\}.\]
Here $C_0^\infty(\R^m)$ denotes the set of all real-valued infinitely differentiable functions on $\R^m$ with compact support.
Define for $u\in\FCE$ and $k\in E$ the following directional derivative
\[\frac{\partial u}{\partial k}(z):=\frac{d}{ds}u(z+sk)\bigg\vert_{s=0},\;z\in E.\]
Observe that if $u=f(l_1,\ldots,l_m)$, then
\[\frac{\partial u}{\partial k}=\sum\limits_{i=1}^m \frac{\partial f}{\partial x_i}(l_1,\ldots,l_m)\lrdual{E}{l_i,k}\in\FCE.\]
This shows that for fixed $u\in\FCE$ and $z\in E$, $h\mapsto\frac{\partial u}{\partial h}(z)$ is
a continuous linear functional on $X$. Define $\nabla_X u(z)\in X$ by
\begin{equation}\label{nabla}
\lrbr{\nabla_X u(z),h}_X=\frac{\partial u}{\partial h}(z),\; h\in X.
\end{equation}

To apply our results on Mosco convergence we fix a sequence $\{\mu_n\}$ of Borel probability measures on $E$
and a Borel probability measure $\mu$ on $E$.
\begin{assumption}\label{mnconv}
$\mu_n$ converges weakly to $\mu$.
\end{assumption}
Note that by Prokhorov's theorem, see~\cite[Theorem 5.2]{Bil},
$\{\mu_n\}$ is tight. 
\begin{assumption}\label{qinv}
$\mu$ is $h$-quasi-invariant for every $h\in X$, i.e.,
$\mu\circ T_h^{-1}$ is absolutely continuous with respect to $\mu$ where $T_h(z):=z-h$. More generally, a measure $\mu'$ is called quasi-invariant along $k\in X$ if $\mu'\circ T_{tk}^{-1}$ is absolutely continuous with respect to $\mu'$ for every $t\in\R$.
\end{assumption}
This implies that $\supp\mu=E$.

Now consider the following sequence of Hilbert spaces $\{H_n\}=\{L^2(E;\mu_n)\}$. Set $C:=\FCE$ and $H:=L^2(E;\mu)$.
Note that $C\subset H$ densely
by the Hahn--Banach theorem and a monotone class argument, see~\cite[Section II.3.a)]{MR} or \cite[Remark 3.1]{AR}.
Let $\Phi_n$ be the identity operator. It is well-defined since $\mu$ has full support. We assume moreover that $\supp\mu_n=E$ for all $n\in\N$, so that each $\Phi_n$ is injective.
Set $\Hscr:=\disjoint_n H_n\dot\cup H$. $\{H_n\}$ converges to $H$ by weak convergence of measures.

Now consider a weakly convergent sequence $h_n\to h$ of vectors from $X$, i.e.,
\[\lrbr{h_n,g}_X\to\lrbr{h,g}_X,\]
for every $g\in X$. Suppose additionally $\lrnorm{h_n}_X=1$ for every $n$, $\norm{h}_X=1$. Under this normalization, weak convergence $h_n\to h$ in $X$ is equivalent to convergence in norm.
Fix some $l \in E'$ such that $\lrbr{l,h_n}_X\not=0$ for every $n$ and $\lrbr{l,h}_X\neq0$, which exists by Baire's theorem, since, by the Hahn--Banach theorem, each of the sets $\{l\in E'\;\vert\;\lrbr{l,h_n}_X=0\}$, $n\in\N$, and $\{l\in E'\;\vert\;\lrbr{l,h}_X=0\}$ is a closed proper linear subspace, hence nowhere dense, in the Banach space $E'$.
Denote by $P_{h_n}$ the projection $P_{h_n}:E\to E_0:=\{x\in E\;\vert\;l(x)=0\}$
\[P_{h_n}(z):=z-\frac{\lrdual{E}{l,z}}{\lrbr{l,h_n}_X}h_n,\; z\in E.\]
It is well known that every probability measure $\mu_n$ has a conditional distribution $\rho_n(x,\cdot)$ on the real line such that, with $\nu_n:=\mu_n\circ (P_{h_n})^{-1}$, the following disintegration formula holds
\begin{equation}\label{disint}
\int_E u(z)\mu_n(dz)=\int_{E_0}\int_\R u(x+sh_n)\rho_n(x,ds)\nu_n(dx),
\end{equation}
see~\cite{DM}. 
\begin{assumption}\label{hnqinv}
Every $\mu_n$ is quasi-invariant along $h_n$ and $h_n\to h$ weakly.
\end{assumption}
It follows from this assumption that the $\rho_n(x,ds)$ have densities w.r.t. the Lebesgue measure, i.e.,
$\rho_n(x,ds)=\rho_n(x+sh_n)ds$ for $\nu_n$-a.e. $x\in E_0$.

\begin{assumption}\label{mutildebd}
The following sequence of Borel measures
\[\tilde{\mu}_n^{h_n}:=\frac{ds}{\rho_n (x+sh_n)}\nu_n(dx)\]
is uniformly bounded on all sets of the type $E_0^N:=\{z\in E\;\vert\;\lrabs{l(z)}\le N\}$. That is, the sequence
\[\int_E u(z)d\tilde{\mu}_n^{h_n}:=\int_{E_0}\int_\R u(x+sh_n)\frac{ds}{\rho_n(x+sh_n)}\nu_n(dx)\]
is bounded for every bounded Borel function $u:E\to\R$ with support in $E_0^N$.
\end{assumption}
In particular, for $\nu_n$-almost every $x$, and hence
$P_{h_n}(\mu_n)$-a.e., the function $s\mapsto(\rho_n(x+sh_n))^{-1}$ is locally integrable. We define $P_h$, $\nu:=\mu\circ P_h^{-1}$, $\rho$ and $\tilde\mu^h$ for $\mu$ and $h$ in the same way, recalling that $\lrbr{l,h}_X\neq0$. Since $\mu$ is $h$-quasi-invariant, $\rho(x,ds)=\rho(x+sh)ds$ for $\nu$-a.e. $x\in E_0$, and we assume that $s\mapsto(\rho(x+sh))^{-1}$ is locally integrable for $\nu$-a.e. $x\in E_0$.

For every $h_n$ consider the following partial form
\[\Dscr(\Escr_{\mu_n,0}^{h_n}):=\FCE,\]
and
\[\Escr_{\mu_n,0}^{h_n}(u,v):=\int_E \pard{u}{h_n}\pard{v}{h_n}d\mu_n,\quad u,v\in\FCE.\]

A sufficient condition for the closability of $\Escr_{\mu_n,0}^{h_n}$ is
that $\mu_n$ \textit{admits a logarithmic derivative along} $h_n$, i.e., there exists a measurable function $\beta_{h_n}^{\mu_n}\in L^2(E;\mu_n)$ such that
\[\int_E \pard{\phi}{h_n}d\mu_n=-\int_E \phi\beta_{h_n}^{\mu_n}d\mu_n\]
for every $\phi\in\FCE$, see~\cite[Section II.3.a)]{MR}.

Denote the ``minimal closed extension'' by $(\Escr_{\mu_n,0}^{h_n},\Dscr(\Escr_{\mu_n,0}^{h_n}))$, which is by definition the
closure of $(\Escr_{\mu_n,0}^{h_n},\FCE)$ in $L^2(E;\mu_n)$.
Furthermore, define a ``maximal'' partial form by
\begin{multline*}
\Dscr(\Escr_{\mu_n}^{h_n}):=\left\{u\in L^2(E;\mu_n)\;\bigg\vert\;\text{for}\;\nu_n\text{-a.e.}\;x\in E_0,\; s\mapsto u(x+sh_n)\;\text{has an absolutely}\right.\\
\left.\text{continuous}\;(ds)\text{-version}\;\tilde{u}_x\;\text{and}\;\pard{u}{h_n}:=\left(\frac{d\tilde{u}(x+sh_n)}{ds}\right)\in L^2(E;\mu_n)\right\},
\end{multline*}
and
\[\Escr_{\mu_n}^{h_n}(u,v):=\int_E \pard{u}{h_n}\pard{v}{h_n}d\mu_n,\quad u,v\in\Dscr(\Escr_{\mu_n}^{h_n}).\]
It follows from the local integrability of $s\mapsto(\rho_n(x+sh_n))^{-1}$ that $\Escr_{\mu_n}^{h_n}$ is closed, see~\cite{AR}.

\begin{lem}\label{minmaxpard}
Let $\mu$, $h$, $\nu$ and $\rho$ be as above for the limit measure $\mu$, so that $\mu$ satisfies Assumption~\ref{qinv} and $s\mapsto(\rho(x+sh))^{-1}$ is locally integrable for $\nu$-a.e. $x\in E_0$.
For the minimal and maximal partial forms of $\mu$ along $h$, defined as above, the following statement holds
\[(\Escr_{\mu,0}^{h},\Dscr(\Escr_{\mu,0}^{h}))=(\Escr_{\mu}^{h},\Dscr(\Escr_{\mu}^{h})).\]
\end{lem}
\begin{proof}
See~\cite[Lemma 3.3]{K2}.
\end{proof}

Now we formulate the main result for the Mosco convergence of partial Dirichlet forms.

\begin{thm}\label{partialmosco}
Let Assumptions~\ref{mnconv}--\ref{mutildebd} hold.
Suppose that there exist disintegrations $\mu_n=\rho_n(x+sh_n)ds\;\nu_n(dx)$ such that
\begin{enumerate}[(1)]
\item $\nu_n\to\nu$ weakly,
\item there exists an increasing sequence of numbers $\{n_i\}$, $n_i\uparrow\infty$ such that $\{\mathbbm{1}_{E_0^{n_i}}\tilde{\mu}_n^{h_n}\}$
is tight and, moreover,
\[\mathbbm{1}_{E_0^{n_i}}\tilde{\mu}_n^{h_n}\to\mathbbm{1}_{E_0^{n_i}}\tilde{\mu}^h\]
weakly for every $n_i$ as $n\to\infty$.
\end{enumerate}
Then $\Escr_{\mu_n}^{h_n}\to\Escr_\mu^h$ Mosco.
\end{thm}
\begin{proof}
See~\cite[Theorem 3.4]{K2}.
\end{proof}

We shall now deal with the gradient case. For this purpose we fix an orthonormal basis $\{e_i\;\vert\;i\in\N\}$ of the separable Hilbert space $X$.
Furthermore, we fix a sequence of vectors $h_n\in X$, $\lrnorm{h_n}_X=1$ weakly converging to $e_1$. We now construct a double sequence $\{e_n^i\}$ such
that
\begin{enumerate}[(1)]
\item for fixed $n$ every $\{e_n^i\;\vert\;i\in\N\}$ is an orthonormal basis of $X$,
\item for fixed $i$ every $\{e_n^i\}$ converges weakly to $e_i$.
\end{enumerate}
To do so, note first that $\norm{h_n}_X=\norm{e_1}_X=1$ and $h_n\to e_1$ weakly imply $\norm{h_n-e_1}_X\to0$. Let $R_n$ be the unitary operator on $X$ which is the identity on the orthogonal complement of $\operatorname{lin}\{e_1,h_n\}$ and maps $e_1$ to $h_n$, namely the rotation in the plane $\operatorname{lin}\{e_1,h_n\}$ if $h_n\notin\operatorname{lin}\{e_1\}$, and $R_n:=\pm\Id$ on $\operatorname{lin}\{e_1\}$ if $h_n=\pm e_1$. For $h_n=-e_1$, $R_n$ is a reflection rather than a rotation, which is sufficient. Then $\norm{R_n-\Id}_{\Lscr(X)}=\norm{h_n-e_1}_X\to0$, and
\[e_n^i:=R_ne_i,\quad i\in\N,\]
defines orthonormal bases of $X$ with $e_n^1=h_n$ and $\sup_i\norm{e_n^i-e_i}_X\to0$ as $n\to\infty$.

Recall that a sum of finitely many closed forms is closed, see~\cite[Proposition I.3.7]{MR}. A countable sum of non-negative closed forms is closed as well, since a countable sum of non-negative lower semicontinuous functionals is lower semicontinuous. We will consider a sequence of forms
\[\Escr_{\mu_n}=\sum\limits_{i=1}^\infty \Escr_{\mu_n}^{e_n^i}\]
with the domain of definition
\[\Dscr(\Escr_{\mu_n})=\left\{u\in\bigcap\limits_{i=1}^\infty \Dscr(\Escr_{\mu_n}^{e_n^i})\;\bigg\vert\;\sum_{i=1}^\infty\Escr_{\mu_n}^{e_n^i}(u)<\infty\right\},\]
which is the ``maximal extension'' of a gradient form. The ``minimal extension'' of $(\Escr_{\mu_n},\FCE)$, denoted
by $(\Escr_{\mu_n,0},\Dscr(\Escr_{\mu_n,0}))$, is the closure of $(\Escr_{\mu_n},\FCE)$ in
$L^2(E;\mu_n)$. Note that
$\FCE \subset \Dscr(\Escr_{\mu_n})$ since
\[\suml_{i=1}^\infty \lrbr{h,e_n^i}_X^2=\lrnorm{h}_X^{2}<\infty\]
for every $h\in X$, see~\cite[Proof of Proposition II.3.4]{MR}.

This type of form coincides with the infinite dimensional Dirichlet integral in the sense of the gradient introduced in \eqref{nabla}, since
by Parseval's identity
\begin{align*}
\Escr_{\mu_n}(u,v)=&\sum\limits_{i=1}^\infty \int_E \pard{u}{e_n^i}\pard{v}{e_n^i}d\mu_n
=\int_E \sum\limits_{i=1}^\infty \lrbr{\nabla_X u,e_n^i}_X \lrbr{e_n^i,\nabla_X v}_X d\mu_n\\
=&\int_E \lrbr{\nabla_X u,\nabla_X v}_X d\mu_n,\quad u,v\in \FCE.
\end{align*}

\begin{thm}\label{gradientmosco}
Suppose that, for every $i\in\N$, the hypotheses of Theorem~\ref{partialmosco}, including Assumptions~\ref{hnqinv} and \ref{mutildebd}, the local integrability of $s\mapsto(\rho(x+se_i))^{-1}$ for the limit measure $\mu$, and conditions (1)--(2), hold with $h_n:=e_n^i$, $h:=e_i$ and a suitable $l=l_i\in E'$ with $\lrbr{l_i,e_n^i}_X\neq0$ for all $n$ and $\lrbr{l_i,e_i}_X\neq0$. Suppose that
$(\Escr_{\mu,0},\Dscr(\Escr_{\mu,0}))=(\Escr_{\mu},\Dscr(\Escr_{\mu}))$, i.e., $\FCE$ is dense in
$\Dscr(\Escr_\mu)$ w.r.t. the norm $(\Escr_{\mu})_1^{1/2}$. Then $\Escr_{\mu_n}\to\Escr_\mu$ Mosco.
\end{thm}
\begin{proof}
Condition (M1) follows from the fact that (M1) is fulfilled for every sequence of partial forms $\{\Escr_{\mu_n}^{e_n^i}\}$, by Fatou's lemma for the sum over $i$. Let
us verify (M2).
Since $\FCE$ is dense in $(\Dscr(\Escr_\mu),(\Escr_\mu)_1^{1/2})$, Lemma~\ref{densemosco} implies
that it suffices to show that $\Escr_{\mu_n}(f)\to\Escr_\mu(f)$ for every $f\in\FCE$. It is clear that
\begin{equation}\label{base}
\supl_n\suml_{i=1}^\infty \lrbr{h,e_n^i}_X^2=\supl_n\lrnorm{h}_X^{2}=\lrnorm{h}_X^{2}<\infty
\end{equation}
for all $h\in X$. Take $f=\phi(l_1,\ldots,l_d)$, $\phi\in C_0^\infty (\R^d)$, $l_1,\ldots,l_d\in E'$. Then
\[\Escr_{\mu_n}(f)=\suml_{i=1}^\infty\int_E\suml_{j_1,j_2=1}^d\pard{\phi}{x_{j_1}}(l_1,\ldots,l_d)
\pard{\phi}{x_{j_2}}(l_1,\ldots,l_d)\lrbr{l_{j_1},e_n^i}_X \lrbr{l_{j_2},e_n^i}_X d\mu_n.\]
The claim follows from the Cauchy--Schwarz inequality, weak convergence $\mu_n\to\mu$ and \eqref{base}. More directly, by Parseval's identity, $\Escr_{\mu_n}(f)=\int_E\norm{\nabla_X f}_X^2d\mu_n$ does not depend on the choice of the orthonormal basis, and \[z\mapsto\norm{\nabla_X f(z)}_X^2=\sum_{j_1,j_2=1}^d\pard{\phi}{x_{j_1}}\pard{\phi}{x_{j_2}}(l_1(z),\ldots,l_d(z))\lrbr{l_{j_1},l_{j_2}}_X\] is bounded and continuous, so that $\Escr_{\mu_n}(f)\to\Escr_\mu(f)$ by weak convergence of $\mu_n\to\mu$.
\end{proof}

\subsubsection{Logarithmic derivatives}\label{sss:logder}

This subsection is taken from~\cite{K2} with slight changes. Its result will be used to prove convergence in the non-symmetric case.

Now consider the same situation as above, but assume only that $\mu_n\to\mu$ weakly. It is clear by Prokhorov's theorem, see~\cite[Theorem 5.2]{Bil}, that $\{\mu_n\}$ is tight.
Recall the well-known fact that
$L^2$-convergence of the logarithmic derivatives of measures implies strong convergence of the corresponding semigroups. We fix some $h\in X$
such that every $\mu_n$ has the logarithmic derivative $\beta_h^{\mu_n}\in L^2(E;\mu_n)$ along $h$ and consider the sequence of partial
forms $\{\Escr_{\mu_n}^h\}$ defined by
\[\Escr_{\mu_n}^h(f,g)=\int_E \pard{f}{h}\pard{g}{h}d\mu_n\]
for $f,g\in\FCE$. The condition $\beta_h^{\mu_n}\in L^2(E;\mu_n)$ implies the closability of these forms. In this case,
$h$ is called \textit{well-$\mu_n$-admissible}, see~\cite[Section II.3.a)]{MR}. As usual, the
closure of $(\Escr_{\mu_n}^h,\FCE)$ is considered, which coincides with the maximal extension by~\cite{RZ3}. It was proved in~\cite{RZ3} that $\FCE$ is dense
in $(\Dscr(\Escr_\mu^h),(\Escr_\mu^h)_1^{1/2})$ for every partial form $\Escr_\mu^h$ if $\mu$ admits a logarithmic derivative along $h$.

\begin{prop}\label{weakbeta}
Let $\sup_n\lrnorm{\beta_h^{\mu_n}}_{L^2(E;\mu_n)}<\infty$. Then $\mu$ possesses a logarithmic derivative along $h$ and $\{\Escr_{\mu_n}^h\}$
$\Gamma$-converges to $\Escr_\mu^h$. If, in addition, $\lrnorm{\beta_h^{\mu_n}}_{L^2(E;\mu_n)}\to\lrnorm{\beta_h^\mu}_{L^2(E;\mu)}$, where $\beta_h^\mu$ denotes the logarithmic derivative of $\mu$ along $h$, then
$\Escr_{\mu_n}^h\to\Escr_{\mu}^h$ Mosco.
\end{prop}
\begin{proof}
Conditions (M2) and (G2) can be verified as in Lemma~\ref{densemosco}. Let us verify condition (G1). Extract from $\{\beta_h^{\mu_n}\}$ an $\Hscr$-weakly
convergent subsequence, denoted in the following again by $\{\beta_h^{\mu_n}\}$, such that
$\beta_h^{\mu_n}\to\beta_h\in L^2(E;\mu)$ $\Hscr$-weakly. Then by the properties of weak convergence in $\Hscr$
\[\int_E \phi\beta_h\;d\mu=\liml_n \int_E \phi \beta_h^{\mu_n} d\mu_n=-\liml_n \int_E \pard{\phi}{h} d\mu_n=-\int_E \pard{\phi}{h} d\mu\]
for every $\phi\in\FCE$. Hence $\mu$ has the logarithmic derivative $\beta_h^\mu:=\beta_h\in L^2(E;\mu)$ and, moreover,
$\beta_h^{\mu_n}\to\beta_h^\mu$ $\Hscr$-weakly along the full sequence, since the weak limit is uniquely determined, see Corollary~\ref{exstrongconv} (iv). Now let $f_n\to f$ strongly in $\Hscr$ with $\sup_n\Escr^h_{\mu_n}(f_n)<\infty$, since
otherwise (G1) is trivial. W.l.o.g. $f_n\in\Dscr(\Escr_{\mu_n}^h)$. Let $K\subset E$ be a compact set. Obviously, the tightness of the measures $\{\mu_n\}$
and the Cauchy--Schwarz inequality
\[\left(\int_{E\setminus K}\lrabs{\pard{f_n}{h}}d\mu_n\right)^2\le \mu_n(E\setminus K)\int_E\left(\pard{f_n}{h}\right)^2 d\mu_n\]
imply that the sequence of measures $\displaystyle{\{\nu_n\}=\left\{{\pard{f_n}{h}}\mu_n\right\}}$ is tight. Passing to a subsequence, we may assume that $\Escr^h_{\mu_n}(f_n)$ converges to $\varliminf_n\Escr^h_{\mu_n}(f_n)$. Extract a weakly convergent
subsequence $\nu_n\to \nu$, denoted again by $\{\nu_n\}$. Then
\begin{multline*}
\int_E \phi\;d\nu=\liml_n\int_E \phi\pard{f_n}{h}d\mu_n\\
=\lim_n\left[-\int_E\pard{\phi}{h} f_n d\mu_n-\int_E \phi f_n \beta_h^{\mu_n} d\mu_n\right]\\
=-\int_E \pard{\phi}{h} f\;d\mu -\int_E \phi f \beta_h^\mu d\mu,
\end{multline*}
$\phi\in\FCE$. Moreover, by the Cauchy--Schwarz inequality, for $\phi\in\FCE$,
\[\lrabs{\int_E \phi\,d\nu}=\lim_n \lrabs{\int_E \phi \frac{\partial f_n}{\partial h}d\mu_n}
\le\left(\lim_n\int_E\left(\frac{\partial f_n}{\partial h}\right)^2 d\mu_n\right)^{1/2}\cdot\left(\int_E \phi^2 d\mu\right)^{1/2}.\]
Since $\FCE$ is dense in $L^2(E;\mu)$, the functional $\phi\mapsto\int_E\phi\,d\nu$ extends to a bounded linear functional on $L^2(E;\mu)$, so that $\nu=g\mu$ with some $g\in L^2(E;\mu)$ and
\[\int_E g^2d\mu\le\lim_n\int_E\left(\frac{\partial f_n}{\partial h}\right)^2 d\mu_n=\varliminf_n\Escr^h_{\mu_n}(f_n).\]
By the identity above, $g$ is the weak derivative of $f$ along $h$. Hence $f\in\Dscr(\Escr^h_\mu)$ and $\Escr^h_\mu(f)\le\varliminf_n\Escr^h_{\mu_n}(f_n)$, which proves (G1).

It is easily seen from the proof that the stronger assumption $\lrnorm{\beta_h^{\mu_n}}_{L^2(E;\mu_n)}\to\lrnorm{\beta_h^\mu}_{L^2(E;\mu)}$
implies by Lemma~\ref{dualconv} (2) that $\beta_h^{\mu_n}\to\beta_h^\mu$ $\Hscr$-strongly and (M1) is fulfilled so that $\Escr_{\mu_n}^h\to\Escr_\mu^h$ Mosco.
\end{proof}

\begin{rem}\label{weakderi}
Theorem~\ref{gradientmosco} also holds in the case of one common orthonormal basis $\{e_i\}$ if we assume
that each $\mu_n$, $n\in\N$, $\mu$ admits a logarithmic derivative along each $e_i$, $i\in\N$ denoted by $\beta_i^{\mu_n}$
and that $\beta_i^{\mu_n}\to\beta_i^\mu$ $\Hscr$-strongly, and if, as in Theorem~\ref{gradientmosco}, the minimal and maximal extensions of the limiting gradient form coincide, i.e., $(\Escr_{\mu,0},\Dscr(\Escr_{\mu,0}))=(\Escr_{\mu},\Dscr(\Escr_{\mu}))$. The latter is needed for (M2) via Lemma~\ref{densemosco} and does not follow from the existence of the logarithmic derivatives along the single directions $e_i$. The proof is essentially the same as above, see also \textup{\cite[Proof of Proposition II.3.4]{MR}}. Condition (M1) follows from (M1) for the partial forms by Fatou's lemma for the sum over $i$.
\end{rem}

\section{Examples of general convergence}\label{sec:genconv}

\subsection{Finite dimensional $a_{i,j}$-forms}\label{ss:finaij}

Let $\sigma_n dx=\mu_n$, $n\in\N$, and $\sigma dx=\mu$ be as in Subsection~\ref{ss:finsym}. Let Assumption~\ref{speed} hold.

Let $d\ge 3$. For $n\in\N$ let $a_{i,j}^n,b_i^n,d_i^n,c^n\in L^1_\loc(\R^d;dx)$, $1\le i,j\le d$ and define for
$u,v\in C_0^\infty(\R^d)$
\begin{align*}
\Escr^n(u,v):=&\sum_{i,j=1}^d\int a_{i,j}^n\pard{u}{x_i}\pard{v}{x_j}dx+\sum_{i=1}^d\int b_i^n \pard{u}{x_i}v\;dx\\
&+\sum_{i=1}^d\int d_i^n u\pard{v}{x_i}dx+\int c^n uv\;dx.
\end{align*}
Then for each $n$ $(\Escr^n,C_0^\infty(\R^d))$ is a densely defined bilinear form on $L^2(\R^d;\sigma_n dx)$.
Set $\tilde{a}_{i,j}^n:=\frac{1}{2}(a_{i,j}^n+a_{j,i}^n)$, $\check{a}_{i,j}^n:=\frac{1}{2}(a_{i,j}^n-a_{j,i}^n)$,
$b^n:=(b_1^n,\ldots,b_d^n)$ and $d^n:=(d_1^n,\ldots,d_d^n)$.

\begin{thm}\label{closability}
Suppose that
\begin{enumerate}[(i)]
\item we have
\[\sum_{i,j=1}^d \tilde{a}_{i,j}^n\xi_i\xi_j\ge \abs{\xi}^2\]
$dx$-a.e. for all $\xi=(\xi_1,\ldots,\xi_d)\in\R^d$ and for every $n$,
\item for some constant $M_n>0$ we have
$\abs{\check{a}_{i,j}^n}\le M_n$ for all $1\le i,j\le d$ and for every $n$,
\item $\abs{b^n+d^n}\in L^d_\loc(\R^d;dx)$, $c^n\in L^{d/2}_\loc(\R^d;dx)$ for every $n$,
\item $\abs{b^n-d^n}\in L^d(\R^d;dx)$ or $\frac{\abs{b^n-d^n}}{\sqrt{\sigma_n}}\in L^\infty(\R^d;dx)$ for every $n$,
\item for some $\a_0^n>0$, $(c^n+\a_0^n\sigma_n)dx-\sum_{i=1}^d\pard{d_i^n}{x_i}$ is a positive measure
on $\Bscr(\R^d)$ for every $n$,
\item $b^n=\beta^n+\gamma^n$ such that $\abs{\beta^n},\abs{\gamma^n}\in L^1_\loc(\R^d;dx)$ and
$(c^n+\a_0^n\sigma_n)dx-\sum_{i=1}^d\pard{\gamma_i^n}{x_i}$ is a positive measure on $\Bscr(\R^d)$ and
$\abs{\beta^n}\in L^d(\R^d;dx)$ or $\frac{\beta^n}{\sqrt{\sigma_n}}\in L^\infty(\R^d;dx)$ for every $n$.
\end{enumerate}
Then for every $n\in\N$ there exists $\a_n\in (0,\infty)$ such that $(\Escr_{\a_n}^n,C_0^\infty(\R^d))$ is closable on
$L^2(\R^d;\sigma_n dx)$ and its closure $(\Escr_{\a_n}^n,\Dscr(\Escr_{\a_n}^n))$ is a semi-Dirichlet form satisfying the weak sector condition with some sector constant $K_n\ge 1$.
\end{thm}
\begin{proof}
Apply~\cite[Theorem 1.2]{RS2} to each $n$ separately.
\end{proof}

Theorem~\ref{closability} is a statement about each single form. The numbers $\a_n$ and $K_n$ depend on the constants in (ii) and (v) and on the coefficients in (iii), (iv) and (vi) through the relative form bounds used in~\cite{RS2}. Membership of the coefficients in the respective function spaces for every $n$ does not provide a common lower bound for the forms $\Escr^n$, as the following example shows.

\begin{expl}\label{nounifshift}
Let $d\ge 3$, $\sigma_n\equiv 1$, $a_{i,j}^n=\delta_{i,j}$, $c^n=d^n=0$, $\gamma^n=0$, and choose $b\in C_0^\infty(\R^d;\R^d)$ such that $\div b\ge c_0>0$ on some ball $B$. Let $x_n\in\R^d$ with $\abs{x_n}\to\infty$ and set $b^n(x):=\beta^n(x):=nb(x-x_n)$. For each fixed $n$, the conditions (i)--(vi) of Theorem~\ref{closability} are satisfied, and $b^n\to 0$ in $L^2_\loc(\R^d;dx)$, since the supports move to infinity. Let $0\neq u\in C_0^\infty(\R^d)$ with $\supp u\subset B$ and set $u_n(x):=u(x-x_n)$. Integration by parts yields
\[\Escr^n(u_n,u_n)=\int_{\R^d}\abs{\nabla u}^2dx-\frac{n}{2}\int_{\R^d}(\div b)u^2dx,\]
so that, for every fixed $\a>0$, $\Escr^n_\a(u_n,u_n)<0$ for all sufficiently large $n$. Hence there is no $\a>0$ such that $\Escr^n_\a(v,v)\ge 0$ for all $v\in C_0^\infty(\R^d)$ and all $n$.
\end{expl}

We are grateful to one of the referees for this example. It shows that a uniform lower bound has to be assumed explicitly. Recall that the framework of Subsection~\ref{ss:nonsym} requires a common lower bound for the forms, see Remark~\ref{lambdalealpha}.

\begin{assumption}[Uniform lower bound]\label{unifshift}
The numbers $\a_n$ in Theorem~\ref{closability} can be chosen such that $\sup_n\a_n<\infty$.
\end{assumption}

If Assumption~\ref{unifshift} holds and $\a\ge\sup_n\a_n$, then for every $n$ the form $(\Escr^n_\a,C_0^\infty(\R^d))$ is closable on $L^2(\R^d;\sigma_n dx)$ and its closure $(\Escr^n_\a,\Dscr(\Escr^n_\a))$ is a semi-Dirichlet form satisfying the weak sector condition. Indeed, $\Escr^n_\a=\Escr^n_{\a_n}+(\a-\a_n)\deli{\cdot}{\cdot}_{H_n}$ differs from $\Escr^n_{\a_n}$ by a bounded, non-negative and symmetric form, which preserves closability, the semi-Dirichlet property and the non-negativity of the symmetric part. Moreover, since $\Escr^n_{\a_n+1}(u,u)\le\Escr^n_{\a+1}(u,u)$ and $\norm{u}_{H_n}^2\le\Escr^n_{\a+1}(u,u)$, the weak sector condition holds with a sector constant $K_n+\a-\a_n$. Assumption~\ref{unifshift} is satisfied, for example, if the constants $M_n$, $\a_0^n$ and the quantities controlling the relative form bounds of the coefficients in the proof of~\cite[Theorem 1.2]{RS2} are bounded uniformly in $n$. Note that only a common lower bound is needed in the sequel. The sector constants $K_n$ do not have to be uniformly bounded, see the proof of Theorem~\ref{convergenceaij} below.

Let $a_{i,j},b_i,d_i,c\in L^1_\loc(\R^d;dx)$, $1\le i,j\le d$ satisfy the conditions of Theorem~\ref{closability} for a single form with some $\a_\infty\in(0,\infty)$ in place of $\a_n$. For $u,v\in C_0^\infty(\R^d)$ define
\begin{align*}
\Escr(u,v):=&\sum_{i,j=1}^d\int a_{i,j}\pard{u}{x_i}\pard{v}{x_j}dx+\sum_{i=1}^d\int b_i \pard{u}{x_i}v\;dx\\
&+\sum_{i=1}^d\int d_i u\pard{v}{x_i}dx+\int cuv\;dx.
\end{align*}
Let Assumption~\ref{unifshift} hold and fix $\a\ge\max\{\sup_n\a_n,\a_\infty\}$. By Theorem~\ref{closability} and the preceding discussion, $(\Escr_\a,C_0^\infty(\R^d))$ is closable on $L^2(\R^d;\sigma dx)$. Its closure we denote by $(\Escr_\a,\Dscr(\Escr_\a))$.
The following theorem generalizes the results of R\"ockner and Zhang~\cite{RZ2} to the case of changing speed measures.
\begin{thm}\label{convergenceaij}
Let Assumption~\ref{unifshift} hold and let $\a\ge\max\{\sup_n\a_n,\a_\infty\}$ be as above. Assume that, for each $1\le i,j\le d$,
\[\frac{a_{i,j}^n}{\sqrt{\sigma_n}}\to \frac{a_{i,j}}{\sqrt{\sigma}},\quad\pard{a_{i,j}^n}{x_j}(\sqrt{\sigma_n})^{-1}\to\pard{a_{i,j}}{x_j}(\sqrt{\sigma})^{-1}\quad\text{strongly in $L^2_\loc(\R^d;dx)$},\]
\[\frac{b_i^n}{\sqrt{\sigma_n}}\to \frac{b_i}{\sqrt{\sigma}}\quad\text{strongly in $L^2_\loc(\R^d;dx)$},\]
\[\frac{d_i^n}{\sqrt{\sigma_n}}\to \frac{d_i}{\sqrt{\sigma}},\quad\pard{d_i^n}{x_i}(\sqrt{\sigma_n})^{-1}\to\pard{d_i}{x_i}(\sqrt{\sigma})^{-1}\quad\text{strongly in $L^2_\loc(\R^d;dx)$},\]
\[\frac{c^n}{\sqrt{\sigma_n}}\to \frac{c}{\sqrt{\sigma}}\quad\text{strongly in $L^2_\loc(\R^d;dx)$},\]
where the derivatives are taken in the distributional sense, and assume that, for every $n\in\N$, all the quotients $a_{i,j}^n/\sqrt{\sigma_n}$, $\pard{a_{i,j}^n}{x_j}(\sqrt{\sigma_n})^{-1}$, $b_i^n/\sqrt{\sigma_n}$, $d_i^n/\sqrt{\sigma_n}$, $\pard{d_i^n}{x_i}(\sqrt{\sigma_n})^{-1}$, $c^n/\sqrt{\sigma_n}$ belong to $L^2_\loc(\R^d;dx)$, and likewise for the limiting coefficients.

Assume that $\sqrt{\sigma_n}\to\sqrt{\sigma}$ weakly in $L^2_\loc(\R^d;dx)$.

Note that $C_0^\infty(\R^d)$ is dense in $\Dscr(\Escr_\a)$ w.r.t. $\tilde{\Escr}_{\a+1}^{1/2}$, since $\Dscr(\Escr_\a)$ was defined as the closure of $C_0^\infty(\R^d)$.

Finally, assume that the symmetric parts satisfy condition (M1), i.e., for every $\Hscr$-weakly convergent sequence $u_n\to u$,
\[\tilde{\Escr}_{\a}(u)\le\varliminf_n\tilde{\Escr}^n_{\a}(u_n),\]
where $\tilde{\Escr}^n_\a(u_n):=+\infty$ if $u_n\notin\Dscr(\Escr^n_\a)$ and $\tilde{\Escr}_\a(u):=+\infty$ if $u\notin\Dscr(\Escr_\a)$.

Then $\Escr^n_\a\to\Escr_\a$ in the generalized sense, more precisely, conditions (F1b) and (F2) of Definition~\ref{Fs} hold.
\end{thm}

Before we prove Theorem~\ref{convergenceaij}, we prove two lemmas. The first provides a useful observation about
``local'' strong and weak convergence, and the second explains the underlying structure of approximation of our forms.

\begin{lem}\label{phi-int}
Let $\phi\in C_0^\infty(\R^d)$, $u_n\in L^2(\R^d;\sigma_n dx)$, $n\in\N$, $u\in L^2(\R^d;\sigma dx)$. If $u_n\to u$ strongly in $\Hscr$,
then $\phi u_n\to \phi u$ strongly in $\Hscr$, too. The same holds for weak convergence.
\end{lem}
\begin{proof}
Note that clearly $\phi\in L^2(\R^d;\sigma_n dx)$ for any $n$ and $\phi\in L^2(\R^d;\sigma dx)$, with corresponding
$dx$-classes respected, if $\phi\in C_0^\infty(\R^d)$.
Let $\{\phi_m\}\subset C_0^\infty(\R^d)$ such that $\phi_m\to u$ in $L^2(\R^d;\sigma dx)$ and
\[\lim_m\varlimsup_n\norm{\phi_m-u_n}_{L^2(\R^d;\sigma_n dx)}=0.\]
Then by H\"older's inequality
\[\norm{\phi_m\phi-u\phi}_{L^2(\R^d;\sigma dx)}\le\norm{\phi}_{L^\infty(\R^d;dx)}\norm{\phi_m-u}_{L^2(\R^d;\sigma dx)}\to 0\]
as $m\to\infty$ and
\[\norm{\phi_m\phi-u_n\phi}_{L^2(\R^d;\sigma_n dx)}\le\norm{\phi}_{L^\infty(\R^d;dx)}\norm{\phi_m-u_n}_{L^2(\R^d;\sigma_n dx)}\to 0\]
as $m,n\to\infty$. The case of weak convergence follows from Lemma~\ref{weakmosco} since $\phi\cdot\psi\in C_0^\infty(\R^d)$ for any
$\phi,\psi\in C_0^\infty(\R^d)$.
\end{proof}

\begin{lem}\label{VtoV}
Consider $(\Escr_\a,\Dscr(\Escr_\a))$ as above. Suppose that the
conditions of Theorem~\ref{convergenceaij} hold. Then
$(L^2(\R^d;\sigma_n dx),\tilde{\Escr}^n_\a)\to(L^2(\R^d;\sigma dx),\tilde{\Escr}_\a)$ in the sense of Definition~\ref{Dconv}.
\end{lem}
\begin{proof}
Condition (1) is trivial, and (2) follows from the assumption. For (3), observe that for $\phi\in C_0^\infty(\R^d)$, $\phi\sqrt{\sigma_n}\to\phi\sqrt{\sigma}$ weakly in $L^2_\loc(\R^d;dx)$
by Lemma~\ref{sigmastrongweak} and the assumption that $\sqrt{\sigma_n}\to\sqrt{\sigma}$ weakly in $L^2_\loc(\R^d;dx)$. Also note that
$\frac{\tilde{a}_{i,j}^n}{\sqrt{\sigma_n}}\to\frac{\tilde{a}_{i,j}}{\sqrt{\sigma}}$ strongly in $L^2_\loc(\R^d;dx)$ for $1\le i,j\le d$ by linearity.
Then clearly, by the product rule,
\begin{align*}
\tilde{\Escr}^n_\a(\phi,\phi)
=&\sum_{i,j=1}^d\int \tilde{a}_{i,j}^n \pard{\phi}{x_i}\pard{\phi}{x_j}dx
+\frac{1}{2}\sum_{i=1}^d\int b_i^n\pard{\phi^2}{x_i}dx+\frac{1}{2}\sum_{i=1}^d\int d_i^n\pard{\phi^2}{x_i}dx\\
&+\int c^n \phi^2 dx+\a\int\phi^2\sigma_n dx\\
=&\sum_{i,j=1}^d\int \frac{\tilde{a}_{i,j}^n}{\sqrt{\sigma_n}} \pard{\phi}{x_i}\pard{\phi}{x_j}\sqrt{\sigma_n}dx
+\frac{1}{2}\sum_{i=1}^d\int \frac{b_i^n}{\sqrt{\sigma_n}}\pard{\phi^2}{x_i}\sqrt{\sigma_n}dx\\
&+\frac{1}{2}\sum_{i=1}^d\int \frac{d_i^n}{\sqrt{\sigma_n}}\pard{\phi^2}{x_i}\sqrt{\sigma_n}dx
+\int \frac{c^n}{\sqrt{\sigma_n}} \phi^2 \sqrt{\sigma_n}dx+\a\int\phi^2\sigma_n dx\\
\underset{n\to\infty}{\longrightarrow}&\sum_{i,j=1}^d\int \frac{\tilde{a}_{i,j}}{\sqrt{\sigma}} \pard{\phi}{x_i}\pard{\phi}{x_j}\sqrt{\sigma}dx
+\frac{1}{2}\sum_{i=1}^d\int \frac{b_i}{\sqrt{\sigma}}\pard{\phi^2}{x_i}\sqrt{\sigma}dx\\
&+\frac{1}{2}\sum_{i=1}^d\int \frac{d_i}{\sqrt{\sigma}}\pard{\phi^2}{x_i}\sqrt{\sigma}dx+\int \frac{c}{\sqrt{\sigma}} \phi^2 \sqrt{\sigma}dx
+\a\int\phi^2\sigma dx\\
=&\tilde{\Escr}_\a(\phi,\phi).
\end{align*}
The proof is complete.
\end{proof}

\begin{proof}[Proof of Theorem~\ref{convergenceaij}]
By Assumption~\ref{unifshift} and the choice of $\a$, the forms $\Escr^n_\a$, $n\in\N$, and $\Escr_\a$ are coercive closed forms with sector constants $K_n$ that need not be uniformly bounded, so that Theorem~\ref{Hino} applies to them. Since (F1b) implies (F1) by Lemma~\ref{PsiEquiv} (ii), it suffices to verify (F1b) and (F2). Let us first prove (F1b).
Let $u_n\in H_n$, $u\in H$, $u_n\to u$ $\Hscr$-weakly with $\varliminf_n(\tilde{\Escr}_{\a+1}^n)^{1/2}(u_n)<\infty$. By the assumed condition (M1) for the symmetric parts,
\[\tilde{\Escr}_\a(u)\le\varliminf_n\tilde{\Escr}^n_\a(u_n)\le\varliminf_n\tilde{\Escr}^n_{\a+1}(u_n)<\infty,\]
hence $u\in\Dscr(\Escr_\a)$, which is (F1b). Together with Lemmas \ref{VtoV} and \ref{densemosco}, we even obtain that $\tilde{\Escr}^n_\a\to\tilde{\Escr}_\a$ Mosco.

Next we prove (F2). Let $\phi\in C_0^\infty(\R^d)$. We choose $\Phi_n(\phi)=\phi\to\phi$ as the desired $\Hscr$-strongly convergent sequence.
Then, if $u_n\in\Dscr(\Escr^n_\a)$, $u\in\Dscr(\Escr_\a)$, $u_n\to u$ $\Hscr$-weakly, where no bound on the energies of the $u_n$ is needed,
we have by the assumption that $\sqrt{\sigma_n}\to\sqrt{\sigma}$ weakly in $L^2_\loc(\R^d;dx)$ and Lemma~\ref{sigmastrongweak} that
$u_n\sqrt{\sigma_n}\to u\sqrt{\sigma}$ weakly in $L^2(\R^d;dx)$ and, locally, via partial integration\footnote{The following identities hold for $u_n\in C_0^\infty(\R^d)$. The identity between the first expression $\Escr^n_\a(\phi,u_n)$ and the third expression, which contains no derivatives of $u_n$, extends to $u_n\in\Dscr(\Escr^n_\a)$ by density. Indeed, $\Escr^n_\a(\phi,\cdot)$ is continuous with respect to $(\tilde{\Escr}^n_{\a+1})^{1/2}$ by the weak sector condition, whereas the third expression is a sum of terms of the form $\int g\psi\,u_n\sqrt{\sigma_n}\,dx$ with $g\in L^2_\loc(\R^d;dx)$ and bounded $\psi$ with compact support, hence bounded by $C\norm{u_n}_{L^2(\R^d;\sigma_ndx)}\le C(\tilde{\Escr}^n_{\a+1}(u_n))^{1/2}$. The same density argument identifies the limit expression with $\Escr_\a(\phi,u)$ for $u\in\Dscr(\Escr_\a)$.}:
\begin{align*}
\Escr^n_\a(\phi,u_n)
=&\sum_{i,j=1}^d\int a_{i,j}^n\pard{\phi}{x_i}\pard{u_n}{x_j}dx+\sum_{i=1}^d\int b_i^n\pard{\phi}{x_i}u_n dx\\
&+\sum_{i=1}^d\int d_i^n\phi\pard{u_n}{x_i}dx+\int c^n \phi u_n dx+\a\int\phi u_n\sigma_n dx\\
=&-\sum_{i,j=1}^d\left[\int \pard{a_{i,j}^n}{x_j}\pard{\phi}{x_i} u_n dx+\int a_{i,j}^n \frac{\partial^2\phi}{\partial x_i\partial x_j} u_n dx\right]
+\sum_{i=1}^d\int b_i^n\pard{\phi}{x_i}u_n dx\\
&-\sum_{i=1}^d\left[\int \pard{d_i^n}{x_i}\phi u_n dx+\int d_i^n\pard{\phi}{x_i}u_n dx\right]+\int c^n \phi u_n dx+\a\int\phi u_n\sigma_n dx\\
=&-\sum_{i,j=1}^d\left[\int \pard{a_{i,j}^n}{x_j}(\sqrt{\sigma_n})^{-1}\pard{\phi}{x_i} u_n \sqrt{\sigma_n}dx+\int \frac{a_{i,j}^n}{\sqrt{\sigma_n}}
\frac{\partial^2\phi}{\partial x_i\partial x_j} u_n \sqrt{\sigma_n}dx\right]\\
&+\sum_{i=1}^d\int \frac{b_i^n}{\sqrt{\sigma_n}}\pard{\phi}{x_i}u_n \sqrt{\sigma_n}dx\\
&-\sum_{i=1}^d\left[\int \pard{d_i^n}{x_i}(\sqrt{\sigma_n})^{-1}\phi u_n\sqrt{\sigma_n}dx+\int \frac{d_i^n}{\sqrt{\sigma_n}}
\pard{\phi}{x_i}u_n\sqrt{\sigma_n} dx\right]\\
&+\int \frac{c^n}{\sqrt{\sigma_n}} \phi u_n\sqrt{\sigma_n} dx+\a\int\phi u_n\sigma_n dx\\
\underset{n\to\infty}{\longrightarrow}&-\sum_{i,j=1}^d\left[\int \pard{a_{i,j}}{x_j}(\sqrt{\sigma})^{-1}\pard{\phi}{x_i} u \sqrt{\sigma}dx+\int \frac{a_{i,j}}{\sqrt{\sigma}}
\frac{\partial^2\phi}{\partial x_i\partial x_j} u \sqrt{\sigma}dx\right]\\
&+\sum_{i=1}^d\int \frac{b_i}{\sqrt{\sigma}}\pard{\phi}{x_i}u \sqrt{\sigma}dx\\
&-\sum_{i=1}^d\left[\int \pard{d_i}{x_i}(\sqrt{\sigma})^{-1}\phi u\sqrt{\sigma}dx+\int \frac{d_i}{\sqrt{\sigma}}
\pard{\phi}{x_i}u\sqrt{\sigma} dx\right]\\
&+\int \frac{c}{\sqrt{\sigma}} \phi u\sqrt{\sigma} dx+\a\int\phi u\;\sigma dx\\
=&-\sum_{i,j=1}^d\left[\int \pard{a_{i,j}}{x_j}\pard{\phi}{x_i} u\;dx+\int a_{i,j} \frac{\partial^2\phi}{\partial x_i\partial x_j} u\; dx\right]
+\sum_{i=1}^d\int b_i\pard{\phi}{x_i}u\;dx\\
&-\sum_{i=1}^d\left[\int \pard{d_i}{x_i}\phi u\;dx+\int d_i\pard{\phi}{x_i}u\;dx\right]+\int c \phi u\;dx+\a\int\phi u\;\sigma dx\\
=&\sum_{i,j=1}^d\int a_{i,j}\pard{\phi}{x_i}\pard{u}{x_j}dx+\sum_{i=1}^d\int b_i\pard{\phi}{x_i}u\;dx\\
&+\sum_{i=1}^d\int d_i\phi\pard{u}{x_i}dx\\
&+\int c \phi u\;dx+\a\int\phi u\;\sigma dx=\Escr_\a(\phi,u)
\end{align*}
The proof is complete.
\end{proof}

\begin{rem}\label{M1remark}
In~\cite{Toe1}, (F1b) was derived without assuming (M1) for the symmetric parts, by identifying the weak limit of $\{u_{n_k}\}$ in $\Hscr_{\tilde{\Escr}_\a}$ with its weak limit in $\Hscr$. This identification is not justified in general. By Proposition~\ref{weakweakmosco}, it is a sufficient condition for (M1) for the symmetric parts. Therefore, (M1) is now assumed explicitly. It holds, for instance, if weak convergence in $\Hscr_{\tilde{\Escr}_\a}$ implies weak convergence in $\Hscr$, see Proposition~\ref{weakweakmosco}. It also holds if $b^n+d^n=0$ and $c^n=0$ for all $n$ and the conditions of Theorem~\ref{Kolconvfinite} are satisfied for the symmetric coefficients $\tilde{a}^n_{i,j}$, that is, $\sup_n\int_\Omega\sigma_n^{-1}dx<\infty$, the local bounds on the inverse matrices and the vague convergence of the coefficients and of the inverse matrices hold, including the coincidence of the minimal and maximal extensions of the limiting form. In this case, $\tilde{\Escr}^n_\a$ is a symmetric $a_{i,j}$-form as in Subsection~\ref{ss:finsym}, shifted by $\a\deli{\cdot}{\cdot}_{H_n}$, and this shift preserves (M1) by Lemma~\ref{dualconv} (2).
\end{rem}

\begin{rem}\label{naturalness}
The assumptions of Theorem~\ref{convergenceaij} are sufficient conditions designed to fit the abstract Theorem~\ref{Hino}. We do not claim that they are optimal. The convergence of the coefficients divided by $\sqrt{\sigma_n}$ in $L^2_\loc(\R^d;dx)$ is natural in the present weighted setting, since it is exactly what is needed to pass to the limit in the pairings with $u_n\sqrt{\sigma_n}$, which converge weakly in $L^2(\R^d;dx)$ by Lemma~\ref{sigmastrongweak}. It is in the spirit of the conditions of Kolesnikov~\cite{K1} in the symmetric case.
\end{rem}

\subsection{Infinite dimensional $a_{i,j}$-forms}\label{ss:infaij}

Let $E,X,\{\mu_n\},\mu$ be as in Section~\ref{idsc}, i.e.,
\[E'\subset X'\equiv X\subset E\quad\text{densely and continuously}.\]

In this situation, we work with $\Fscr C_b^\infty(E)$
instead of $\FCE$. $\FCb$ is defined to be the linear space
\[\Fscr C_b^\infty(E):=\{f(l_1,\ldots,l_m)\;\vert\;m\in\N,\;f\in C_b^\infty(\R^m),\;l_1,\ldots,l_m\in E'\},\]
where $C_b^\infty(\R^m)$ denotes the space of all real-valued infinitely differentiable and bounded functions on
$\R^m$ with all partial derivatives bounded.
This does not affect the general setting, since this choice is standard. Note, however, that
$\FCE$ functions are needed for the results of A.V. Kolesnikov in~\cite{K2}. 
Let Assumptions~\ref{mnconv} and \ref{qinv} hold. Then $H_n:=L^2(E;\mu_n)\to L^2(E;\mu)=:H$ in the sense of convergent Hilbert spaces, where we set
$C:=\FCb$. As usual,
set $\Hscr:=\disjoint_n H_n\dot\cup H$.

Assume that there exists an orthonormal basis $\{e_i\;\vert\;i\in\N\}$ of $X$
whose elements are well-$\mu_n$-admissible for each $n\in\N$ in the sense of Subsection~\ref{sss:logder}. Denote by $\beta_i^{\mu_n}$ the logarithmic derivative of $\mu_n$
along $e_i$, $n\in\N,i\in\N$.

The next lemma is analogous to Lemma~\ref{phi-int}.

\begin{lem}\label{phiuninfty}
Let $\phi\in \Fscr C_b^\infty(E)$, $u_n\in L^2(E;\mu_n)$, $n\in\N$, $u\in L^2(E;\mu)$. If $u_n\to u$ strongly in $\Hscr$,
then $\phi u_n\to \phi u$ strongly in $\Hscr$, too. The same holds for weak convergence.
\end{lem}
\begin{proof}
Note that clearly $\phi\in L^2(E;\mu_n)$ for any $n$ and $\phi\in L^2(E;\mu)$, with corresponding
$L^2$-classes respected, if $\phi\in\Fscr C_b^\infty(E)$.
Let $\{\phi_m\}\subset \Fscr C_b^\infty(E)$ such that $\phi_m\to u$ in $L^2(E;\mu)$ and
\[\lim_m\varlimsup_n\norm{\phi_m-u_n}_{L^2(E;\mu_n)}=0.\]
Then by H\"older's inequality
\[\norm{\phi_m\phi-u\phi}_{L^2(E;\mu)}\le\norm{\phi}_{L^\infty(E;\mu)}\norm{\phi_m-u}_{L^2(E;\mu)}\to 0\]
as $m\to\infty$ and
\[\norm{\phi_m\phi-u_n\phi}_{L^2(E;\mu_n)}\le\norm{\phi}_{L^\infty(E;\mu_n)}\norm{\phi_m-u_n}_{L^2(E;\mu_n)}\to 0\]
as $m,n\to\infty$, since $\sup_n\norm{\phi}_{L^\infty(E;\mu_n)}=\sup_n {{\esssup}_{z\in E}^{(\mu_n)}}\abs{\phi(z)}<\infty$.
The case of weak convergence follows from Lemma~\ref{weakmosco} since $\phi\cdot\psi\in\Fscr C_b^\infty(E)$ for any
$\phi,\psi\in\Fscr C_b^\infty(E)$.
\end{proof}

Fix $n\in\N$. Let $A^{(n)}:E\to\Lscr(X)$ be such that $z\mapsto\lrbr{A^{(n)}(z)h_1,h_2}_X$ is $\Bscr(E)$-measurable for all $h_1,h_2\in X$, and set $a_{i,j}^{(n)}(z):=\lrbr{A^{(n)}(z)e_j,e_i}_X$, $i,j\in\N$, so that
\[A^{(n)}(z)h=\sum_{i=1}^\infty\Big(\sum_{j=1}^\infty a_{i,j}^{(n)}(z)\lrbr{h,e_j}_X\Big) e_i,\quad h\in X,\]
with convergence in $X$.\footnote{In infinite dimensions, a summability condition such as $\sum_{i,j}\abs{a_{i,j}^{(n)}(z)}<\infty$, which would make $A^{(n)}(z)$ a trace class operator, is incompatible with the uniform ellipticity \eqref{linopA} below, since a trace class operator satisfies $\lrbr{A^{(n)}(z)e_i,e_i}_X\to0$ as $i\to\infty$. Such a condition is therefore not assumed.} Double sums $\sum_{i,j=1}^\infty$ are understood as limits of $\sum_{i,j=1}^N$ as $N\to\infty$. Furthermore,
assume that there exists a uniform ellipticity constant $c>0$ not depending on $n$ with
\begin{equation}\label{linopA}
\lrbr{A^{(n)}(z)h,h}_X=\sum_{i=1}^\infty\sum_{j=1}^\infty a_{i,j}^{(n)}(z)\lrbr{h,e_i}_X \lrbr{h,e_j}_X\ge c\norm{h}_X^2
\end{equation}
for all $h\in X$.

Set  $\tilde{A}^{(n)}:=\frac{1}{2}(A^{(n)}+\hat{A}^{(n)})$, $\check{A}^{(n)}:=\frac{1}{2}(A^{(n)}-\hat{A}^{(n)})$ where
$\hat{A}^{(n)}(z)$ denotes the adjoint of $A^{(n)}(z)$, $z\in E$. One easily observes that $\tilde{A}$ and $\check{A}$ can
be constructed as in \eqref{linopA} with $\tilde{a}_{i,j}^{(n)}:=\frac{1}{2}(a_{i,j}^{(n)}+a_{j,i}^{(n)})$ and
$\check{a}_{i,j}^{(n)}:=\frac{1}{2}(a_{i,j}^{(n)}-a_{j,i}^{(n)})$, respectively.
Assume that $\norm{\tilde{A}^{(n)}}_{\Lscr(X)}\in L^1(E;\mu_n)$,
$\norm{\check{A}^{(n)}}_{\Lscr(X)}\in L^\infty(E;\mu_n)$. In Theorem~\ref{infaij} below, this follows from assumption (T1).

Let $c^{(n)}\in L^\infty(E;\mu_n)$ and $b^{(n)},d^{(n)}\in L^\infty(E\to X;\mu_n)$ such that
\begin{align*}
&\int_E \left(\lrbr{b^{(n)},\nabla_X u}_X+c^{(n)} u\right)d\mu_n\ge 0,\\
&\int_E \left(\lrbr{d^{(n)},\nabla_X u}_X+c^{(n)} u\right)d\mu_n\ge 0
\end{align*}
for all $u\in\FCb$, $u\ge 0$.

Define for $n\in\N$, $u,v\in\FCb$
\begin{align*}
\Escr^n(u,v)=&\int_E\lrbr{\nabla_X u(z),A^{(n)}(z)\nabla_X v(z)}_X\mu_n(dz)\\
&+\int_E \lrbr{b^{(n)}(z),\nabla_X u(z)}_X v(z)\mu_n(dz)\\
&+\int_E u(z)\lrbr{d^{(n)}(z),\nabla_X v(z)}_X \mu_n(dz)\\
&+\int_E c^{(n)}(z)u(z)v(z)\mu_n(dz).
\end{align*}
Then for $n\in\N$ each $(\Escr^n,\FCb)$ is a densely defined positive bilinear form on $L^2(E;\mu_n)$ which is closable and
whose closure $(\Escr^n,\Dscr(\Escr^n))$ is a Dirichlet form. Each $\Escr^n$ satisfies the weak sector condition with a sector constant $\ge1$ which is bounded in terms of the ellipticity constant $c$ and of
\[\kappa_n:=\sup_{z\in E}\left[\norm{\check{A}^{(n)}(z)}_{\Lscr(X)}\vee\norm{d^{(n)}(z)-b^{(n)}(z)}_X\right]\]
by the results of~\cite[Section II.3.e)]{MR}, which also show that $\Escr^n$ has the semi-Dirichlet property. We denote this sector constant by $K_n$.

Define
\[Q^n(u,v):=\int_E\lrbr{\nabla_X u(z),\nabla_X v(z)}_X \mu_n(dz),\quad u,v\in\FCb.\]
Since each $e_i$, $i\in\N$ is well-$\mu_n$-admissible, $(Q^n,\FCb)$ is closable and its closure $(Q^n,\Dscr(Q^n))$ is
a symmetric Dirichlet form, see~\cite[Sections II.3.a) and II.3.b)]{MR}. Let us define $(Q,\Dscr(Q))$ in the same way.

For each $n\in\N$, $i,j\in\N$ set $b_i^{(n)}:=\lrbr{b^{(n)},e_i}_X$, $d_i^{(n)}:=\lrbr{d^{(n)},e_i}_X$ and assume that
\[a_{i,j}^{(n)},b_i^{(n)},d_i^{(n)},c^{(n)}\in L^2(E;\mu_n)\]
Again, in Theorem~\ref{infaij} below, this follows from (T1), since $\mu_n$ is a probability measure.

\begin{thm}\label{infaij}

Suppose that for $i,j\in\N$ there are given $a_{i,j},b_i,d_i,c\in L^2(E;\mu)$ such that
\begin{align*}
&a_{i,j}^{(n)}\to a_{i,j}\quad\Hscr\text{-strongly},\\
&b_{i}^{(n)}\to b_{i}\quad\Hscr\text{-strongly},\\
&d_{i}^{(n)}\to d_{i}\quad\Hscr\text{-strongly},\\
&c^{(n)}\to c\quad\Hscr\text{-strongly},
\end{align*}
for $i,j\in\N$ and that the limiting coefficients fulfill the same conditions as above.

Assume also that for $i\in\N$, $\sup_n\norm{\beta_i^{\mu_n}}_{L^2(E;\mu_n)}<\infty$, so that $\mu$ admits a logarithmic derivative $\beta_i^\mu$ along $e_i$ by Proposition~\ref{weakbeta}, and
\[\norm{\beta_i^{\mu_n}}_{L^2(E;\mu_n)}\to\norm{\beta_i^\mu}_{L^2(E;\mu)}.\]
Note that there exists a uniform constant $c_0>0$, namely $c_0:=\min\{c,1\}$, such that
\[c_0 Q_1^n(u,u)\le\tilde{\Escr}^n_1(u,u)\]
for every $u\in\Dscr(\Escr^n)\subset\Dscr(Q^n)$, $n\in\N$, and $c_0Q_1(u,u)\le\tilde{\Escr}_1(u,u)$ for every $u\in\Dscr(\Escr)\subset\Dscr(Q)$.

Indeed, for $u\in\FCb$ we have $u^2\in\FCb$, $u^2\ge0$, and hence, by the positivity conditions on $b^{(n)}$, $c^{(n)}$, $d^{(n)}$,
\begin{multline*}\int_E\lrbr{b^{(n)},\nabla_Xu}_Xu\,d\mu_n+\int_Eu\lrbr{d^{(n)},\nabla_Xu}_Xd\mu_n+\int_Ec^{(n)}u^2d\mu_n\\=\frac12\int_E\left(\lrbr{b^{(n)},\nabla_Xu^2}_X+c^{(n)}u^2\right)d\mu_n+\frac12\int_E\left(\lrbr{d^{(n)},\nabla_Xu^2}_X+c^{(n)}u^2\right)d\mu_n\ge0,\end{multline*}
so that $\Escr^n(u,u)\ge cQ^n(u,u)$ by \eqref{linopA}. This extends to the closures, and the same argument applies to the limiting form.

Assume the following uniform bounds and uniform tail conditions, where $P_N$ denotes the orthogonal projection of $X$ onto $\operatorname{lin}\{e_1,\ldots,e_N\}$:
\begin{enumerate}
\item[(T1)] $\sup_n\esssup_{z}\left[\norm{A^{(n)}(z)}_{\Lscr(X)}+\norm{b^{(n)}(z)}_X+\norm{d^{(n)}(z)}_X+\abs{c^{(n)}(z)}\right]<\infty$, with the essential suprema taken w.r.t. $\mu_n$, and the same bound holds for the limiting coefficients w.r.t. $\mu$;
\item[(T2)] for every $l\in E'$,
\[\lim_{N\to\infty}\sup_n\esssup_{z}\left[\norm{(\Id-P_N)\hat{A}^{(n)}(z)l}_X+\norm{(\Id-P_N)d^{(n)}(z)}_X\right]=0,\]
and the same holds for the limiting coefficients.
\end{enumerate}
Moreover, assume that the symmetric parts satisfy condition (M1), i.e., $\tilde{\Escr}(u)\le\varliminf_n\tilde{\Escr}^n(u_n)$ for every $\Hscr$-weakly convergent sequence $u_n\to u$, with the usual convention $+\infty$ outside the domains.

Assume that $(\Escr,\FCb)$ defined as follows
\begin{align*}
\Escr(u,v):=&\sum_{i,j=1}^\infty\int_E a_{i,j}\pard{u}{e_i}\pard{v}{e_j}d\mu+\sum_{i=1}^\infty\int_E b_i\pard{u}{e_i} v\;d\mu\\
&+\sum_{i=1}^\infty\int_E u d_i\pard{v}{e_i}d\mu+\int_E cuv\;d\mu,\quad u,v\in\FCb,
\end{align*}
is closable with closure $(\Escr,\Dscr(\Escr))$ and that $\FCb$ is dense in $\Dscr(\Escr)$ w.r.t. $\tilde{\Escr}_1^{1/2}$.

Then $\Escr^n\to\Escr$ in the generalized sense.
More precisely, (F1b) and (F2'b) hold. Moreover, $(\Escr,\Dscr(\Escr))$ is a semi-Dirichlet form, which follows from the assumptions on the limiting coefficients by~\cite[Section II.3.e)]{MR}.
\end{thm}
Note that, by Corollary~\ref{Knbdd} (i), the conditions (F1b) and (F2'b) imply the strong convergence of the resolvents also without the assumption $\sup_n K_n<\infty$, which, in Theorem~\ref{infaij}, follows from (T1) anyway. In contrast to Theorem~\ref{closability}, no shift is needed here, since the forms $\Escr^n$ are non-negative by the assumptions on $b^{(n)}$, $c^{(n)}$ and $d^{(n)}$.
For example, (T1), (T2) and the positivity conditions are satisfied if $A^{(n)}(z)=c\,\Id+K^{(n)}(z)$, where the operators $K^{(n)}(z)$ satisfy $\lrbr{K^{(n)}(z)h,h}_X\ge0$, have range in $\operatorname{lin}\{e_1,\ldots,e_M\}$ for some fixed $M\in\N$ and are bounded uniformly in $n$ and $z$, and if $b^{(n)}=d^{(n)}=0$ and $c^{(n)}\ge0$ are bounded uniformly in $n$. The conditions on the reference measures and their logarithmic derivatives hold trivially if $\mu_n=\mu$ for all $n$, where $\mu$ admits logarithmic derivatives in $L^2(E;\mu)$ along all $e_i$, for instance a suitable Gaussian measure. For genuinely changing reference measures, we refer to~\cite{K2}. Condition (M1) for the symmetric parts has to be verified separately, e.g. by means of Proposition~\ref{weakweakmosco}. Lemma~\ref{okydoky} below is not needed for the proof of Theorem~\ref{infaij}. Together with condition (M1) for the symmetric parts and Lemma~\ref{densemosco}, it shows that the symmetric parts $\tilde{\Escr}^n$ even converge to $\tilde{\Escr}$ in the sense of Mosco.

\begin{lem}\label{okydoky}
Let the conditions of Theorem~\ref{infaij} hold. Then
$(L^2(E;\mu_n),\tilde{\Escr}^n)\to(L^2(E;\mu),\tilde{\Escr})$ in the sense of Definition~\ref{Dconv}.
\end{lem}
\begin{proof}
Condition (1) is trivial, and (2) is a condition on $\Escr$. To prove (3), let $\phi\in\FCb$. Then clearly
\begin{align*}
\tilde{\Escr}^n(\phi,\phi)
=&\sum_{i,j=1}^\infty\int_E \tilde{a}_{i,j}^{(n)} \pard{\phi}{e_i}\pard{\phi}{e_j}d\mu_n
+\frac{1}{2}\sum_{i=1}^\infty\int_E b_i^{(n)}\pard{\phi^2}{e_i}d\mu_n\\
&+\frac{1}{2}\sum_{i=1}^\infty\int_E d_i^{(n)}\pard{\phi^2}{e_i}d\mu_n+\int_E c^{(n)} \phi^2 d\mu_n\\
\underset{n\to\infty}{\longrightarrow}&\sum_{i,j=1}^\infty\int_E \tilde{a}_{i,j} \pard{\phi}{e_i}\pard{\phi}{e_j}d\mu
+\frac{1}{2}\sum_{i=1}^\infty\int_E b_i\pard{\phi^2}{e_i}d\mu\\
&+\frac{1}{2}\sum_{i=1}^\infty\int_E d_i\pard{\phi^2}{e_i}d\mu+\int_E c \phi^2 d\mu\\
=&\tilde{\Escr}(\phi,\phi),
\end{align*}
where we have used Lemma~\ref{phiuninfty} and the product rule for G{\^a}teaux derivatives, and where the interchange of the limit and the infinite sums is justified by (T1) and (T2) as in the proof of Theorem~\ref{infaij} below.
Note that all sums, understood as limits over squares $\{i,j\le N\}$, converge by (T1) and (T2), cf.~\cite[Section II.3.a)]{MR}.
\end{proof}

\begin{proof}[Proof of Theorem~\ref{infaij}]
We first prove (F2'b). Note that by (T1) the sector constants $K_n$ are uniformly bounded. Set $\tilde{C}:=C=\FCb$. By assumption
$C\subset\Dscr(\Escr)$ densely w.r.t. $\tilde{\Escr}_1^{1/2}$. Let $n_k\in\N$, $n_k\uparrow \infty$, $n_{k+1}>n_k$. Let $u_k\in L^2(E;\mu_{n_k})$,
$u\in\Dscr(\Escr)$ such that $u_k\to u$ $\Hscr$-weakly and $\sup_k(\tilde{\Escr}_1^{n_k})^{1/2}(u_k)<\infty$. Without loss of generality, $u_k\in\Dscr(\Escr^{n_k})$.
By $c_0 Q_1^n(u,u)\le\Escr_1^n(u,u)$ for each $u\in\Dscr(\Escr^n)\subset\Dscr(Q^n)$ we have that $\sup_k(Q_1^{n_k})^{1/2}(u_k)<\infty$. Clearly
each $u_{k}$ admits a weak derivative in each direction $e_i$ denoted by $\pard{u_{k}}{e_i}$, and $\sup_k\sum_i\norm{\pard{u_k}{e_i}}^2_{L^2(E;\mu_{n_k})}=\sup_kQ^{n_k}(u_k)<\infty$. Hence, arguing as in the proof of Proposition~\ref{weakbeta} for each direction $e_i$, using that $\beta_i^{\mu_n}\to\beta_i^\mu$ $\Hscr$-strongly by Lemma~\ref{dualconv} (2), and by a diagonal argument,
we can extract a subsequence $\{u_{k_l}\}$, $u_{k_l}\in\Dscr(\Escr^{n_{k_l}})$ such that
$u_{k_l}\to u$ $\Hscr$-weakly and $u$ admits a weak derivative in each direction $e_i$ denoted by $\pard{u}{e_i}$ such that $\pard{u_{k_l}}{e_i}\to\pard{u}{e_i}$
weakly in $\Hscr$ for every $i$. Since $u\in\Dscr(\Escr)\subset\Dscr(Q)$, these are the weak derivatives of $u$ as an element of $\Dscr(Q)$. By Lemma~\ref{phiuninfty}, we have for $\phi\in\FCb$
\begin{align*}
\Escr^{n_{k_l}}(\phi,u_{k_l})
=&\sum_{i,j=1}^\infty\int_E {a}_{i,j}^{(n_{k_l})} \pard{\phi}{e_i}\pard{u_{k_l}}{e_j}d\mu_{n_{k_l}}
+\sum_{i=1}^\infty\int_E b_i^{(n_{k_l})}\pard{\phi}{e_i}u_{k_l}d\mu_{n_{k_l}}\\
&+\sum_{i=1}^\infty\int_E d_i^{(n_{k_l})}\phi\pard{u_{k_l}}{e_i}d\mu_{n_{k_l}}+\int_E c^{(n_{k_l})} \phi u_{k_l} d\mu_{n_{k_l}}\\
\underset{l\to\infty}{\longrightarrow}&\sum_{i,j=1}^\infty\int_E {a}_{i,j} \pard{\phi}{e_i}\pard{u}{e_j}d\mu
+\sum_{i=1}^\infty\int_E b_i\pard{\phi}{e_i}u\;d\mu\\
&+\sum_{i=1}^\infty\int_E d_i\phi\pard{u}{e_i}d\mu+\int_E c \phi u \;d\mu\\
=&\Escr(\phi,u).
\end{align*}
Here each single term with fixed indices converges by Lemma~\ref{phiuninfty} and the strong convergence of the coefficients, and the interchange of the limit with the infinite sums is justified as follows. For $\phi=f(l_1,\ldots,l_m)\in\FCb$ we have $\nabla_X\phi=\sum_{r=1}^m\partial_rf(l_1,\ldots,l_m)l_r$, hence $\sup_z\norm{(\Id-P_N)\nabla_X\phi(z)}_X\to0$ as $N\to\infty$. Writing the first three terms as $\int\lrbr{\nabla_X u_{k_l},\hat{A}^{(n_{k_l})}\nabla_X\phi}_X d\mu_{n_{k_l}}+\int\lrbr{b^{(n_{k_l})},\nabla_X\phi}_Xu_{k_l}d\mu_{n_{k_l}}+\int\phi\lrbr{d^{(n_{k_l})},\nabla_X u_{k_l}}_Xd\mu_{n_{k_l}}$, the parts of the sums with $i>N$ or $j>N$ are bounded by
\begin{multline*}C_\phi\esssup_z\Big[\norm{(\Id-P_N)\nabla_X\phi(z)}_X\\+\sum_{r=1}^m\norm{(\Id-P_N)\hat{A}^{(n_{k_l})}(z)l_r}_X+\norm{(\Id-P_N)d^{(n_{k_l})}(z)}_X\Big]\left(Q^{n_{k_l}}_1(u_{k_l})\right)^{1/2},\end{multline*}
with a constant $C_\phi$ depending on $\phi$ and on the bounds in (T1), by the Cauchy--Schwarz inequality. By (T2) and $\sup_lQ_1^{n_{k_l}}(u_{k_l})<\infty$, this tends to $0$ as $N\to\infty$, uniformly in $l$, where the essential supremum is taken w.r.t. $\mu_{n_{k_l}}$. The same holds for the limiting form, with the essential supremum w.r.t. $\mu$. Hence the limit can be taken termwise. Finally, the coordinate formula for $\Escr(\phi,u)$, $u\in\Dscr(\Escr)$, used above follows from the one for $u\in\FCb$ by density: the left-hand side is continuous in $u$ w.r.t. $\tilde{\Escr}_1^{1/2}$ by the weak sector condition, and the right-hand side is bounded by $C_\phi Q_1(u)^{1/2}\le C_\phi c_0^{-1/2}\tilde{\Escr}_1(u)^{1/2}$ by (T1), the above tail estimate and the Cauchy--Schwarz inequality.
Since this is the limit of a subsequence, we immediately get
\[\varliminf_k\Escr^{n_k}(\phi,u_k)\le\Escr(\phi,u).\]
This proves (F2'b).

To prove (F1b), let $u_n\in L^2(E;\mu_n)$, $u\in L^2(E;\mu)$, $u_n\to u$ $\Hscr$-weakly with $\varliminf_n(\tilde{\Escr}^n_1)^{1/2}(u_n)<\infty$.
By the assumed condition (M1) for the symmetric parts, $\tilde{\Escr}(u)\le\varliminf_n\tilde{\Escr}^n(u_n)<\infty$, hence $u\in\Dscr(\Escr)$, which is (F1b), cf. Remark~\ref{M1remark}.

The proof is complete.
\end{proof}

\section{Applications to stochastics}\label{sec:stoch}

\subsection{Convergence of laws}\label{ss:laws}

Let $E$ be either a separable real Banach space with dual space $E'$, which we call the infinite dimensional case, or $E=\R^d$, which we call the finite dimensional case. In the first case, define $\FCb$ as in the previous section and set $C:=\FCb$. In the second case, set $C:=C_0^\infty(\R^d)$. In both cases $E$ is a Polish space.
Denote by $\Bscr_b(E)$ the set
of all real-valued bounded Borel-measurable functions on $E$ and by $C_b(E)$ the set of all real-valued bounded and continuous
functions on $E$.

Let $\Omega:=C([0,\infty);E)$ be the space of continuous paths in $E$, equipped with the topology of uniform convergence on compact subsets of $[0,\infty)$. Then $\Omega$ is a Polish space. Let $X_t(\omega):=\omega(t)$, $t\ge 0$, $\omega\in\Omega$, be the coordinate process and let $\Mscr:=\sigma(X_t\;\vert\;t\ge 0)$, which coincides with the Borel $\sigma$-algebra of $\Omega$.
Let $\mu$ be a Borel probability measure on the Borel $\sigma$-algebra $\Bscr(E)$ of $E$. Let
\[\M=(\Omega,\Mscr,(X_t)_{t\ge 0},(\P_z)_{z\in E})\]
be a right process with lifetime $\zeta\equiv+\infty$, that is, the process is \textit{conservative}, and with continuous sample paths, realized canonically on $(\Omega,\Mscr)$.
Define a probability measure $\P$ on $(\Omega,\Mscr)$ by
\[\P(\Gamma):=\int_E \P_z(\Gamma)\mu(dz),\quad\Gamma\in\Mscr.\]
Let $J:=\{t_1,\ldots,t_k\}\subset [0,\infty)$, $0\le t_1<\cdots< t_k$, $k\in\N$.
The finite dimensional distribution $\P^J$ of $\P$ on $(E^J,\Bscr(E^J))$ is defined by
\[\P^J(A):=\P((X_{t_1},\ldots,X_{t_k})\in A),\quad A\in\Bscr(E^J).\]
Since $\Mscr$ is generated by the cylinder sets, which form a $\cap$-stable system, a probability measure on $(\Omega,\Mscr)$ is uniquely determined by its finite dimensional distributions.

Let $m_n$, $n\in\N$, and $m$ be fully supported Borel measures on $E$ such that $m_n\to m$ weakly, or vaguely if $E=\R^d$. In the infinite dimensional case, $m_n$ and $m$ are assumed to be finite.
Then
\[L^2(E;m_n)=:H_n\to H:=L^2(E;m)\]
in the sense of convergent Hilbert spaces, where $C\subset L^2(E;m)$ is dense, see e.g.~\cite[Section II.3.a)]{MR}.
As usual, set $\Hscr:=\disjoint_n H_n\dot\cup H$.

Assume that for every $n\in\N$ there exists a conservative right process
\[\M^{(n)}=(\Omega,\Mscr,(X_t^{(n)})_{t\ge 0},(\P^{(n)}_z)_{z\in E})\]
with continuous sample paths, realized canonically on the same space $(\Omega,\Mscr)$, i.e., $X_t^{(n)}=X_t$ is the coordinate process, and with common state space $E$. Define
$\P^{(n)}$, $\P^{(n),J}$, similarly w.r.t. some Borel probability measure $\mu_n$, $n\in\N$ on $\Bscr(E)$.

\begin{assumption}\label{stoch1}
\begin{enumerate}[(1)]
\item $\mu_n\to\mu$ weakly.
\item Each $\mu_n$, $n\in\N$, and $\mu$ is absolutely continuous w.r.t. $m_n$, $n\in\N$, and $m$, respectively, with
positive densities $h_n$, $n\in\N$, and $h$.
\item $h_n\in L^2(E;m_n)$, $n\in\N$, $h\in L^2(E;m)$, and $\sup_n\norm{h_n}_{L^2(E;m_n)}<\infty$.
\end{enumerate}
\end{assumption}
Note that, in this section, $\mu_n$ and $\mu$ denote the initial distributions and $h_n$, $h$ their densities, whereas the reference measures, which are denoted by $\mu_n$, $\mu$ in Sections~\ref{sec:mosco} and \ref{sec:genconv}, are denoted by $m_n$, $m$. When the results of Section~\ref{sec:genconv} are applied, $m_n$ and $m$ are the reference measures of Section~\ref{sec:genconv}.

By weak convergence, clearly
\[\lim_n\int_E \phi h_n\;dm_n=\lim_n\int_E\phi\;d\mu_n=\int_E\phi\;d\mu=\int_E\phi h\;dm\]
for any $\phi\in\Fscr C_b^\infty(E)\subset C_b(E)$. Hence we get by Lemma~\ref{weakmosco} that
$h_n\to h$ weakly in $\Hscr$. The argument is similar in the finite dimensional case with $C=C_0^\infty(\R^d)\subset C_b(\R^d)$. Note that only the weak convergence $\mu_n\to\mu$ is used here, whereas the vague convergence $m_n\to m$ is used for the convergence $H_n\to H$.

Here, a right process $\M$ with transition function $(p_t)_{t\ge0}$ is called \emph{associated} with a possibly generalized Dirichlet form $\Escr$ on $L^2(E;m)$ with $L^2$-semigroup $(T_t)_{t\ge0}$ if $p_tf$ is an $m$-version of $T_tf$ for every $f\in\Bscr_b(E)\cap L^2(E;m)$ and $t>0$, cf.~\cite{MR}, \cite{S1}. In particular, $p_t$ respects $m$-null sets. Moreover, Lemmas~\ref{phi-int} and \ref{phiuninfty} are used below in $L^2(E;m_n)$ for general fully supported measures $m_n$. Their proofs carry over verbatim, since they only use that $\phi\in C$ is bounded and that $\phi\cdot\psi\in C$ for $\phi,\psi\in C$.

\begin{assumption}\label{stoch2}
$\{\P^{(n)}\}_{n\in\N}$ is tight as a family of probability measures on the Polish space $\Omega$.
\end{assumption}

Assumption~\ref{stoch2} can be verified in many cases with the help of the Lyons--Zheng decomposition \cite{LyZh}, see also~\cite{FOT}, \cite{LZa}, \cite{LZb}, \cite{LZ}, \cite{RZ} and \cite{T}. For the non-sectorial case of generalized Dirichlet forms, see~\cite{Tr2}, and for the stochastic calculus of generalized Dirichlet forms, see~\cite{Tr}. Roughly speaking, the Lyons--Zheng decomposition represents $f(X_t)-f(X_0)$, for suitable functions $f$ in the domain of the form, for instance $f\in C$, as the sum of a forward martingale and a time-reversed backward martingale, plus, in the non-symmetric case, a term involving the antisymmetric part. The quadratic variations of the martingales are controlled by the energy of $f$. If the energies of a countable family of such functions, which separates the points of $E$, are bounded uniformly in $n$, together with a uniform control of the antisymmetric part such as a uniform sector condition, this yields uniform moment bounds for the increments of $f(X_t)$, hence tightness of the laws of $f(X)$ by Kolmogorov's criterion. Since these estimates are obtained for the processes started with the reference measures $m_n$, they have to be transferred to $\P^{(n)}$, e.g. by means of uniform bounds on the densities $h_n$. Tightness of $\{\P^{(n)}\}$ then follows if, in addition, a compact containment condition holds, cf.~\cite[Section 3.9]{EK}. See~\cite{LZb,LZ,RZ} for details and~\cite{Tr2} for the non-sectorial case.

\begin{thm}\label{lawconv}
Let $\M^{(n)}$, $n\in\N$, and $\M$ be right processes as above, that is, conservative, with continuous sample paths and realized canonically on $(\Omega,\Mscr)$, associated with
possibly generalized Dirichlet forms $\Escr^{(n)}$, $n\in\N$, and $\Escr$ on
$L^2(E;m_n)$, $n\in\N$, and $L^2(E;m)$, respectively. Denote by $(T_t^{(n)})_{t\ge 0}$, $n\in\N$, and $(T_t)_{t\ge 0}$
the associated $L^2$-semigroups. Let Assumptions~\ref{stoch1} and \ref{stoch2} hold.
If $T_t^{(n)}\to T_t$ strongly in $\Hscr$ for any $t\ge 0$, then
$\P^{(n)}\to\P$ weakly as probability measures on $\Omega$.
\end{thm}
\begin{proof}
Let $J=\{t_1,\ldots,t_k\}$, $0\le t_1<\cdots< t_k$, $k\in\N$ be a subset of $[0,\infty)$. Denote by $p_t(\cdot,dx)$, $t\ge 0$
the transition semigroup of $(X_t)_{t\ge 0}$. Let $g(x)=g_1(x_1)\cdots g_k(x_k)$, $x=(x_1,\ldots,x_k)\in E^J$, with $g_i\in C$, $1\le i\le k$. Then
\begin{align*}
\int g\;d\P^J=&\int g_1(X_{t_1})\cdots g_k(X_{t_k})d\P\\
=&\iint g_1(X_{t_1})\cdots g_k(X_{t_k})d\P_z\;\mu(dz)\\
=&\int\idotsint g_1(x_1)p_{t_1}(z,dx_1)g_2(x_2)p_{t_2-t_1}(x_1,dx_2)\cdots\\
&\cdots g_k(x_k)p_{t_k-t_{k-1}}(x_{k-1},dx_k)\mu(dz)\\
=&\int p_{t_1}(g_1 p_{t_2-t_1}(g_2\cdots p_{t_k-t_{k-1}} g_k\cdots))(z)\mu(dz).
\end{align*}
$p_t f$ is an $m$-version of $T_t f$ for every $f\in\Bscr_b(E)\cap L^2(E;m)$.\footnote{For a general right process, $p_tf$ is only universally measurable, which is harmless here, since all functions are integrated against $\mu\ll m$.} Then clearly
\[p_{t_{k-1}-t_{k-2}}(g_{k-1}p_{t_k-t_{k-1}}g_k)\]
is an $m$-version of 
\[T_{t_{k-1}-t_{k-2}}(g_{k-1}T_{t_k-t_{k-1}}g_k)\]
since $g_{k-1}p_{t_k-t_{k-1}}g_k$ is bounded and contained in the $L^2$-class of $g_{k-1}T_{t_k-t_{k-1}}g_k$.
By an induction argument
\[p_{t_1}(g_1 p_{t_2-t_1}(g_2\cdots p_{t_k-t_{k-1}} g_k\cdots))\]
is an $m$-version of
\[T_{t_1}(g_1 T_{t_2-t_1}(g_2\cdots T_{t_k-t_{k-1}} g_k\cdots)).\]
A similar statement holds for every $n\in\N$. Since $h_n\to h$ $\Hscr$-weakly as a consequence of Assumption~\ref{stoch1} and $T_t^{(n)}\to T_t$, $t\ge 0$
strongly in the sense of Definition~\ref{Operatorconvergence}, we get
\begin{align*}
&\lim_n\int g\;d\P^{(n),J}\\
=&\lim_n\int T_{t_1}^{(n)}(g_1 T_{t_2-t_1}^{(n)}(g_2\cdots T_{t_k-t_{k-1}}^{(n)} g_k\cdots))d\mu_n\\
=&\lim_n\int T_{t_1}^{(n)}(g_1 T_{t_2-t_1}^{(n)}(g_2\cdots T_{t_k-t_{k-1}}^{(n)} g_k\cdots)) h_n\;dm_n\\
=&\int T_{t_1}(g_1 T_{t_2-t_1}(g_2\cdots T_{t_k-t_{k-1}} g_k\cdots)) h\;dm\\
=&\int T_{t_1}(g_1 T_{t_2-t_1}(g_2\cdots T_{t_k-t_{k-1}} g_k\cdots))d\mu
\end{align*}
by an induction argument, since
\[T_{t_k-t_{k-1}}^{(n)}g_k\to T_{t_k-t_{k-1}} g_k\]
$\Hscr$-strongly by strong convergence of semigroups and
\[g_{k-1}T_{t_k-t_{k-1}}^{(n)}g_k\to g_{k-1}T_{t_k-t_{k-1}} g_k\]
$\Hscr$-strongly by Lemma~\ref{phiuninfty}, or by Lemma~\ref{phi-int} in the finite dimensional case.

Hence we obtain
\[\lim_n\int g\;d\P^{(n),J}=\int g\;d\P^J\]
for every product $g$ as above. If some of the factors $g_i$ are replaced by the constant function $1$, then $g$ is a product of the same type for the set $J'\subset J$ of the remaining times, and $\int g\,d\P^{(n),J}=\int g\,d\P^{(n),J'}$, $\int g\,d\P^{J}=\int g\,d\P^{J'}$. Since all measures involved are probability measures, we conclude that
\[\lim_n\int g\;d\P^{(n),J}=\int g\;d\P^J\quad\text{for every}\;g\in\Pscr_J,\]
where
\[\Pscr_J:=\operatorname{lin}\{g\;\vert\;g(x)=g_1(x_1)\cdots g_k(x_k),\;g_i\in C\cup\{1\},\;1\le i\le k\}.\]

Now assume that $\P^{(n)}\not\to\P$ weakly. Then there exist $f\in C_b(\Omega)$, $\eps>0$ and a subsequence $\{\P^{(n_k)}\}$ such that $\abs{\int f\,d\P^{(n_k)}-\int f\,d\P}\ge\eps$ for all $k$. By Assumption~\ref{stoch2} and Prokhorov's theorem, see~\cite[Theorem 5.1]{Bil}, there exists a further subsequence, again denoted by $\{\P^{(n_k)}\}$, which converges weakly to some probability measure $\P_1$ on $(\Omega,\Mscr)$, and necessarily $\P_1\neq\P$. Since the projection $\pi_J:\Omega\to E^J$, $\omega\mapsto(\omega(t_1),\ldots,\omega(t_k))$, is continuous on $\Omega=C([0,\infty);E)$, the continuous mapping theorem yields
\[\lim_k\int f\;d\P^{(n_k),J}=\int f\;d\P^{J}_1\]
for every $f\in C_b(E^J)$. As $\Pscr_J\subset C_b(E^J)$, we get
\[\int f\;d\P_1^J=\int f\;d\P^J\quad\text{for every}\;f\in\Pscr_J.\]
We show that $\sigma(\Pscr_J)=\Bscr(E^J)$. Clearly, $\Pscr_J$ is an algebra of bounded continuous functions which contains the constants. It separates the points of $E^J$, which we show by adapting the argument of~\cite[Section IV.4.b), (4.2)]{MR} for the points of $E$: Let $x\neq y$ in $E^J$, say $x_j\neq y_j$. In the infinite dimensional case, by the Hahn--Banach theorem, see~\cite[Theorem III.5]{ReSi1}, there exists $l\in E'$ with $l(x_j)\neq l(y_j)$, and $g_j:=\arctan\circ\, l\in\FCb$ satisfies $g_j(x_j)\neq g_j(y_j)$. In the finite dimensional case, choose $g_j\in C_0^\infty(\R^d)$ with $g_j(x_j)\neq g_j(y_j)$. Setting $g_i:\equiv 1$ for $i\neq j$, the function $g(z):=g_1(z_1)\cdots g_k(z_k)$ belongs to $\Pscr_J$ and separates $x$ and $y$. Hence
\[(E^J\times E^J)\setminus d=\bigcup_{g\in\Pscr_J}\{(x,y)\in E^J\times E^J\;\vert\;g(x)\neq g(y)\},\]
where $d$ denotes the diagonal in $E^J\times E^J$. Since every open subset of the separable metric space $E^J\times E^J$ is Lindel\"of, this open cover has a countable subcover, i.e., there exists a countable set $K\subset\Pscr_J$ which separates the points of $E^J$. Since $E^J$ is Polish, it follows from~\cite[Lemma 18, p. 108]{Schw} that $\sigma(K)=\Bscr(E^J)$, and hence $\sigma(\Pscr_J)=\Bscr(E^J)$, because $\sigma(K)\subset\sigma(\Pscr_J)\subset\Bscr(E^J)$. Now the monotone class theorem, see e.g.~\cite[Chapter I]{DM}, applied to the monotone vector space
\[H_0:=\left\{f\;\bigg\vert\;f\in\Bscr_b(E^J)\;\text{such that}\;\int f\;d\P_1^J=\int f\;d\P^J\right\},\]
which contains the algebra $\Pscr_J$, yields
\[\int f\;d\P_1^J=\int f\;d\P^J\]
for every $f\in\Bscr_b(E^J)$, $J$ as above.
Since the finite dimensional distributions determine a probability measure on $(\Omega,\Mscr)$ uniquely, we get $\P_1=\P$, which is a contradiction.
Thus $\P^{(n)}\to\P$ weakly.
\end{proof}

\begin{rem}\label{cadlag}
Theorem~\ref{lawconv} is formulated for processes with continuous sample paths, which is the relevant case for the diffusions associated with the forms of Sections~\ref{sec:mosco} and \ref{sec:genconv}, provided these processes exist, are conservative and have continuous sample paths. This has to be verified separately, for instance via quasi-regularity, see~\cite{MR}, \cite{S1}, and the same applies to Assumption~\ref{stoch2}. If $\M^{(n)}$, $\M$ are right processes with c\`adl\`ag sample paths, one can take $\Omega:=D([0,\infty);E)$ with the Skorokhod $J_1$-topology, which is Polish, and $\Mscr:=\sigma(X_t\;\vert\;t\ge 0)$, which again coincides with the Borel $\sigma$-algebra, see~\cite[Section 3.7]{EK}. In this case the evaluation map $\omega\mapsto\omega(t)$ is not continuous at paths having a jump at $t$, so that the continuous mapping theorem gives $\P^{(n_k),J}\to\P_1^J$ weakly only for $J\subset T(\P_1):=\{t\ge 0\;\vert\;\P_1(X_t\neq X_{t-})=0\}$, whose complement in $[0,\infty)$ is at most countable, see~\cite[Section 3.7]{EK}. The above proof then shows that $\P_1^J=\P^J$ for all finite $J\subset T(\P_1)$. Since $T(\P_1)$ is dense and the paths are right-continuous, the finite dimensional distributions along $T(\P_1)$ determine a probability measure on $(\Omega,\Mscr)$, hence $\P_1=\P$ and the conclusion of Theorem~\ref{lawconv} remains valid, cf.~\cite[Section 3.7]{EK}.
\end{rem}

\begin{rem}
If we had assumed instead that $h_n\to h$ $\Hscr$-strongly, which in our case would have followed from $\norm{h_n}_{L^2(E;m_n)}\to\norm{h}_{L^2(E;m)}$
by Lemma~\ref{normmosco}, weak convergence of semigroups $T_t^{(n)}\to T_t$, $t\ge 0$ would be sufficient for the same result, using Lemma~\ref{phiuninfty} or Lemma~\ref{phi-int}, respectively.
This complements the above result via Theorem~\ref{Hino}, which yields the weak convergence of the co-semigroups, and thus we obtain a sufficiently strong
type of convergence of
the dual processes $\hat{\M}^{(n)}$, $n\in\N$ to $\hat{\M}$, provided that dual right processes associated with the co-semigroups exist, have continuous sample paths and that their laws are tight.
\end{rem}

\section{Concluding remarks}\label{sec:concl}

The abstract theory of Section~\ref{sec:general} characterizes the convergence of non-symmetric generalized forms defined on varying Hilbert spaces by the strong convergence of the associated resolvents and semigroups, by the weak convergence of the co-resolvents and by the strong graph convergence of the generators. The applications in Sections~\ref{sec:genconv} and \ref{sec:stoch} show that it is intended for approximation problems in which both the reference or speed measures and the coefficients of non-symmetric elliptic operators change, in finite and infinite dimensions, and in which one is ultimately interested in the convergence of the laws of the associated Markov processes. Theorem~\ref{lawconv} yields this convergence, provided the associated processes exist, are conservative, have continuous sample paths and their laws are tight.

The theory has been developed further in the author's doctoral dissertation~\cite{Toe2}. There, the framework of Kuwae and Shioya is extended to varying Banach spaces and to Mosco convergence of convex functionals on varying Banach spaces. It is shown that strong convergence generates a sequential topology~\cite[Theorem 5.32]{Toe2}. Moreover, versions of Theorem~\ref{Kato}, including the equivalence with the strong graph convergence of the generators, and of Theorem~\ref{Hino} are obtained there by a different method, see~\cite[Proposition 7.2, Theorem 7.10 and Theorem 7.15]{Toe2}.

\appendix

\section{Remaining proofs from Section 2}

In this appendix we give the proofs of Lemmas~\ref{basicconv}, \ref{normmosco}, \ref{weakmosco}, \ref{dualconv}, \ref{opnormconv}, \ref{compactB} and of Theorem~\ref{KSConv}, which are taken from~\cite{KS} and have been completed or rewritten, and the proof of Lemma~\ref{metricdoubletrouble}, which is taken from~\cite{A}. The proofs of the results which were new in~\cite{Toe1} or have been substantially modified, in particular of Theorems~\ref{Kolmetric}, \ref{Kato} and \ref{Hino}, are given in Section~\ref{sec:general}.

Recall that we use the convention $H:=H_\infty$, $\Phi_n:=\Phi_{\infty,n}$ and
$C:=\Cscr_\infty$.

\begin{proof}[Proof of Lemma~\ref{basicconv}]
\begin{enumerate}[(1):]
\item For the ``if'' part, assume that $\norm{u_n}_{H_n}\to 0$. Set $\phi_m=0\in C$ for every $m\in\N$. Then
\[\lim_m\varlimsup_n\norm{u_n-\Phi_n(0)}_{H_n}=\lim_m\varlimsup_n\norm{u_n}_{H_n}=0,\]
hence clearly $u_n\to 0$ $\Hscr$-strongly.

For the ``only if'' part, assume that $u_n\to 0\in H$ strongly in $\Hscr$. There exists a sequence $\{\phi_m\}\subset C$
with $\phi_m\to 0$ in $H$ and 
\[\lim_m\varlimsup_n\norm{u_n-\Phi_n(\phi_m)}_{H_n}=0.\]
Clearly,
\[\varlimsup_n\norm{u_n}_{H_n}\le\varlimsup_n\norm{u_n-\Phi_n(\phi_m)}_{H_n}+\varlimsup_n\norm{\Phi_n(\phi_m)}_{H_n}\to 0\]
as $m\to\infty$.
\item Obvious.
\item Let $u_n\in H_n$, $n\in\N$, $u\in H$ such that $u_n\to u$ $\Hscr$-strongly. Choose $\{\phi_m\}\subset C$ with
$\lim_m\norm{\phi_m-u}_H=0$ and $\lim_m\varlimsup_n\norm{\Phi_n(\phi_m)-u_n}_{H_n}=0$. Evidently,
\[
\abs{\norm{u_n}_{H_n}-\norm{u}_H}\le\norm{u_n-\Phi_n(\phi_m)}_{H_n}+\abs{\norm{\Phi_n(\phi_m)}_{H_n}-\norm{u}_H}.
\]
The first term tends to zero as $m,n\to\infty$. Taking $n\to\infty$ in the second term, one obtains
\[\abs{\norm{\phi_m}_{H}-\norm{u}_H}\le\norm{\phi_m-u}_H,\]
which tends to $0$ as $m\to\infty$.
\item Let $\{\tilde{u}_m\}\subset C$ with $\tilde{u}_m\to u$ in $H$ and $\{\tilde{v}_m\}\subset C$ with $\tilde{v}_m\to v$ in $H$.
Then by the Cauchy--Schwarz inequality
\begin{multline}\label{uvuvuvuv}
\lrabs{(u_n,v_n)_{H_n}-(u,v)_H}\\
=\bigl|(u_n-\Phi_n(\tilde{u}_m),v_n)_{H_n}+(\Phi_n(\tilde{u}_m),v_n-\Phi_n(\tilde{v}_m))_{H_n}\\
+(\Phi_n(\tilde{u}_m),\Phi_n(\tilde{v}_m))_{H_n}-(u,v)_H\bigr|\\
\le\norm{v_n}_{H_n}\norm{u_n-\Phi_n(\tilde{u}_m)}_{H_n}\\+\norm{\Phi_n(\tilde{u}_m)}_{H_n}\norm{v_n-\Phi_n(\tilde{v}_m)}_{H_n}\\
+\lrabs{(\Phi_n(\tilde{u}_m),\Phi_n(\tilde{v}_m))_{H_n}-(u,v)_{H}}.
\end{multline}
By (3) $\{\norm{v_n}_{H_n}\}$ is uniformly bounded in $n$, and $\varlimsup_n\norm{\Phi_n(\tilde{u}_m)}_{H_n}=\norm{\tilde u_m}_H$ is bounded uniformly in $m$. So taking $\varlimsup_n$ and then the limit $m\to\infty$, the first two terms on the right-hand side of \eqref{uvuvuvuv} tend to $0$. But by the polarization identity and linearity of each $\Phi_n$
\begin{multline*}
(\Phi_n(\tilde{u}_m),\Phi_n(\tilde{v}_m))_{H_n}
=\frac{1}{4}\left[\norm{\Phi_n(\tilde{u}_m)+\Phi_n(\tilde{v}_m)}_{H_n}^2-\norm{\Phi_n(\tilde{u}_m)-\Phi_n(\tilde{v}_m)}_{H_n}^2\right]\\
\longrightarrow\frac{1}{4}\left[\norm{u+v}_{H}^2-\norm{u-v}_{H}^2\right]=(u,v)_H
\end{multline*}
as first $n\to\infty$ and then $m\to\infty$, thus the last term on the right-hand side of \eqref{uvuvuvuv} tends to $0$ as well.
\end{enumerate}
\end{proof}

\begin{proof}[Proof of Lemma~\ref{normmosco}]
The ``only if'' part follows from Lemma~\ref{basicconv} (3) and (4) combined with the fact
that $\Phi_n(\phi)\to\phi$ strongly for every $\phi\in C$.

To prove the ``if'' part, let $\phi_m\to u$ in $H$, $\{\phi_m\}\subset C$.
\begin{multline*}
\lim_m\varlimsup_n\norm{u_n-\Phi_n(\phi_m)}_{H_n}=\lim_m\varlimsup_n\lrabs{\norm{u_n}_{H_n}^2-2(u_n,\Phi_n(\phi_m))_{H_n}+\norm{\Phi_n(\phi_m)}_{H_n}^2}^{1/2}\\
=\lim_m\lrabs{\norm{u}_H^2-2(u,\phi_m)_H+\norm{\phi_m}_H^2}^{1/2}=\lim_m\norm{u-\phi_m}_H=0.
\end{multline*}
\end{proof}

\begin{proof}[Proof of Lemma~\ref{weakmosco}]
The ``only if'' part follows from the fact that $\Phi_n(\phi)\to\phi$ strongly, together with the bound $\sup_n\norm{u_n}_{H_n}<\infty$ from Lemma~\ref{dualconv} (2), whose proof only uses Lemma~\ref{basicconv} (1).

To prove the ``if'' part, let
$\{u_n\}$, $u$ as in the assertion and $(u_n,\Phi_n(\phi))_{H_n}\to(u,\phi)_H$ for all $\phi\in C$. Take $v_n\to v$ strongly in $\Hscr$.
By strong convergence there exists a sequence $\{\phi_m\}$,
$\phi_m\in C$ with $\norm{\phi_m-v}_H\to 0$. We have to prove $(u_n,v_n)_{H_n}\to (u,v)_H$. By the Cauchy--Schwarz inequality one obtains
\begin{align*}
&\lrabs{(u_n,v_n)_{H_n}-(u,v)_H}\\
\le&\lrabs{(u_n,v_n-\Phi_n(\phi_m))_{H_n}}+\lrabs{(u_n,\Phi_n(\phi_m))_{H_n}-(u,v)_H}\\
\le&\norm{u_n}_{H_n}\norm{v_n-\Phi_n(\phi_m)}_{H_n}+\lrabs{(u_n,\Phi_n(\phi_m))_{H_n}-(u,v)_H}.
\end{align*}
$\{\norm{u_n}_{H_n}\}$ is uniformly bounded in $n$ by assumption, and $\norm{v_n-\Phi_n(\phi_m)}_{H_n}$ tends to zero
as $m,n\to\infty$ by strong convergence. Furthermore, by assumption, $(u_n,\Phi_n(\phi_m))_{H_n}\to (u,\phi_m)_H$ as $n\to\infty$, so the
last expression yields by the Cauchy--Schwarz inequality
\[\lrabs{(u,\phi_m)_H-(u,v)_H}=\lrabs{(u,\phi_m-v)_H}\le\norm{u}_H\norm{v-\phi_m}_H\to 0\]
as $m\to\infty$.
\end{proof}

\begin{proof}[Proof of Lemma~\ref{dualconv}]
\begin{enumerate}[(1):]
\item Pick a complete orthonormal basis $\{e_k\}$ of $H$. By uniform boundedness
of $\{u_n\}$, passing to a subsequence, the limit
\[\lim_n(u_n,\Psi_n(e_1))_{H_n}=:a_1\in\R\]
exists. Similarly, for every $k$ there is a further subsequence along which
\[\lim_n(u_n,\Psi_n(e_k))_{H_n}=:a_k\in\R\]
exists. By a diagonal argument we can find a common subsequence $n_l\uparrow +\infty$ such
that
\[\lim_l (u_{n_l},\Psi_{n_l}(e_k))_{H_{n_l}}=a_k\]
for every $k\in\N$.

Fix $N\in\N$. Let $\Lscr_n^N:=\operatorname{lin}\{\Psi_n(e_k)\;\vert\;k=1,\ldots,N\}\subset H_n$ and
$P_{\Lscr^N_n}:H_n\to\Lscr^N_n$ be the orthogonal projection on the
finite dimensional linear subspace $\Lscr^N_n$ of $H_n$.
$\{\Psi_n(e_k)\;\vert\;k=1,\ldots,N\}$ clearly is an orthonormal basis of
$\Lscr_n^N$, since the $\Psi_n$'s are unitary operators.
Clearly for every $n\in\N$,
\[P_{\Lscr_n^N}(u_n)=\sum_{k=1}^N(u_n,\Psi_n(e_k))_{H_n}\Psi_n(e_k),\]
and therefore, by the orthonormality of the $\Psi_n(e_k)$'s,
\[\sum_{k=1}^N\abs{a_k}^2=\lim_l\sum_{k=1}^N\lrabs{(u_{n_l},\Psi_{n_l}(e_k))_{H_{n_l}}}^2\le\varlimsup_n\norm{P_{\Lscr_n^N}(u_n)}_{H_n}^2\le\varlimsup_n\norm{u_n}_{H_n}^2<\infty\]
for any $N\in\N$. This gives us the existence of
\[u:=\sum_{k=1}^\infty a_k e_k\in H.\]
Let $\phi\in\bigcup_{N\ge 1}\operatorname{lin}\{e_k\;\vert\;k=1,\ldots,N\}=:C_0\subset H$.
It suffices to prove that
\[\lim_l(u_{n_l},\Psi_{n_l}(\phi))_{H_{n_l}}=(u,\phi)_H\]
for every such $\phi$, noting that $C_0\subset H$ is dense. Indeed, if $v_l\to v$ strongly along $\{n_l\}$, then $\norm{v_l-\Psi_{n_l}v}_{H_{n_l}}\to0$ by Theorem~\ref{Kolmetric}, and for $\phi\in C_0$
\begin{multline*}\abs{(u_{n_l},v_l)_{H_{n_l}}-(u,v)_H}\le\sup_n\norm{u_n}_{H_n}\left(\norm{v_l-\Psi_{n_l}v}_{H_{n_l}}+\norm{v-\phi}_H\right)\\+\abs{(u_{n_l},\Psi_{n_l}\phi)_{H_{n_l}}-(u,\phi)_H}+\norm{u}_H\norm{\phi-v}_H,\end{multline*}
which becomes arbitrarily small for large $l$ and $\phi$ close to $v$.

So, let $\phi=\sum_{k=1}^N (\phi,e_k)_H e_k$, i.e., $\phi$ depends only on the
first $N$ coordinates. Then,
\begin{align*}
(u,\phi)_H&=\sum_{k=1}^N a_k(\phi,e_k)_H\\
&=\lim_l\sum_{k=1}^N (u_{n_l},\Psi_{n_l}(e_k))_{H_{n_l}}(\Psi_{n_l}(e_k),\Psi_{n_l}(\phi))_{H_{n_l}}\\
&=\lim_l(P_{\Lscr_{n_l}^N}(u_{n_l}),\Psi_{n_l}(\phi))_{H_{n_l}}\\
&=\lim_l (u_{n_l},\Psi_{n_l}(\phi))_{H_{n_l}}.
\end{align*}
\item Suppose $\sup_n\norm{u_n}_{H_n}=\infty$. Then there exists a subsequence $\{u_{n_k}\}$ of $\{u_n\}$ such that $\norm{u_{n_k}}_{H_{n_k}}\ge k$.
Setting
\[v_k:=\frac{1}{k}\cdot\frac{u_{n_k}}{\norm{u_{n_k}}_{H_{n_k}}},\]
one has $\norm{v_k}_{H_{n_k}}={1/k}\to 0$. Extending $\{v_k\}$ by $v_n:=0$ for $n\notin\{n_k\;\vert\;k\in\N\}$, we obtain a sequence with $\norm{v_n}_{H_n}\to0$, hence by Lemma~\ref{basicconv} (1) $v_n\to 0$ in $\Hscr$, which implies
$(u_{n_k},v_k)_{H_{n_k}}\to (u,0)_H=0$ by weak convergence of $\{u_n\}$. On the other hand,
\[(u_{n_k},v_k)_{H_{n_k}}=\frac{1}{k}\norm{u_{n_k}}_{H_{n_k}}\ge 1.\]
This is a contradiction and thus we obtain $\sup_n\norm{u_n}_{H_n}<\infty$.

Let $\{v_n\}$ be a sequence with $v_n\in H_n$ which strongly converges to $u$. Such a sequence exists by
Corollary~\ref{exstrongconv}. Then, $(u_n,v_n)_{H_n}\to (u,u)_H$. Hence,
\begin{align*}
0&\le\varliminf_n\norm{u_n-v_n}_{H_n}^2\\
&=\varliminf_n(\norm{u_n}_{H_n}^2-2(u_n,v_n)_{H_n}+\norm{v_n}_{H_n}^2)\\
&=\varliminf_n\norm{u_n}_{H_n}^2-\norm{u}_H^2.
\end{align*}
This completes the proof of the first assertion. The second follows from Lemma~\ref{basicconv} (3) and Lemma~\ref{normmosco}.
\item The ``only if'' part is trivial. We prove the ``if'' part. The assumption implies that $u_n\to u$ weakly. Setting $v_n:=u_n$ and $v:=u$
in the assumption, we have $\norm{u_n}_{H_n}\to\norm{u}_H$, which proves the assertion by statement (2) of this lemma.
\end{enumerate}
\end{proof}

\begin{proof}[Proof of Lemma~\ref{opnormconv}]
\begin{enumerate}[(1):]
\item For any $\eps>0$ there is a unit vector $u\in H$ such that $\norm{Bu}_H>\norm{B}_{\Lscr(H)}-\eps$. Let
$u_n\in H_n$ be strongly converging to $u$. Note that $\norm{u_n}_{H_n}\to 1$. Since $B_n\to B$ strongly, we have
$\norm{B_n u_n}_{H_n}\to\norm{Bu}_H$ and therefore,
\[\varliminf_n\norm{B_n}_{\Lscr(H_n)}\ge\varliminf_n\frac{\norm{B_n u_n}_{H_n}}{\norm{u_n}_{H_n}}=\norm{Bu}_H>\norm{B}_{\Lscr(H)}-\eps,\]
which gives the desired statement.
\item Choose a subsequence $\{n_k\}$ with $\norm{B_{n_k}}_{\Lscr(H_{n_k})}\to\varlimsup_n\norm{B_n}_{\Lscr(H_n)}$ and unit vectors $u_k\in H_{n_k}$ such that $\abs{\norm{B_{n_k}}_{\Lscr(H_{n_k})}-\norm{B_{n_k} u_k}_{H_{n_k}}}\to 0$. Passing to a further subsequence, we may assume by Lemma~\ref{dualconv} (1) that $u_k$ converges weakly along $\{n_k\}$ to a vector $u\in H$ with $\norm{u}_H\le 1$. Extending $\{u_k\}$ to a weakly convergent full sequence by Corollary~\ref{exstrongconv} (iii) and using the compact convergence $B_n\to B$ and Corollary~\ref{exstrongconv} (iv), we obtain $B_{n_k}u_k\to Bu$ strongly along $\{n_k\}$. Hence
\[\norm{B}_{\Lscr(H)}
\ge\norm{Bu}_H=\lim_{k}\norm{B_{n_k} u_{k}}_{H_{n_k}}=\varlimsup_n\norm{B_n}_{\Lscr(H_n)},\]
which together with (1) completes the proof.
\end{enumerate}
\end{proof}

\begin{proof}[Proof of Lemma~\ref{compactB}]
Note that the compactness of $B$ is equivalent to that of $\hat{B}$. Let $\{u_m\}\subset H$ be a sequence of vectors in $H$ weakly converging to a vector
$u\in H$. It suffices to prove that $\hat{B}u_m\to\hat{B}u$ $H$-strongly. We easily see that $\hat{B}u_m\to\hat{B}u$ $H$-weakly. For each
$m\in\N$, there is a sequence $u_{m,n}$ with $u_{m,n}\in H_n$ such that $\lim_n u_{m,n}=u_m$ $\Hscr$-strongly. Since $\hat{B}_n\to\hat{B}$ strongly,
we have $\lim_n\hat{B}_n u_{m,n}=\hat{B}u_m$ $\Hscr$-strongly. Since $\norm{u_{m,n}-\Psi_n u_m}_{H_n}\to0$ and $\norm{\hat{B}_n u_{m,n}}_{H_n}\to\norm{\hat{B}u_m}_H$ as $n\to\infty$ for every $m$, we can choose a strictly increasing sequence of natural numbers $n_m\uparrow\infty$ such that $\norm{u_{m,n_m}-\Psi_{n_m}u_m}_{H_{n_m}}<1/m$ and
\begin{gather}\label{Bumnm}
\abs{\norm{\hat{B}_{n_m}u_{m,n_m}}_{H_{n_m}}-\norm{\hat{B}u_m}_H}<\frac1m.
\end{gather}
Then $u_{m,n_m}\to u$ weakly along $\{n_m\}$. Indeed, let $v_m\to v$ strongly along $\{n_m\}$. Then $\norm{\Psi_{n_m}^{-1}v_m-v}_H=\norm{v_m-\Psi_{n_m}v}_{H_{n_m}}\to0$ by Theorem~\ref{Kolmetric}. The number $M:=\sup_m\norm{u_m}_H+\sup_m\norm{v_m}_{H_{n_m}}$ is finite. Let $\eps>0$, and choose $m_0\in\N$ such that $1/m_0\le\eps$, $\norm{\Psi_{n_m}^{-1}v_m-v}_H\le\eps$ and $\abs{(u_m,v)_H-(u,v)_H}\le\eps$ for all $m\ge m_0$. Since $\Psi_{n_m}$ is unitary, we obtain for all $m\ge m_0$
\begin{multline*}
\lrabs{(u_{m,n_m},v_m)_{H_{n_m}}-(u,v)_H}\le\norm{u_{m,n_m}-\Psi_{n_m}u_m}_{H_{n_m}}\norm{v_m}_{H_{n_m}}\\
+\lrabs{(u_m,\Psi_{n_m}^{-1}v_m-v)_H}+\lrabs{(u_m,v)_H-(u,v)_H}\le(2M+1)\eps.
\end{multline*}
Hence $(u_{m,n_m},v_m)_{H_{n_m}}\to(u,v)_H$.
The compact convergence $\hat{B}_n\to\hat{B}$ from Corollary~\ref{Badjoint} (2), applied along $\{n_m\}$ by means of Corollary~\ref{exstrongconv} (iii) and (iv), shows that $\lim_m\hat{B}_{n_m}u_{m,n_m}=\hat{B}u$ $\Hscr$-strongly.
Hence, by \eqref{Bumnm}, $\norm{\hat{B}u_m}_H\to\norm{\hat{B}u}_H$ and so $\hat{B}u_m\to\hat{B}u$ $H$-strongly. This completes the proof.
\end{proof}

\begin{proof}[Proof of Theorem~\ref{KSConv}]
The proof uses only our notions of convergence and basic facts about symmetric closed forms, which can be found e.g. in~\cite{MR} or \cite{FOT}. In particular, the proof of Proposition~\ref{BetaMR}, which is stated after this theorem, does not
depend on the results of this paper and can be found in~\cite[Chapter I.2]{MR}. We also use Lemma~\ref{metricdoubletrouble}, whose proof
is self-contained and can be found later in this Appendix.

(2) $\Leftrightarrow$ (3) is a special case of Theorem~\ref{Kato}.

Let us prove (1) $\Rightarrow$ (2). Let $\a>0$, $\{u_n\}$, $u_n\in H_n$, $u_n\to u\in H$ strongly. Define $z_n:=G_\a^n u_n$, $z:=G_\a u$.
The vector $z_n$ is characterized as the unique minimizer of $v\mapsto\Escr^n(v,v)+\a(v,v)_{H_n}-2(u_n,v)_{H_n}$ over $H_n$, see~\cite[Proof of Theorem I.2.6]{MR} and also~\cite[Theorem I.2.8]{MR}.
Since $\sup_n\norm{G_\a^n}_{\Lscr(H_n)}\le\a^{-1}$, we can extract a subsequence of $\{z_n\}$ by Lemma~\ref{dualconv} (1),
still denoted by $\{z_n\}$, which converges weakly to some $\tilde{z}\in H$. For an arbitrary $v\in H$, by condition (M2), we
can find a sequence $\{v_n\}$, $v_n\in H_n$, $v_n\to v$ strongly such that $\lim_n\Escr^n(v_n,v_n)=\Escr(v,v)$. Since for every $n$,
\begin{equation}\label{minimizer}
\Escr^n(z_n,z_n)+\a(z_n,z_n)_{H_n}-2(u_n,z_n)_{H_n}\le\Escr^n(v_n,v_n)+\a(v_n,v_n)_{H_n}-2(u_n,v_n)_{H_n},
\end{equation}
by taking condition (M1), Lemma~\ref{basicconv} (3) and Lemma~\ref{dualconv} (2) into account, we find in the limit
\[\Escr(\tilde{z},\tilde{z})+\a(\tilde{z},\tilde{z})_{H}-2(u,\tilde{z})_{H}\le\Escr(v,v)+\a(v,v)_H-2(u,v)_H.\]
Therefore $\tilde{z}=G_\a u$. Here (M1) is applied along the subsequence, which is justified as in the proof of Lemma~\ref{gammamosco}, and (M2) along a subsequence is trivial, since one can restrict the recovery sequence to the subsequence. The same extraction and identification argument applies to every subsequence of $\{z_n\}$. By the uniqueness of such $\tilde{z}$ and Corollary~\ref{exstrongconv} (iv), it follows that $z_n\to z$ weakly. We now prove that
$(z_n,z_n)_{H_n}\to(z,z)_H$. By (M2) we choose $\{v_n\}$, $v_n\in H_n$, $v_n\to z$ strongly such that $\lim_n\Escr^n(v_n,v_n)=\Escr(z,z)$,
therefore, by rewriting \eqref{minimizer} as
\[\Escr^n(z_n,z_n)+\a\lrnorm{z_n-\frac{u_n}{\a}}_{H_n}^2\le\Escr^n(v_n,v_n)+\a\lrnorm{v_n-\frac{u_n}{\a}}_{H_n}^2\]
we get in the limit again by condition (M1) and Lemma~\ref{basicconv} (3)
\[\varlimsup_n\lrnorm{z_n-\frac{u_n}{\a}}_{H_n}^2\le\lrnorm{z-\frac{u}{\a}}_H^2,\]
hence $\lrnorm{z_n-\frac{u_n}{\a}}_{H_n}^2\to\lrnorm{z-\frac{u}{\a}}_H^2$, and this together with Lemma~\ref{basicconv} (3) and Lemma~\ref{dualconv} (2)
concludes the proof.

Let us now prove (2) $\Rightarrow$ (1). Suppose $G_\a^n$ converges strongly to $G_\a$ for every $\a>0$. We first prove (M1).
For $\a>0$, $n\in\N$ define approximate forms
\[\Escr^{n,(\a)}(u,v):=\a(u-\a G_\a^n u,v)_{H_n},\quad u,v\in H_n\]
see Proposition~\ref{BetaMR} and~\cite[Chapter I.2]{MR}.
Now let $\{v_n\}$, $v_n\in H_n$, $v_n\to u\in H$ weakly, and let $\{u_n\}$, $u_n\in H_n$, be a sequence converging strongly to $u$, which exists by Corollary~\ref{exstrongconv}. By the strong convergence of $G_\a^n$ we have
\[\lim_n\Escr^{n,(\a)}(u_n,u_n)=\Escr^{(\a)}(u,u).\] First observe that for every $n$, $\a>0$,
\begin{multline*}
\Escr^{n,(\a)}(v_n,v_n)-\Escr^{n,(\a)}(u_n,u_n)-2\Escr^{n,(\a)}(u_n,v_n-u_n)\\
=\Escr^{n,(\a)}(v_n,v_n)+\Escr^{n,(\a)}(u_n,u_n)-2\Escr^{n,(\a)}(u_n,v_n)\\
=\Escr^{n,(\a)}(u_n-v_n,u_n-v_n)\ge 0.
\end{multline*}
So for every $n$, $\a>0$,
\[
\Escr^n(v_n,v_n)\ge\Escr^{n,(\a)}(v_n,v_n)\ge\Escr^{n,(\a)}(u_n,u_n)+2\Escr^{n,(\a)}(u_n,v_n-u_n).
\]
Clearly, $v_n-u_n\to 0$ weakly. Thus, by the strong convergence of $\{G_\a^n\}$, we get
$2\Escr^{n,(\a)}(u_n,v_n-u_n)=2\a(u_n-\a G_\a^n u_n,v_n-u_n)_{H_n}\to 0$. Taking $\varliminf_n$ we get for every $\a>0$
\[\varliminf_n\Escr^{n}(v_n,v_n)\ge\Escr^{(\a)}(u,u).\]
Since this holds for any $\a>0$, we conclude
\[\varliminf_n\Escr^{n}(v_n,v_n)\ge\Escr(u,u).\]
This proves (M1).

To prove (M2), let first $u\in H\setminus\Dscr(\Escr)$. Then (M2) holds for every strongly convergent sequence $u_n\to u$ by (M1). Now let $u\in\Dscr(\Escr)$ and let $u_n\to u$ strongly. By strong convergence of $\{G_\a^n\}$ we have
\[\Escr(u,u)=\lim_{\a\to\infty}\lim_n\Escr^{n,(\a)}(u_n,u_n)\]
and
\[\lim_{\a\to\infty}\lim_n \a G_\a^n u_n=u.\]
Since $\Hscr$ is metrizable and $\Escr(u,u)<\infty$, we can apply Corollary~\ref{doubletrouble} with $\a=m\in\N$ and pick a sequence of natural numbers $\a_n\uparrow\infty$ with
\[\Escr(u,u)=\lim_n\Escr^{n,(\a_n)}(u_n,u_n)\quad\text{and}\quad\lim_n \a_n G_{\a_n}^n u_n=u.\]
Set $w_n:=\a_n G_{\a_n}^n u_n$. Observe that $w_n\to u$ strongly. Using Proposition~\ref{BetaMR} (ii) we get
\[\Escr^{n,(\a_n)}(u_n,u_n)=\Escr^n(w_n,w_n)+\a_n \norm{u_n-w_n}_{H_n}^2\ge\Escr^n(w_n,w_n)\]
for every $n$. Taking $\varlimsup_n$ and using (M1) proves the assertion.
\end{proof}

\begin{proof}[Proof of Lemma~\ref{metricdoubletrouble}]
The proof is taken from~\cite[Lemma 1.15 and Corollary 1.18]{A}. First, let $\{a_{n,m}\}\subset\ol{\R}$, $n\in\N$, $m\in\N$ be a double indexed sequence.
We will prove that there exists a nondecreasing mapping $n\mapsto m(n)$ tending to $+\infty$, such that
\begin{equation}\label{infinfinf}
\varliminf_{n\to+\infty}a_{n,m(n)}\ge\varliminf_{m\to+\infty}\varliminf_{n\to+\infty}a_{n,m}.
\end{equation}
Let $a_{m}:=\varliminf_{n\to+\infty}a_{n,m}$ and $a:=\varliminf_{m\to+\infty}a_m$. If $a=-\infty$, there is nothing to prove. Assume therefore
that $a>-\infty$, and define a sequence $(\tilde{a}_p)_{p\in\N}$ of real numbers by $\tilde{a}_p:=a-2^{-p}$ if $a<+\infty$, and by $\tilde{a}_p:=p$ if $a=+\infty$.
By definition of $a$, there exists a strictly increasing sequence $(m_p)_{p\in\N}$, $m_p\uparrow +\infty$ such that
\[a_m>\tilde{a}_p\quad\text{for all}\quad m\ge m_p.\]
In other words,
\[a_m>(a-2^{-p})\wedge p\quad\text{for all}\quad m\ge m_p.\]
In the same way, there exists a strictly increasing sequence $(n_p)_{p\in\N}$, $n_p\uparrow +\infty$ such that
\[a_{n,m_p}>(a_{m_p}-2^{-p})\wedge p\quad\text{for all}\quad n\ge n_p.\]
Set $m(n):=m_p$ if $n_p\le n<n_{p+1}$, and $m(n):=m_1$ for $n<n_1$, and verify that \eqref{infinfinf} is satisfied: when $n_p\le n<n_{p+1}$, we get by the above
\[a_{n,m(n)}>(a_{m_p}-2^{-p})\wedge p\ge(((a-2^{-p})\wedge p)-2^{-p})\wedge p.\]
Thus, for all $n\ge n_p$
\[a_{n,m(n)}>(((a-2^{-p})\wedge p)-2^{-p})\wedge p.\]
It follows that
\[\varliminf_{n\to+\infty}a_{n,m(n)}\ge(((a-2^{-p})\wedge p)-2^{-p})\wedge p.\]
This being true for any $p\in\N$, using the fact that for any $a\in\ol{\R}$,
\[(((a-2^{-p})\wedge p)-2^{-p})\wedge p\]
increases to $a$ as $p$ goes to $+\infty$, we get:
\[\varliminf_{n\to+\infty}a_{n,m(n)}\ge a=\varliminf_{m\to+\infty}\varliminf_{n\to+\infty}a_{n,m}.\]
This proves our first assertion. Clearly,
replacing $a_{n,m}$ by $-a_{n,m}$ we obtain
\begin{equation}\label{supsupsup}
\varlimsup_{n\to+\infty}a_{n,m(n)}\le\varlimsup_{m\to+\infty}\varlimsup_{n\to+\infty}a_{n,m}.
\end{equation}

Now we prove the lemma. Let $\{x_{n,m}\}$, $\{x_m\}$, $x$ be as above. Set $a_{n,m}:=d(x_{n,m},x)\in\R$. By the above
there exists a nondecreasing map $n\mapsto m(n)$ such that \eqref{supsupsup} holds. By definition of $a_{n,m}$,
\[\varlimsup_{n\to+\infty}a_{n,m}=\lim_{n\to+\infty} d(x_{n,m},x)=d(x_m,x)\]
and
\[\varlimsup_{m\to+\infty}\varlimsup_{n\to+\infty}a_{n,m}=\lim_{m\to+\infty} d(x_m,x)=0.\]
So,
\[\varlimsup_{n\to+\infty} d(x_{n,m(n)},x)=0,\]
which proves the assertion.
\end{proof}


\begin{thebibliography}{10}
\providecommand{\url}[1]{{#1}}
\providecommand{\urlprefix}{URL }
\expandafter\ifx\csname urlstyle\endcsname\relax
  \providecommand{\doi}[1]{DOI~\discretionary{}{}{}#1}\else
  \providecommand{\doi}{DOI~\discretionary{}{}{}\begingroup
  \urlstyle{rm}\Url}\fi

\bibitem{AHKS}
Albeverio, S., H{\o}egh-Krohn, R., Streit, L.: Regularization of {H}amiltonians
  and processes.
\newblock J. Math. Phys. \textbf{21}(7), 1636--1642 (1980)

\bibitem{AKS}
Albeverio, S., Kusuoka, S., Streit, L.: Convergence of {D}irichlet forms and
  associated {S}chr{\"o}dinger operators.
\newblock J. Funct. Anal. \textbf{68}(2), 130--148 (1986)

\bibitem{AR}
Albeverio, S., R{\"o}ckner, M.: Classical {D}irichlet forms on topological
  vector spaces -- closability and a {C}ameron--{M}artin formula.
\newblock J. Funct. Anal. \textbf{88}(2), 395--436 (1990)

\bibitem{A}
Attouch, H.: Variational Convergence for Functions and Operators.
\newblock Applicable Mathematics Series. Pitman, Boston (1984)

\bibitem{BvR14}
Barret, F., von Renesse, M.: Averaging principle for diffusion processes via
  {D}irichlet forms.
\newblock Potential Anal. \textbf{41}(4), 1033--1063 (2014)

\bibitem{Bil}
Billingsley, P.: Convergence of Probability Measures, 2nd edn.
\newblock Wiley Series in Probability and Statistics. John Wiley \& Sons, New
  York (1999)

\bibitem{CS}
Camar-Eddine, M., Seppecher, P.: Closure of the set of diffusion functionals
  with respect to the {M}osco-convergence.
\newblock Math. Models Methods Appl. Sci. \textbf{12}(8), 1153--1176 (2002)

\bibitem{Can}
Cannon, J.T.: Convergence criteria for a sequence of semi-groups.
\newblock Applicable Anal. \textbf{5}(1), 23--31 (1975)

\bibitem{CG10}
Conrad, F., Grothaus, M.: Construction, ergodicity and rate of convergence of
  {$N$}-particle {L}angevin dynamics with singular potentials.
\newblock J. Evol. Equ. \textbf{10}(3), 623--662 (2010)

\bibitem{CG11}
Conrad, F., Grothaus, M.: {$N/V$}-limit for {L}angevin dynamics in continuum.
\newblock Rev. Math. Phys. \textbf{23}(1), 1--51 (2011)

\bibitem{DaMa}
Dal~Maso, G.: An Introduction to {$\Gamma$}-Convergence, \emph{Progress in
  Nonlinear Differential Equations and their Applications}, vol.~8.
\newblock Birkh{\"a}user, Boston, MA (1993)

\bibitem{DGF}
De~Giorgi, E., Franzoni, T.: Su un tipo di convergenza variazionale.
\newblock Atti Accad. Naz. Lincei Rend. Cl. Sci. Fis. Mat. Natur. (8)
  \textbf{58}(6), 842--850 (1975)

\bibitem{DGS}
De~Giorgi, E., Spagnolo, S.: Sulla convergenza degli integrali dell'energia per
  operatori ellittici del secondo ordine.
\newblock Boll. Un. Mat. Ital. (4) \textbf{8}, 391--411 (1973)

\bibitem{DM}
Dellacherie, C., Meyer, P.A.: Probabilities and Potential, \emph{North-Holland
  Mathematics Studies}, vol.~29.
\newblock North-Holland, Amsterdam (1978)

\bibitem{Eb}
Eberle, A.: Uniqueness and Non-Uniqueness of Semigroups Generated by Singular
  Diffusion Operators, \emph{Lecture Notes in Mathematics}, vol. 1718.
\newblock Springer, Berlin (1999)

\bibitem{EK}
Ethier, S.N., Kurtz, T.G.: Markov Processes: Characterization and Convergence.
\newblock Wiley Series in Probability and Mathematical Statistics. John Wiley
  \& Sons, New York (1986)

\bibitem{FOT}
Fukushima, M., Oshima, Y., Takeda, M.: Dirichlet Forms and Symmetric {M}arkov
  Processes, \emph{de Gruyter Studies in Mathematics}, vol.~19, 2nd edn.
\newblock Walter de Gruyter, Berlin (2011)

\bibitem{GN20}
Grothaus, M., Nonnenmacher, A.: Overdamped limit of generalized stochastic
  {H}amiltonian systems for singular interaction potentials.
\newblock J. Evol. Equ. \textbf{20}(2), 577--605 (2020)

\bibitem{GW24}
Grothaus, M., Wittmann, S.: Mosco convergence of gradient forms with non-convex
  interaction potential.
\newblock Integral Equations Operator Theory \textbf{96}(4), Paper No. 29
  (2024)

\bibitem{GW25}
Grothaus, M., Wittmann, S.: Mosco convergence of gradient forms with non-convex
  potentials {II}.
\newblock Potential Anal. \textbf{63}(3), 1593--1643 (2025)

\bibitem{H}
Hino, M.: Convergence of non-symmetric forms.
\newblock J. Math. Kyoto Univ. \textbf{38}(2), 329--341 (1998)

\bibitem{HM22}
Hinz, M., Meinert, M.: Approximation of partial differential equations on
  compact resistance spaces.
\newblock Calc. Var. Partial Differential Equations \textbf{61}(1), Paper No.
  19 (2022)

\bibitem{Ito}
Ito, H.: Discrete approximation of semi-{D}irichlet forms with non-symmetric
  perturbations of resistance forms.
\newblock Preprint, arXiv:2605.27974 (2026).
\newblock \urlprefix\url{https://arxiv.org/abs/2605.27974}

\bibitem{Kas}
Kasue, A.: Convergence of metric measure spaces and energy forms.
\newblock In: Selected Papers on Differential Equations and Analysis,
  \emph{Amer. Math. Soc. Transl. Ser. 2}, vol. 215, pp. 79--96. American
  Mathematical Society, Providence, RI (2005)

\bibitem{Ka}
Kato, T.: Perturbation Theory for Linear Operators, \emph{Grundlehren der
  Mathematischen Wissenschaften}, vol. 132.
\newblock Springer, New York (1966)

\bibitem{K1}
Kolesnikov, A.V.: Convergence of {D}irichlet forms with changing speed measures
  on {$\mathbb{R}^d$}.
\newblock Forum Math. \textbf{17}(2), 225--259 (2005)

\bibitem{K2}
Kolesnikov, A.V.: Mosco convergence of {D}irichlet forms in infinite dimensions
  with changing reference measures.
\newblock J. Funct. Anal. \textbf{230}(2), 382--418 (2006)

\bibitem{K3}
Kolesnikov, A.V.: Weak convergence of diffusion processes on {W}iener space.
\newblock Probab. Theory Related Fields \textbf{140}(1-2), 1--17 (2008)

\bibitem{KS}
Kuwae, K., Shioya, T.: Convergence of spectral structures: a functional
  analytic theory and its applications to spectral geometry.
\newblock Comm. Anal. Geom. \textbf{11}(4), 599--673 (2003)

\bibitem{KS2}
Kuwae, K., Shioya, T.: Variational convergence over metric spaces.
\newblock Trans. Amer. Math. Soc. \textbf{360}(1), 35--75 (2008)

\bibitem{KU2}
Kuwae, K., Uemura, T.: Weak convergence of symmetric diffusion processes {II}.
\newblock In: S.~Watanabe, M.~Fukushima, Y.V. Prohorov, A.N. Shiryaev (eds.)
  Probability Theory and Mathematical Statistics ({T}okyo, 1995), pp. 256--265.
  World Scientific, River Edge, NJ (1996)

\bibitem{KU1}
Kuwae, K., Uemura, T.: Weak convergence of symmetric diffusion processes.
\newblock Probab. Theory Related Fields \textbf{109}(2), 159--182 (1997)

\bibitem{LM}
Lancaster, P., Tismenetsky, M.: The Theory of Matrices, 2nd edn.
\newblock Academic Press, Orlando, FL (1985)

\bibitem{LRV16}
Lancia, M.R., Regis~Durante, V., Vernole, P.: Asymptotics for {V}enttsel'
  problems for operators in non divergence form in irregular domains.
\newblock Discrete Contin. Dyn. Syst. Ser. S \textbf{9}(5), 1493--1520 (2016)

\bibitem{Lob}
L{\"o}bus, J.U.: Mosco type convergence of bilinear forms and weak convergence
  of {$n$}-particle systems.
\newblock Potential Anal. \textbf{43}(2), 241--267 (2015)

\bibitem{LZa}
Lyons, T.J., Zhang, T.S.: Note on convergence of {D}irichlet processes.
\newblock Bull. London Math. Soc. \textbf{25}(4), 353--356 (1993)

\bibitem{LZb}
Lyons, T.J., Zhang, T.S.: Decomposition of {D}irichlet processes and its
  application.
\newblock Ann. Probab. \textbf{22}(1), 494--524 (1994)

\bibitem{LZ}
Lyons, T.J., Zhang, T.S.: Convergence of non-symmetric {D}irichlet processes.
\newblock Stochastics Stochastics Rep. \textbf{57}(3-4), 159--167 (1996)

\bibitem{LyZh}
Lyons, T.J., Zheng, W.: A crossing estimate for the canonical process on a
  {D}irichlet space and a tightness result.
\newblock Ast\'erisque \textbf{157--158}, 249--271 (1988)

\bibitem{MR}
Ma, Z.M., R{\"o}ckner, M.: Introduction to the Theory of (Non-Symmetric)
  {D}irichlet Forms.
\newblock Universitext. Springer, Berlin (1992)

\bibitem{Mat}
Mataloni, S.: On a type of convergence for non-symmetric {D}irichlet forms.
\newblock Adv. Math. Sci. Appl. \textbf{9}(2), 749--773 (1999)

\bibitem{Mei20}
Meinert, M.: Partial differential equations on fractals. {E}xistence,
  uniqueness and approximation results.
\newblock Ph.D. thesis, Universit{\"a}t Bielefeld (2020).
\newblock
  \urlprefix\url{https://nbn-resolving.org/urn:nbn:de:0070-pub-29447129}

\bibitem{Me}
Mertens, P.: Konvergenz von koerziven abgeschlossenen {F}ormen.
\newblock Diploma thesis, Universit{\"a}t Bonn (1994)

\bibitem{MSS18}
Mimica, A., Sandri{\'c}, N., Schilling, R.L.: Markov chain approximation of
  pure jump processes.
\newblock Acta Appl. Math. \textbf{158}(1), 167--206 (2018)

\bibitem{M69}
Mosco, U.: Convergence of convex sets and of solutions of variational
  inequalities.
\newblock Adv. Math. \textbf{3}, 510--585 (1969)

\bibitem{M}
Mosco, U.: Composite media and asymptotic {D}irichlet forms.
\newblock J. Funct. Anal. \textbf{123}(2), 368--421 (1994)

\bibitem{OTT}
Ogura, Y., Tomisaki, M., Tsuchiya, M.: Convergence of local type {D}irichlet
  forms to a non-local type one.
\newblock Ann. Inst. H. Poincar{\'e} Probab. Statist. \textbf{38}(4), 507--556
  (2002)

\bibitem{Osh}
Oshima, Y.: Semi-{D}irichlet Forms and {M}arkov Processes, \emph{de Gruyter
  Studies in Mathematics}, vol.~48.
\newblock Walter de Gruyter, Berlin (2013)

\bibitem{Po}
Posilicano, A.: Convergence of distorted {B}rownian motions and singular
  {H}amiltonians.
\newblock Potential Anal. \textbf{5}(3), 241--271 (1996)

\bibitem{PoZ}
Posilicano, A., Zhang, T.S.: Convergence of symmetric diffusions on {W}iener
  spaces.
\newblock Acta Math. Appl. Sin. Engl. Ser. \textbf{20}(1), 19--24 (2004)

\bibitem{ReSi1}
Reed, M., Simon, B.: Methods of Modern Mathematical Physics. {I}. {F}unctional
  Analysis.
\newblock Academic Press, New York (1972)

\bibitem{RS2}
R{\"o}ckner, M., Schmuland, B.: Quasi-regular {D}irichlet forms: examples and
  counterexamples.
\newblock Canad. J. Math. \textbf{47}(1), 165--200 (1995)

\bibitem{RZ3}
R{\"o}ckner, M., Zhang, T.S.: Uniqueness of generalized {S}chr{\"o}dinger
  operators and applications.
\newblock J. Funct. Anal. \textbf{105}(1), 187--231 (1992)

\bibitem{RZ4}
R{\"o}ckner, M., Zhang, T.S.: Uniqueness of generalized {S}chr{\"o}dinger
  operators. {II}.
\newblock J. Funct. Anal. \textbf{119}(2), 455--467 (1994)

\bibitem{RZ}
R{\"o}ckner, M., Zhang, T.S.: Finite dimensional approximation of diffusion
  processes on infinite dimensional spaces.
\newblock Stochastics Stochastics Rep. \textbf{57}(1-2), 37--55 (1996)

\bibitem{RZ2}
R{\"o}ckner, M., Zhang, T.S.: Convergence of operator semigroups generated by
  elliptic operators.
\newblock Osaka J. Math. \textbf{34}(4), 923--932 (1997)

\bibitem{Schw}
Schwartz, L.: Radon Measures on Arbitrary Topological Spaces and Cylindrical
  Measures.
\newblock Oxford University Press, London (1973)

\bibitem{S2}
Stannat, W.: Dirichlet forms and {M}arkov processes: a generalized framework
  including both elliptic and parabolic cases.
\newblock Potential Anal. \textbf{8}(1), 21--60 (1998)

\bibitem{S1}
Stannat, W.: The theory of generalized {D}irichlet forms and its applications
  in analysis and stochastics.
\newblock Mem. Amer. Math. Soc. \textbf{142}(678), viii+101 (1999)

\bibitem{Str}
Stroock, D.W.: Diffusion semigroups corresponding to uniformly elliptic
  divergence form operators.
\newblock In: S{\'e}minaire de {P}robabilit{\'e}s {XXII}, \emph{Lecture Notes
  in Mathematics}, vol. 1321, pp. 316--347. Springer, Berlin (1988)

\bibitem{Sun}
Sun, W.: Mosco convergence of quasi-regular {D}irichlet forms.
\newblock Acta Math. Appl. Sinica (English Ser.) \textbf{15}(3), 225--232
  (1999)

\bibitem{Suz}
Suzuki, K.: Convergence of non-symmetric diffusion processes on {RCD} spaces.
\newblock Calc. Var. Partial Differential Equations \textbf{57}(5), Paper No.
  120 (2018)

\bibitem{Suz22}
Suzuki, K.: Regularity and stability of invariant measures for diffusion
  processes under synthetic lower {R}icci curvature bounds.
\newblock Ann. Sc. Norm. Super. Pisa Cl. Sci. (5) \textbf{23}, 745--808 (2022)

\bibitem{T}
Takeda, M.: On a martingale method for symmetric diffusion processes and its
  applications.
\newblock Osaka J. Math. \textbf{26}(3), 605--623 (1989)

\bibitem{Toe1}
T{\"o}lle, J.M.: Convergence of non-symmetric forms with changing reference
  measures.
\newblock Diploma thesis, Universit{\"a}t Bielefeld (2006).
\newblock
  \urlprefix\url{https://bibos.math.uni-bielefeld.de/preprints/E06-09-234.pdf}.
\newblock BiBoS Preprint E06-09-234

\bibitem{Toe2}
T{\"o}lle, J.M.: Variational convergence of nonlinear partial differential
  operators on varying {B}anach spaces.
\newblock Ph.D. thesis, Universit{\"a}t Bielefeld (2010).
\newblock \urlprefix\url{https://nbn-resolving.org/urn:nbn:de:hbz:361-16758}

\bibitem{Tr}
Trutnau, G.: Stochastic calculus of generalized {D}irichlet forms and
  applications to stochastic differential equations in infinite dimensions.
\newblock Osaka J. Math. \textbf{37}(2), 315--343 (2000)

\bibitem{Tr2}
Trutnau, G.: A short note on {L}yons--{Z}heng decomposition in the
  non-sectorial case.
\newblock RIMS K{\^o}ky{\^u}roku Bessatsu \textbf{B6}, 237--245 (2008)

\bibitem{Uem}
Uemura, T.: On weak convergence of diffusion processes generated by energy
  forms.
\newblock Osaka J. Math. \textbf{32}(4), 861--868 (1995)

\bibitem{Wit23}
Wittmann, S.: Large scale asymptotics for {M}arkov processes in the analytic
  framework of {M}osco--{K}uwae--{S}hioya.
\newblock Ph.D. thesis, Rheinland-Pf{\"a}lzische Technische Universit{\"a}t
  Kaiserslautern-Landau (2023).
\newblock
  \urlprefix\url{https://nbn-resolving.org/urn:nbn:de:hbz:386-kluedo-71900}

\bibitem{Y}
Yosida, K.: Functional Analysis, \emph{Grundlehren der Mathematischen
  Wissenschaften}, vol. 123, 5th edn.
\newblock Springer, Berlin (1978)

\end{thebibliography}
\end{document}